\documentclass{article}
\usepackage{amsmath,amsfonts,amsthm,amssymb,amscd,color,xcolor,mathrsfs,verbatim,microtype}
\usepackage[notref,color]{showkeys}
\usepackage{framed}
\usepackage{graphicx,eurosym}
\usepackage{comment}
\usepackage{amsfonts}
\usepackage{color}
\usepackage{bbm}
\usepackage{bm}
\usepackage{hyperref}
\usepackage{stmaryrd} 

\usepackage[applemac]{inputenc}

\usepackage[cyr]{aeguill}

\colorlet{darkblue}{blue!50!black}

\hypersetup{
    colorlinks,%
    citecolor=blue,%
    filecolor=red,%
    linkcolor=darkblue,%
    urlcolor=blue,%
    pdfnewwindow=true,%
    pdfstartview={FitH}
}

\usepackage{graphicx,amscd,mathrsfs,wrapfig,mathrsfs,lipsum}
\usepackage{eufrak}
\usepackage{float}
\usepackage{tikz}
\usepackage{multicol}
\usepackage{caption}
\usetikzlibrary{arrows}
\usepackage{capt-of}

\colorlet{darkblue}{blue!50!black}

\usepackage{marginnote}

\newcommand{\p}{\partial}
\newcommand{\e}{\varepsilon}

\newcommand{\R}{{\mathbb R}}
\newcommand{\Z}{{\mathbb Z}}

\newcommand{\I}{{\mathbb I}}
\newcommand{\E}{{\mathbb E}}

\newcommand{\ty}{\infty}

\newcommand{\mW}{{\mathbb W}}
\newcommand{\mH}{{\mathbb  H}}

\newcommand{\mT}{{\mathbb T}}

\newcommand{\aA}{{\cal A}}

\newcommand{\KK}{{\cal K}}

\def\mE{{\mathbb E}}
\def\dif{{\mathord{{\rm d}}}}
\def\cF{{\mathcal F}}
\def\cC{{\mathcal C}}
\def\cZ{{\mathcal Z}}
\def\cK{{\mathcal K}}
\def\mP{{\mathbb P}}
\def\mN{{\mathbb N}}
\def\mS{{\mathbb S}}
\def\mR{{\mathbb R}}
\def\mZ{{\mathbb Z}}
\def\cB{{\mathcal B}}
\def\cD{{\mathcal D}}
\def\cM{{\mathcal M}}
\def\cA{{\mathcal A}}
\def\cS{{\mathcal S}}

\def\cR{{\mathcal R}}
\def\cH{{\mathcal H}}
\def\eps{\varepsilon}
\def\cL{{\mathcal L}}
\def\mI{{\mathbb I}}

\newcommand{\lag}{\langle}
\newcommand{\rag}{\rangle}

\newcommand{\dd}{{\textup d}}

\def\[{{\big[}}
\def\]{{\big]}}

\theoremstyle{plain}

\newtheorem*{lemma*}{Lemma}
\newtheorem{theorem}{Theorem}[section]
\newtheorem{lemma}[theorem]{Lemma}
\newtheorem{proposition}[theorem]{Proposition}

\theoremstyle{definition}

\newtheorem{condition}[theorem]{Condition}

\theoremstyle{remark}
 \newtheorem{remark}[theorem]{Remark}

\allowdisplaybreaks

\numberwithin{equation}{section}

\begin{document}

\title{\Large\bf Ergodicity   for  2D Navier-Stokes equations with
a degenerate pure jump noise
  }

\author{
   {Xuhui Peng}$^a$\footnote{E-mail: xhpeng@hunnu.edu.cn}\;
\ \ \ {Jianliang Zhai}$^b$\footnote{Email: zhaijl@ustc.edu.cn}\;
\ \ \
{Tusheng Zhang}$^{b,c}$\footnote{E-mail: tushengz@ustc.edu.cn
}
\\
 \small a. MOE-LCSM, School of Mathematics and Statistics, Hunan Normal University,\\
  \small  Changsha, Hunan 410081, P. R.China\\
 \small b. School of Mathematical Sciences, University of Science and Technology of China,
 \\ \small  Hefei,  Anhui 230026, P.R.China   \\
\small c. Department of Mathematics, University of Manchester,
\\  \small Oxford Road, Manchester, M13
9PL, UK.
 }\,
  \date{\today}
\maketitle

%
%

\begin{abstract}

In this paper, we establish the ergodicity  for stochastic  2D Navier-Stokes equations driven by
a highly degenerate  pure jump L\'evy noise. The noise could appear in as few as four directions.  This gives an affirmative anwser to a longstanding problem. The case of Gaussian noise was treated in Hairer and  Mattingly [\emph{Ann. of Math.}, 164(3):993--1032, 2006]. To obtain the uniqueness of  invariant  measure, we use
Malliavin calculus and anticipating stochastic calculus to establish the equi-continuity of the semigroup, the so-called {\em e-property}, and prove some weak irreducibility of the solution process.

\vspace{0.5cm}

\smallskip


\smallskip


\end{abstract}

%
%
%

\section{Introduction and Main results}
\subsection{Introduction}

In the theory of
turbulence, the study
of the equations of fluid mechanics driven by degenerate  noise forcing, that is, the driving noise does not act directly on all  the determining modes of the flow,  is ubiquitous; see e.g., \cite{Eyi96,Nov65,Sta88,VKF79}.
And in the physics literature when discussing
the behavior of stochastic fluid dynamics  in the turbulent regime, the main assumptions usually made are ergodicity and statistical
translational invariance of the stationary state.  The uniqueness of an invariant  measure and the ergodicity of the randomly forced dissipative partial
 differential equations(PDEs) driven by degenerate noise have been the problem of central concern for many years.

Because of the complexity and the difficulty of the problem, it is much less understood and there are only a few works on this topic. In this paper, we confine ourselves to stochastic 2D Navier-Stokes equations. In \cite{HM-2006,HM-2011} the authors studied the 2D Navier-Stokes equations on the torus and the sphere and established the exponential mixing, provided that the random perturbation is white in time and contains several Fourier modes. In \cite{shirikyan-asens2015}, the exponential mixing was established for the 2D Navier-Stokes system perturbed by a space-time localised smooth stochastic forcing. In the paper \cite{KNS-2018} the authors proved a similar result in the situation when random forces are localised in the Fourier space and coloured in time. The problem of mixing for the Navier-Stokes system with a random perturbation acting through the boundary has been studied in \cite{Shi2021}.
The authors in  \cite{KS-book} proved the polynomial mixing of the 2D Navier-Stokes equations driven by a compound Poisson process.
 We remark that the volume of the  intensity measure of compound Poisson process  is finite.

So far, there are no results for the case when the random perturbation is pure jump L\'evy noise of infinite activity, that is, the volume of the  intensity measure of the L\'evy noise is infinite. This is the subject of the present article. 
We point out that there are no results even for the non-degenerate case of pure jump L\'evy noise of infinite activity, i.e., all determining modes of the unforced PDE are directly affected by the noise.

\vskip 0.2cm

Now, let us give a brief introduction to the main result.
The Navier-Stokes (NS) equations on the torus  $\mT^2=[-\pi,\pi]^2$ are given by
 \begin{equation}
 \label{16-4}
\partial_t u + (u\cdot \nabla)u = \nu \Delta  u -  \nabla p + \xi,
\text{ div }u = 0,
 \end{equation}
 where
$u=u(x, t)\in \mR^2$ denotes the value of the velocity field at time
 $t$ and position
$x=(x_1,x_2)$, $p(x, t)$ denotes the pressure, $\nu>0$ is the viscosity  and $\xi=\xi(x, t)$
 is an external force field acting on
the fluid.
Since the velocity and
vorticity formulations are equivalent in this setting, we choose to use the vorticity equation as this simplifies the exposition.
For a divergence-free velocity
field $u$, we define the vorticity $w$ by $w=\nabla\wedge u=\partial_2u_1-\partial_1u_2$. Note that $u$ can be
recovered from $w$ and the condition $\text{ div }u = 0$.  With these notations
we rewrite the NS system \eqref{16-4}  in the vorticity formulation:
   \begin{equation}
  \label{0.1}
    \p_t w =\nu \Delta w+B(\cK w,w)+ \partial_t \eta, \quad  w\big|_{t=0}=w_0,
  \end{equation}
  where $\eta=\eta(x,t)$ is   an external random force,
    $B(u,w)=-( u \cdot \nabla ) w$,
  and  $\KK$ is the Biot--Savart operator
which
 is defined in Fourier space
by $(\cK w)_k =-iw_kk^{\bot}/|k|^2$,  where $k^\bot=(k_1,k_2)^\bot=(-k_2,k_1)$, and  $w_k$ is
the scalar product of $w$ with $\frac{1}{2\pi} e^{i k \cdot x}$.
The Biot--Savart operator has the property that the
divergence of $\cK w$ vanishes and that $w=\nabla\wedge(\cK w)$.

In this paper, we prove that there exists a unique statistically
invariant state of the system (\ref{0.1}).
Roughly  speaking, we establish   the following result. For  rigorous statement { and general version}
of the result, please see  Theorem \ref{16-8} below.

\vspace{0.2cm}
{  Consider the system}  (\ref{0.1}) with noise $$\eta(x,t)=
  b_1\sin(x_1)L_1(t)+b_2\cos(x_1)L_2(t)+b_3\sin(x_1+x_2)L_3(t)+b_4\cos(x_1+x_2)L_4(t),$$ where
  $b_1,\cdots,b_4$ are non-zero constants,
  $L_t=(L_1(t),\cdots,L_4(t))$
  is a $4$-dimensional  pure jump  process with L\'evy measure $\nu_L$:
  \begin{eqnarray}
  \label{23-1}
  \nu_L(\dif z)=\int_0^\infty
  (2\pi u)^{-2} e^{-\frac{|z|^2}{2u}} \nu_{S}(\dif u)\dif z,
  \end{eqnarray}
  here   $\nu_S$ is a measure on $(0,\infty)$ satisfying
  \begin{eqnarray*}
   \nu_S((0,\infty))=\infty \text{ and } \int_0^\infty (e^{\zeta   u}-1)\nu_S(\dif u)<\infty   \text{  for some }
 \zeta  >0
.
\end{eqnarray*}
Then the main results imply that the  Markov semigroup generated by the system  (\ref{0.1}) possesses a unique  invariant measure $\mu^*$ on the space  $H=\Big\{w\in L^2(\mT^2,\R):  \int_{\mT^2} w(x)\dif x=0 \Big\}$.
\vskip 0.2cm

There are now empirical data which shows that L\'evy processes
are more suitable to realistically represent external forces in statistical physics(c.f. \cite{Nov65}), climatology(c.f. \cite{IP06}) and mathematics of finance(c.f. \cite{KT13}).
Therefore,  the mathematical analysis of stochastic partial differential equations  driven by L\'evy processes becomes very important.
This motivates the study of this paper.

To prove the main result,  we will use
Malliavin calculus and anticipating stochastic calculus to establish an equi-continuity of the semigroup, the so-called {\em e-property}, and prove some weak irreducibility for the solution process of the system  (\ref{0.1}).

To deal with the setting of highly degenerate noises, we need some quantitative control of the Malliavin matrix, and it is inevitable to use the ``future information". Hence, some anticipating
stochastic analysis are necessary. However, Malliavin calculus associated with Poisson random measures is much less effective than that of the Wiener case.
In this paper, we assume that the driving noise is a subordinated Brownian motion. Introducing a sort of time change, we borrow the nonadapted
stochastic analysis associated with Brownian motion when dealing with the ``future information".

Because of  the strong intensity  of  the jumps and also the unbounded jumps, we need to introduce a new set of ideas and techniques to establish the uniqueness of invariant
measures in comparison to the case of Gaussian noise (see \cite{HM-2006}  \cite{HM-2011}\cite{FGRT-2015}).  Now we
highlight some of them.

\begin{itemize}
    \item In \cite{HM-2006}, the authors gave preliminary
estimates for the solutions in \cite[Lemma 4.10]{HM-2006}, which plays an essential role in controlling various terms during the proof of the asymptotically strong Feller.
The proof of \cite[Lemma 4.10]{HM-2006} strongly relies on the Girsanov transformation and exponential martingale estimate of the Gaussian noise. However,
 the Girsanov transformation and exponential martingale inequality associated with Poisson random measures
is expressed in terms of complicated exponential type nonlinear transformations, which seems very hard (if not impossible) to control; see \cite[Theorem 5.2.9]{Dav-2009}.
In the  setting of pure jump L\'evy noise, to overcome the above difficulty,
we design a sequence of stopping times $\sigma_k$ (see (\ref{2121-4}) and (\ref{qu-20})), and build new preliminary estimates, see Lemma \ref{qu-2} in this paper, which is totally different from the ones for the Gaussian case, i.e.,
\cite[Lemma 4.10]{HM-2006}. We point out that Lemma \ref{qu-2} seems not possible to be proved by the exponential martingale inequalities and the $\rm It\hat{o}$ formula for Poisson random measures.

With the help of the stopping times $\sigma_k$, we could  identify the  ``bad part" of the sample space $\Omega$, denoted by $\{\omega\in\Omega:\Theta>M\}$; see (\ref{58-1}) and (\ref{1609}) for the definition of the random variable
$\Theta$; and the ``bad part" could be controlled by the strong law of large numbers, which means that $\lim_{M\rightarrow\infty}\mathbb{P}(\{\omega\in\Omega:\Theta>M\})=0$.
On the ``good part" of the sample space, we could obtain something like the asymptotic strong Feller property; see (\ref{180-3}). Combining the two parts, we obtain the $e$-property.

    \item Let $\cM_{0,t}$ be the Malliavin matrix of $w_t$
    and
   $\cS_{\alpha,N}$ be some   subspace of $H$
   (For the definition of $\cM_{0,t}$,   see Section 4.2 in \cite{HM-2006}  or   (\ref{200-3})--(\ref{200-1}) below).
To obtain the ergodicity via Malliavin calculus,
one key  ingredient    is   to show
   \begin{align}
   \label{18-2}
     \mP(\inf_{\phi \in\cS_{\alpha,N} } \langle \cM_{0,1}\phi,\phi\rangle <\eps)\leq C(\|w_0\|) r(\eps), \forall \eps\in (0,1)\text{ and } w_0\in H,
   \end{align}
   where  $\|w_0 \|$ denotes the $L^2$ norm of $w_0$,   $C$ is some function  from $[0,\infty)$ to $ [0,\infty)$  and  $r$ is a function on $(0,1)$ with $\lim_{\eps\rightarrow 0}r(\eps)=0$.
      In the existing literatures   \cite{HM-2006}\cite{HM-2011}\cite{FGRT-2015} etc.,
   the properties of  Gaussian polynomials(see, e.g.,   \cite[Theorem 7.1]{HM-2011}) play very essential roles for the estimate of the left side  of (\ref{18-2}).
  Similar arguments do not work in the case of pure jump processes. In this paper, using the fact that the jump times of  pure jump noise with infinite activity are dense in any time interval $[a,b]$ with $0\leq a<b$,
  we find a new way to get something like (\ref{18-2}).
  First, we prove
  \begin{eqnarray}
  \label{18-4}
  \mP(\inf_{\phi \in\cS_{\alpha,N} } \langle \cM_{0,\sigma_1}\phi,\phi\rangle =0)=0, ~\forall w_0\in H,
  \end{eqnarray}
  where $\sigma_1$ is a positive stopping time.
  Then, with the help of (\ref{18-4}) and the dissipative property of Navier-Stokes system,
 we derive a  weaker version   of   (\ref{18-2}) which   is  sufficient for our purpose.
 In a word, our method of proving something like  (\ref{18-2}) is totally different from that of the Gaussian case.
  See  Proposition \ref{1-66} and Proposition \ref{3-8} in  Section \ref{S:3} for more  details.
\end{itemize}


\vskip 0.2cm
Finally, we point out that
there are not many results on the ergodicity of stochastic partial
differential equations driven by pure jump L\'evy noise. And we list them here for readers who are interested. For the case that the driving noise is non-degenerate, we refer to \cite{PZ2011,PXZ2011,PSXZ2012,Xu13,WXX17,WX18,DXZ-2014,DWX20,BHR16SIAM,FHR16CMP,WYZZ 2207}.
For the case that the driving noise is degenerate, we refer to \cite{SXX19,WYZZ 2209,MR10}.

\subsection{Main results}
\label{S-1-2}

  We consider the system  (\ref{0.1}) in the Hilbert space:
\begin{equation}\label{0.2}
H=\Big\{w\in L^2(\mT^2,\R):  \int_{\mT^2} w(x)\dif x=0 \Big\}
\end{equation}
   endowed    with the   $L^2$-scalar product $\lag \cdot, \cdot\rag$ and the corresponding $L^2$-norm $\|\cdot\|.$

 In order to describe the noise $\eta$,   we introduce the following notation.
 Denote
 $$\mZ_{+}^2=\{(k_1,k_2)\in \mZ^2: k_2>0 \} \cup \{(k_1,0)\in \mZ^2: k_1>0 \}
 $$  and
 $\mZ_{-}^2=\{(k_1,k_2):-(k_1,k_2)\in \mZ_{+}^2\}.$
 For any $k=(k_1,k_2)\in \mZ_*^2:=\mZ^2\backslash (0,0),$ set
  $$
{ e_k}=e_k(x)=
\begin{cases}
 \sin \lag k,x\rag& \text{if } (k_1,k_2)\in \mZ_{+}^2, \\
  \cos\lag k,x\rag   & \text{if } (k_1,k_2)\in \mZ_{-}^2.
 \end{cases}
$$
We   assume that $\eta$ is a white-in-time~noise of the form
\begin{equation}\label{0.4}
\eta(x,t)=  \sum_{k\in \cZ_0}
b_kW_{S_t}^k  e_k(x),
\end{equation}
where $\cZ_0 \subset \Z^2_*$ is a  finite set,
 $ b_k, {k\in \cZ_0} $ are  non-zero constants,
 $W_{S_t}=(W_{S_t}^k)_{k\in \cZ_0}$ is a
$|\cZ_0|$-dimensional subordinated  Brownian motion which will be specified below.
For convenience,  we  always denote $|\cZ_0|$ by $d.$
Assume  the canonical basis of $\mR^d$
is
$\{\theta_k\}_{k\in \cZ_0}$
and the   linear operator {$Q: \mR^d\rightarrow H$ is defined   in the following way: 
  \begin{eqnarray*}
    Qz=\sum_{k\in \cZ_0 }   b_{k} z_k e_{k}, ~\forall z=\sum_{k\in \cZ_0}z_k\theta_k \in \mR^d,
  \end{eqnarray*}
  then, $\eta(t) =QW_{S_t}.$

Now let us give the details for the subordinated Brownian motion
$W_{S_t}$.  Let
$(\mW, \mH, \mP^{\mu_\mW})$
 be the classical Wiener space, i.e.,
 $\mW$ is the space of all continuous
functions from $\mR^{+}$ to $\mR^d$
 with vanishing values at starting point $0$, $\mH \subseteq \mW$
  is the
Cameron-Martin space consisting of all absolutely continuous
functions with square
integrable derivatives,
$\mP^{\mu_\mW}$ is the Wiener measure so that the coordinate process
$W_t(\mathrm{w}  ) :=\mathrm{w}_t$
is a $d$-dimensional standard Brownian motion.
 Let $\mS$ be the space of all $\rm c\grave{a}dl\grave{a}g$
 increasing functions $\ell$ from $\mR^{+}$
 to $\mR^{+}$  with $\ell_0 = 0.$
 Suppose that $\mS$ is endowed with the
Skorohod metric and a probability measure
$\mP^{\mu_\mS}$ so that the coordinate process
$S_t(\ell) := \ell_t$
is a pure jump subordinator
 with L\'evy measure $\nu_{S}$ satisfying
$$
\int_0^\infty (1 \wedge  u)\nu_S(du) <\infty.
$$
Consider now the
following product probability space
$(\Omega, \cF, \mP) :=
 (\mW \times  \mS, \cB(\mW) \times \cB(\mS), \mP^{\mu_\mW} \times \mP^{\mu_\mS}
)$,
and define for $\omega = (\mathrm{w}, \ell) \in  \mW \times  \mS$,
$L_t(\omega):= \mathrm{w}_{\ell_t}.$
Then, $(L_t=W_{S_t})_{t\geq 0}$ is a $d$-dimensional pure jump L\'evy process with
L\'evy measure $\nu_L$ given by
\begin{eqnarray}
\label{1129}
  \nu_L(E)=\int_0^\infty \int_E
  (2\pi u)^{-\frac{d}{2}} e^{-\frac{|z|^2}{2u}}\dif z \nu_{S}(\dif u),
  \quad E\in \cB(\mR^d).
\end{eqnarray}

To formulate the main result, let us recall that a set $\cZ_0 \subset \Z^2_*$
is a generator if any element of~$\Z_*^2$
is a finite linear
combination of elements of~$\cZ_0$ with integer coefficients.
In~what follows, we    assume that the following  two  conditions ~are in place.
\begin{condition}
\label{16-5}
 The set  $\cZ_0 \subset \Z^2_*$ appeared in \eqref{0.4} is a
 finite, symmetric (i.e.,  $-\cZ_0=\cZ_0$)  generator
that contains at least two non-parallel vectors~$m$
and~$n$
such that $|m|\neq |n|$.
\end{condition}
This is   the   condition under which   the  ergodicity
of the NS system is  established  in~\cite{HM-2006, HM-2011}
in the case of a white-in-time noise and in \cite{KNS-2018} in
 the case of a bounded   noise. The set
$$
\cZ_0=\{(1,0), \, (-1,0), \, (1,1), \, (-1,-1)\}\subset \Z^2_*:=\mZ^2 \backslash \{(0,0)\}
$$ is an   example satisfying this condition.

\begin{condition}
\label{14-2}
Assume that  $\nu_S$ satisfies
\begin{eqnarray*}
  \int_0^\infty (e^{\zeta   u}-1)\nu_S(\dif u)<\infty   \text{  for some }
 \zeta  >0
\end{eqnarray*}
and

\begin{eqnarray}\label{Condition nu S}
  \nu_S((0,\infty))=\infty.
\end{eqnarray}
\end{condition}

%

\vskip 0.4cm

\begin{remark}
  If  $\nu_S(\dif u)=u^{-1-\frac{\alpha}{2}}
I_{\{0< u\leq \aleph  \}}\dif u $
for some {$\alpha\in [0,2) $} and  $ \aleph >0,$  then
 condition \ref{14-2} is satisfied.
In this case, $\nu_{L}(\dif z)=\zeta(z)\dif z$ and   $\zeta(z)$ satisfies
\begin{eqnarray}
\label{1121}
\begin{split}
   & \frac{C_1}{|z|^{d+\alpha}}\leq  \zeta(z) \leq \frac{C_2}{|z|^{d+\alpha}},\quad \forall |z|\leq 1,
  \\ & \zeta(z) \leq \frac{C_3}{|z|^{d+\alpha}}\exp\{-\frac{|z|^2}{4\aleph }\},\quad \forall |z|\geq 1,
  \end{split}
\end{eqnarray}
where $C_1,C_2,C_3$ are positive constants depending on $\alpha,d$ and $\aleph $.
Thus, the appearance  of the  small jumps of $L_t$ will behave like $\alpha$-stable processes  and  the  appearance  of big  jumps will be
very rare.
\end{remark}

We denote by $P_t(w_0,\cdot)$ the transition probabilities of the solution of the stochastic Navier-Stokes equation (\ref{0.1}), i.e,
\begin{eqnarray*}
  P_t(w_0,A)=\mP(w(t)\in A\big|w(0)=w_0)
\end{eqnarray*}
for every Borel set $A\subseteq H$ and
\begin{equation*}
P_t f(w_0)=\int_{H}  f(w)P_t(w_0,\dif w), \quad  P_t^*\mu(A)=\int_{H}P_t(w_0,A)\mu(\dif w_0),
\end{equation*}
for every $f:H\rightarrow \mathbb{R}$ and probability measure $\mu$ on $H.$

Before we state the main theorem in this paper, we present two propositions which are the essential ingredients in the proof of the main result.

\begin{proposition}
\label{3-11}
Under  the   Condition \ref{16-5} and Condition \ref{14-2},
 the Markov semigroup $\{P_t\}_{t\geq 0}$
 has the e-property,  i.e.,
  for any bounded and Lipschitz continuous function $f$, $w_0\in H$ and $\eps>0,$ there exists a $\delta>0$ such that
 \begin{eqnarray*}
 |P_tf(w_0')-P_tf(w_0)|<\eps, ~\forall t\geq 0 \text{ and } w_0' \text{ with }
 \|w_0'-w_0\|<\delta.
 \end{eqnarray*}
\end{proposition}


The proof of this proposition is given in Section \ref{54-1} based on Malliavin calculus, which constitutes a major part of the paper.
\vskip 0.3cm
Since there are many constants appearing in the proof, we adopt  the following convention.
Without otherwise specified,  the letters $C,C_1,C_2,\cdots$ are  always  used to denote unessential constants that  may
change from line to line  and  implicitly  depend on the data of the system (\ref{0.1}), i.e.,
 $\nu,\{b_k\}_{k\in \cZ_0},\nu_S$ and $d=|\cZ_0|.$
 Also, we usually do not explicitly indicate the dependencies on the parameters  $\nu,\{b_k\}_{k\in \cZ_0},\nu_S$ and $d=|\cZ_0|$ on every occasion.
The proof of the proposition below is almost the same as that in
 \cite[Lemma 3.1]{EM01}; for the convenience of the readers,  we give its  short proof in Section \ref{54-2}.
\begin{proposition}(Weak Irreducibility)
\label{16-6}
  For any $\mathcal{C},\gamma>0$, there exists a $T=T(\mathcal{C},\gamma)>0$ such that
  \begin{eqnarray*}
    \inf_{\|w_0\|\leq \mathcal{C}}P_{T}(w_0, \cB_\gamma)>0,
  \end{eqnarray*}
  where $\cB_\gamma=\{w\in H:~\|w\|\leq \gamma\}.$
\end{proposition}

After we state Proposition   \ref{3-11} and Proposition  \ref{16-6}, we have the following main result of the paper.

\begin{theorem}
\label{16-8}
 Consider the  2D Navier-Stokes equation  (\ref{0.1}) with a
   degenerate   pure jump noise
  (\ref{0.4}).
  Under the Condition \ref{16-5} and Condition {\ref{14-2}},
   there exists a unique invariant measure $\mu^*$
for the system (\ref{0.1}), i.e., $\mu^*$ is a unique   probability measure on $H$
 such that $P_t^*\mu^*=\mu^*$ for every $t\geq 0.$
\end{theorem}
\begin{proof}
 We first prove the existence.
 By     Lemma \ref{166-3} below, it holds that
 \begin{eqnarray*}
   \frac{\nu}{2}~ \mE \int_0^t\|w_s\|_1^2\dif s\leq \|w_0\|^2+C t,
 \end{eqnarray*}
 here the definition of $\|\cdot\|_1$ is introduced in (\ref{eq 2024 04 15 00}). Following the arguments in the proof of  \cite[Theorem 5.1]{DX-2009} and using Krylov-Bogoliubov criteria,
 we obtain the  existence of invariant measure.

 Now we prove the  uniqueness.
 Assume that there are two distinct invariant probability measures ${ \mu_1}$ and ${\ \mu_2}$ for
 $\{P_t\}_{t\geq0}$.
  By   Proposition \ref{3-11}  and \cite[Theorem 1]{KSS12}
  { (or  \cite[Proposition 1.10]{GL15})}, one has
  \begin{eqnarray}
  \label{16-10}
    \text{Supp } { \mu_1} \cap  \text{Supp } { \mu_2}=\emptyset.
  \end{eqnarray}
  On the other hand,   by (\ref{A-2}) in   Lemma \ref{166-3} below, for every
  invariant measure
  $\mu,$
  the following priori bound
  \begin{eqnarray*}
    \int_H\|w\|^2\mu(\dif w)\leq C
  \end{eqnarray*}
holds(See \cite[Lemma B.2]{EMS-2001}).
Following the arguments in the proof of \cite[Corollary 4.2]{HM-2006} and using Proposition \ref{16-6} ,
for every invariant measure
  $\mu,$ we have $0\in \text{Supp } \mu.$
  This contradicts (\ref{16-10}).
  We complete the proof of uniqueness.
\end{proof}
\vskip 0.5cm
    The rest of the paper is organised as follows.
    In Section \ref{S:2}, we provide some estimates for the solution $w_t$ and introduce the essential ingredients of the Malliavin  calculus for the solution. Moreover, we give all the necessary estimates associated with the Malliavin matrix.
    Section \ref{S:3} is devoted to  the proof of the invertibility of the  Malliavin matrix  of  the solution  $w_t$ which plays a key role in the proof of Proposition \ref{3-11}.
    In Section \ref{54-1}, we give the proof  of  Proposition \ref{3-11}.
The proof  of  Proposition \ref{16-6} is given In Section \ref{54-2}. Some of the technical proofs are put in the Appendix.

\section{Preliminaries}
\label{S:2}

\subsection{Notations}

In this paper, we use the following  notations.
 \smallskip
\noindent
Let ${H}_N=span\{e_j: j\in \mZ_*^2  \text{ and } |j|\leq N\}.$
${P}_N$  denotes the orthogonal projections from $H$ onto ${H}_N$.
      Define   ${Q}_Nu:=u-{P}_Nu, \forall u\in H.$

\smallskip
\noindent
For $\alpha\in \mR$ and a smooth function $w\in H$, we define the norm $\|w\|_\alpha$ by
\begin{eqnarray}\label{eq 2024 04 15 00}
\|w\|_\alpha^2=\sum_{k\in \mZ_*^2}|k|^{2\alpha}w_k^2,
\end{eqnarray}
where $w_k$ denotes the Fourier mode with wavenumber $k.$
 When $\alpha=0,$ as stated in Section \ref{S-1-2},  we also    denote  this norm $\|\cdot \|_\alpha$ by $\|\cdot \|.$
For any $(s_1,s_2,s_3)\in \mR_{+}^3$ with $\sum_{i=1}^3s_i\geq 1$
and $(s_1,s_2,s_3)\neq (1,0,0),(0,1,0),(0,0,1)$, the following relations will be used frequently in this paper(c.f. \cite{CF88}):
\begin{eqnarray}
\label{44-1}
\begin{split}
\langle B(u,v),w\rangle& = -\langle B(u,w),v\rangle, \quad \text{if }\nabla \cdot u=0,
\\
  \big|  \langle B(u,v),w\rangle\big| &\leq
C\|u\|_{s_1}\|v\|_{1+s_2}\|w\|_{s_3},
\\ \|\cK u\|_{\alpha} & = \|u\|_{\alpha-1},
\\  \|w\|_{1/2}^2 &\leq  \|w\|\|w\|_1.
\end{split}
\end{eqnarray}



\smallskip
\noindent
$L^\ty(H)$  is the space of bounded Borel-measurable functions $\psi:H\to\R$   with the norm $\|\psi\|_\ty=\sup_{u\in H}|\psi(u)|$.~$C_b(H)$ is the space of continuous functions. ~$C^1_b(H)$ is the space of functions~$\psi\in C_b(H)$ that are continuously Fr\'echet differentiable with bounded derivatives.
\noindent
$\cL(X,Y)$ is the space of bounded linear operators from Banach spaces $X$ into Banach space $Y$ endowed with the natural norm $\|\cdot\|_{\cL(X,Y)}$.
If there are no confusions, we always write  the operator  norm  $\|\cdot\|_{\cL(X,Y)}$
as $\|\cdot\|$.

\vskip 0.4cm
\noindent
Let $N_{L}(\dif t,\dif z)$ be the Poisson random measure associated with the L\'evy process $L_t=W_{S_t}$, i.e.,
\begin{eqnarray*}
N_{L}((0,t]\times U)&=& \sum_{s\leq t}I_{U}(L_t-L_{t-}), U\in \cB(\mR^d\setminus\{0\}).
\end{eqnarray*}
Let $\tilde N_{L}(\dif t,\dif z)$ denote the compensated Poisson random measure associated with $N_{L}(\dif t,\dif z)$, i.e.,
\begin{eqnarray*}
\tilde N_{L}(\dif t,\dif z)=N_{L}((0,t]\times U)-\dif t \nu_{L}(\dif z).
\end{eqnarray*}
Similar notation also apply to  $N_{S}(\dif t,\dif z)$ and
$\tilde  N_{S}(\dif t,\dif z).$ As the measure $\nu_{L}(\dif z)$ is symmetric, the L\'evy process $L_t$ admits the following representation:
$$L_t=\int_0^t\int_{\mathbb{R}^d\setminus\{0\}}z\tilde N_L (\dif s, \dif z).$$


\smallskip
\noindent
Let $F=F(\mathrm{w},\ell)$ be a random variable on the space $(\Omega,\cF,\mP).$
We use $\mE^{\mu_{\mW}} F$ to denote the expectation of
$F$   when we take  the element $\ell$ as fixed,
i.e,
\begin{align*}
  \mE^{\mu_{\mW}} F =\int_{\mW}F(\mathrm{w},\ell) \mP^{ \mu_{\mW}}(\dif\mathrm{w} ).
\end{align*}
The notation   $\mE^{\mu_{\mS}} F$ has the similar meaning.
We use  $\mE F$ to denote the expectation of
$F$   under the measure $\mP=\mP^{\mu_\mW}  \times \mP^{\mu_\mS}.$

\smallskip
\noindent
The filtration  used in this paper is
\begin{eqnarray*}
  \cF_{t}:=\sigma(W_{S_s},S_s:s\leq t).
\end{eqnarray*}
For  any  fixed $\ell\in \mS$  and positive  number $a=a(\ell)$ which is independent of the Brownian motion $(W_t)_{t\geq 0}$,  the filtration $\cF_{a}^{W}$ is defined by
\begin{eqnarray*}
  \cF_{a}^{W}:=\sigma(W_{s} : s\leq a).
\end{eqnarray*}
If $\tau:\Omega\rightarrow [0,\infty]$ is a stopping time with respect to   the filtration $\cF_t$,
$\cF_\tau$ denotes the past $\sigma$-field defined by
\begin{eqnarray*}
  \cF_{\tau}=\{A\in \cF: \forall t\geq 0,A\cap \{\tau\leq t\}\in \cF_t\}.
\end{eqnarray*}

\subsection{Priori estimates  on the solutions}


\begin{lemma}
\label{166-3}
  Let $w_t$ be the  solution to equation (\ref{0.1}) with initial value $w_0$. Then, there exist  positive    constants   $C_1,C_2$,  which depend on  the parameters $\nu,\{b_j\}_{j\in \cZ_0},\nu_S,d$, such that
  \begin{eqnarray}
  \label{A-2}
  \mE \|w_t\|^2 & \leq &  e^{-\nu t}\|w_0\|^2 +C_1, ~\forall t\geq 0,
  \\ \label{A-3}
  \frac{\nu}{2}~\mE \int_0^t \|w_s\|_1^2 \dif s & \leq&  \|w_0\|^2+C_2t,~\forall t\geq 0.
  \end{eqnarray}
\end{lemma}

\begin{proof}

Let $\nu_{L}$ be the intensity measure defined in (\ref{1129}). We claim that
\begin{eqnarray}
\label{0722-1}
\int_{|z|\leq 1}|z|^2 \nu_{L}(\dif z)+\int_{|z|\geq 1}|z|^n \nu_{L}(\dif z)<\infty,\quad \forall n\geq 2.
\end{eqnarray}
We only prove that  $\int_{|z|\geq 1}|z|^n \nu_{L}(\dif z)<\infty, \forall n\geq 2,$ the first term in (\ref{0722-1}) can be treated similarly.
 By definition, we have
\begin{eqnarray*}
  \int_{|z|\geq  1}|z|^n \nu_{L}(\dif z)
  &=&  \int_{|z|\geq 1}|z|^n \[\int_{0}^{\infty}(2\pi u)^{-d/2} e^{-\frac{|z|^2}{2u}}\nu_S(\dif u) \] \dif z
  \\ &=& C_d \int_{1}^\infty  r^{n+d-1}   \[\int_{0}^{\infty }u^{-d/2} e^{-\frac{r^2}{2u}}\nu_S(\dif u) \] \dif r.
    \\ &= & C_d \int_{0}^\infty \nu_S(\dif u)\int_1^\infty  r^{n+d-1}e^{-\frac{r^2}{2u}} u^{-d/2}  \dif r
    \\ &= &C_{d}\int_{0}^\infty \nu_S(\dif u) \int_{1/(2u)}^\infty (2ux)^{\frac{n+d-1}{2}}e^{-x}u^{-d/2}  \frac{2u}{2\sqrt{2ux}}\dif x
    \\ &\leq & C_{d,n}\int_{0}^\infty u^{n/2} \nu_S(\dif u)<\infty.
\end{eqnarray*}

Now, we prove (\ref{A-2}) and (\ref{A-3}).
  Applying It\^o's formula to $\|w_t\|^2,$ we obtain
   \begin{eqnarray}\label{eq 20231011}
   \begin{split}
    \dif \|w_t\|^2& = -2\nu \|w_t\|_1^2\dif t+2\int_{z\in\mathbb{R}^d\setminus\{0\}}
    \langle w_t,Qz\rangle \tilde N_L (\dif t, \dif z)
    \\ & \quad +\int_{z\in\mathbb{R}^d\setminus\{0\}}\|Qz\|^2
    N_L (\dif t, \dif z).
    \end{split}
  \end{eqnarray}
  Set  $C=\int_{z\in\mathbb{R}^d\setminus\{0\}}\|Qz\|^2
    \nu_{L}(\dif z)$, which is a constant  only  depending  on  $\{b_j\}_{j\in \cZ_0},\nu_S,d$. By (\ref{eq 20231011}) and standard arguments,
  \begin{eqnarray*}
   &&  \mE \|w_t\|^2+2\nu  \int_0^t \mE\|w_s\|_1^2\dif s=\|w_0\|^2 + Ct,
  \end{eqnarray*}
  which yields the desired results  (\ref{A-2}) and (\ref{A-3}).
\end{proof}

Let $\sigma_0=0$ and $\mathfrak{B}_0=\sum_{j\in \cZ_0}b_j^2.$
 For any $\kappa>0$, $k\in \mN$ and $\ell\in\mathbb{S}$, we define
   \begin{eqnarray}
  \label{2121-4}
  \sigma=\sigma(\ell)=\sigma_1(\ell)= \inf\big\{t\geq 0:\nu t-8 \mathfrak{B}_0 \kappa \ell_t >1  \big\}
\end{eqnarray}
and
 \begin{eqnarray}\label{qu-20}
   \sigma_k= \sigma_k(\ell)=\inf\{t\geq \sigma_{k-1},~ \nu (t-\sigma_{k-1})-8\mathfrak{B}_0 \kappa (\ell_t-\ell_{\sigma_{k-1}})> 1\}.
 \end{eqnarray}

For the  solutions to equation  (\ref{0.1}) and   theses  stopping times $\sigma_k$(with respect to  $\cF_t$),
we have the following moment estimates.
\begin{lemma}
\label{qu-2}
There exists   a   constant $\kappa_0\in (0,\nu]$  only depending   on
  $\nu,\{b_j\}_{j\in \cZ_0},\nu_S$ and $d=|\cZ_0|$
 such that the following statements hold:
  \begin{itemize}
    \item[(1)] For any $\kappa\in (0,\kappa_0]$ and the stopping time $\sigma$ defined in (\ref{2121-4}),
        we have
        \begin{eqnarray}
        \label{2727-1}
          \mE^{\mu_{\mS}} \exp\{10 \nu \sigma \}\leq C_\kappa,
        \end{eqnarray}
           where  $C_{\kappa}$ is a constant depending on
  $\kappa$ and  $\nu,\{b_j\}_{j\in \cZ_0},\nu_S,d.$ Hence
  \begin{eqnarray}\label{2727-1-000}
          \mE \exp\{10 \nu \sigma \}\leq C_\kappa.
        \end{eqnarray}
    \item[(2)]
    For any $\kappa\in (0,\kappa_0]$,  almost all  $\ell\in \mS$(under the measure $\mP^{\mu_\mS}$) and the   stopping times     $\sigma_k$ defined in (\ref{2121-4}) and  (\ref{qu-20}), we have
    \begin{eqnarray}
   \label{2626-2-000}
  \begin{split}
 & \mE^{\mu_{\mW}} \Big[ \exp\Big\{ \kappa \|w_{\sigma_k}\|^2 - \kappa \| w_{\sigma_{k-1}} \|^2
   e^{-1}
   \\ & \quad\quad \quad\quad  +
\nu \kappa   \int_{\sigma_{k-1}}^{\sigma_k}  e^{-\nu (\sigma_k-s)+8\mathfrak{B}_0 \kappa (\ell_{{\sigma_k}}-\ell_s) }  \| w_s\|_1^2
    \dif s
 \\ &  \quad\quad \quad\quad -  \kappa \mathfrak{B}_0 (\ell_{\sigma_k}- \ell_{\sigma_{k-1}})      \Big\}\Big|\cF_{\ell_{\sigma_{k-1}}}^{W}\Big]
  \leq C,
    \end{split}
\end{eqnarray}
 where $C$ is  a  constant   only depending   on
  $\nu,\{b_j\}_{j\in \cZ_0},\nu_S$ and $d.$
Moreover, the following statements hold:
    \begin{eqnarray}
   \label{2626-2}
  \begin{split}
 & \mE   \Big[ \exp\Big\{ \kappa \|w_{\sigma_k}\|^2 - \kappa \| w_{\sigma_{k-1}} \|^2
   e^{-1}
   \\ & \quad\quad \quad\quad  +
\nu \kappa   \int_{\sigma_{k-1}}^{\sigma_k}  e^{-\nu (\sigma_k-s)+8\mathfrak{B}_0 \kappa (\ell_{{\sigma_k}}-\ell_s) }  \| w_s\|_1^2
    \dif s
 \\ &  \quad\quad \quad\quad -  \kappa \mathfrak{B}_0 (\ell_{\sigma_k}- \ell_{\sigma_{k-1}})      \Big\}\Big|\cF_{\sigma_{k-1}} \Big]
  \leq C,
    \end{split}
\end{eqnarray}
 where $C$ is the  constant appearing in (\ref{2626-2-000}).

   \item[(3)]  For any $\kappa\in (0,\kappa_0]$
   and $k\in \mN,$ one has
   \begin{eqnarray}
   \label{100-2}
     \mE \Big[ \exp\big\{ \kappa \|w_{\sigma_{k+1} }\|^2\big\}\Big| \cF_{\sigma_{k}}\Big]\leq C_\kappa
     \exp\big\{ \kappa e^{-1}\|w_{\sigma_{k} }\|^2\big\},
   \end{eqnarray}
   where  $C_{\kappa}$ is a constant depending on
  $\kappa$ and  $\nu,\{b_j\}_{j\in \cZ_0},\nu_S,d.$
    \item[(4)]
For any $\kappa\in (0,\kappa_0]$, there exists
a $C_{\kappa}>0$  depending on
  $\kappa$ and  $\nu,\{b_j\}_{j\in \cZ_0},\nu_S,d$
such that for any  $ n\in \mN$ and $w_0\in H$, one has
\begin{eqnarray}
\label{30-7}
 \mE_{w_0} \exp\{\kappa \sum_{i=1}^n \|w_{\sigma_i}\|^2-C_\kappa n \}\leq e^{a\kappa \|w_0\|^2},
\end{eqnarray}
where $a=\frac{1}{1-e^{-1}}.$ In this paper, we
use the notation $\mE_{w_0}$ for expectations under the measure $\mathbb{P}$ with respect to solutions to (\ref{0.1}) with
initial condition $w_0$.
  \end{itemize}
   \item[(5)]
   For any $\kappa\in (0,\kappa_0]$, $w_0\in H$  and $n\in \mN,$ one has
  \begin{eqnarray}
  \label{100-1}
    \mE_{w_0} \sup_{s\in [0,\sigma]}\|w_s\|^{2n}\leq C_{n,\kappa} (1+\|w_0\|^{2n}),
  \end{eqnarray}
  where $C_{n,\kappa}$ is a constant depending on $n,\kappa$ and $ \nu,\{b_j\}_{j\in \cZ_0},\nu_S,d.$

\end{lemma}

 The proof of the above lemma is long, we leave it in Appendix \ref{ssss-1}.
 Throughout this paper, $\kappa_0$ is the constant appearing in Lemma \ref{qu-2}.

\subsection{Elements of Malliavin calculus}

Let   $d=|\cZ_0|$ and denote the canonical basis of $\mR^d$
by
$\{\theta_j\}_{j\in \cZ_0}.$
We have defined the linear operator {$Q: \mR^d\rightarrow H$ in the following way: for any $z=\sum_{j\in \cZ_0}z_j\theta_j \in \mR^d$,
  \begin{eqnarray*}
    Qz=\sum_{j\in \cZ_0 }   b_{j} z_j e_{j}.
  \end{eqnarray*}
  The adjoint of $Q$ is given by  $Q^*: H\rightarrow \mR^d$:
  \begin{eqnarray*}
   Q^*\xi=\big(b_{j}\langle \xi, e_{j}\rangle \big)_{j\in \cZ_0} \in \mR^d,  \text{ for } \xi\in H.
  \end{eqnarray*}

For any $0\le s\le t$ and $\xi\in  H$, let~$J_{s,t}\xi$  be the   solution of the linearised problem:
\begin{align}
\label{10-1}
	\partial_t  J_{s,t}\xi- \nu \Delta J_{s,t}\xi -\tilde{B}(w_t,J_{s,t}\xi)&=0,
  \\  J_{s,s}\xi&=\xi, \nonumber
\end{align}
where $\tilde{B}(w,v)=B(\cK w,v)+B(\cK v,w)$.

  For any $0\le t\le T$ and $\xi\in  H$, let~$K_{t,T}$  be the   adjoint of $J_{t,T}$, i.e., $\varrho_t:=K_{t,T}\xi$ satisfies the following backward  equation
\begin{align}
\label{w-1}
	\partial_t  \varrho _t+\nu \Delta \varrho_t + D\tilde B^*(w_t)\varrho_t&=0,
  \\ \varrho_T &=\xi, \nonumber
\end{align}
where $\langle D \tilde B^*(w)\rho, \psi\rangle=\langle \rho, D \tilde B(w)\psi\rangle$
and $D \tilde  B(w)\psi=B(\cK w,\psi)+B( \cK \psi,w)$.

Denote by $J^{(2)}_{s,t}(\phi,\psi)$ the second derivative of $w_t$ with respect to initial value  $w_0$ in the directions of $\phi$ and $\psi$. Then
\begin{eqnarray}
\label{0927-1}
\left\{
\begin{split}
  & \partial_t J^{(2)}_{s,t}(\phi,\psi)= \nu \Delta J^{(2)}_{s,t}(\phi,\psi)+ B(\cK  J_{s,t}\phi,J_{s,t}\psi)+B(\cK J_{s,t}\psi,J_{s,t}\phi)
  \\ &\quad\quad +B(\cK w_t,J^{(2)}_{s,t}(\phi,\psi))+
  B(\cK J^{(2)}_{s,t}(\phi,\psi),w_t), \quad \text{ for } t>s,
  \\ & J^{(2)}_{s,s}(\phi,\psi)=0.
  \end{split}
  \right.
\end{eqnarray}

For  given $\ell\in \mS,$  $t>0$, let $\Phi(t,W)$ be a ${\cal F}_{\ell_t}^W$-measurable random variable.  For $v\in L^2([0,\ell_t];\R^d),$
the Malliavin derivative of~$\Phi$ in the direction~$v$ is defined by
$$
\cD^v\Phi(t,W)=\lim_{\eps \to  0}\frac{1}{\e}
  \left(\Phi(t,w_0,W+\eps \int_0^\cdot v\dd  s)-\Phi(t,w_0,W)\right), \quad
$$
  where the limit   holds almost surely (e.g., see the book~\cite{nualart2006} for finite-dimensional setting or the papers~\cite{MP-2006, HM-2006, HM-2011, FGRT-2015}   for   Hilbert space~case).
 {In the definition of  Malliavin derivative,  the element $\ell$ is   taken as fixed.}
  Then, $\cD^v w_t$ satisfies the following equation:
  \begin{eqnarray*}
    \dif \cD^v w_t =\nu \Delta \cD^v w_t\dif t+\tilde B(\cD^v w_t,w_t)\dif t+
    Q  \dif \left(\int_0^{\ell_t} v_s\dif s\right).
  \end{eqnarray*}
  By the Riesz representation theorem, there is a linear operator $  \cD:L^2(\Omega,   H)\to L^2(\Omega; L^2([0,\ell_t ];\R^d)\otimes  H)$ such that
\begin{equation}\label{2.3}
	  \cD^v w_t  =\lag    \cD w,v \rag_{L^2([0,\ell_t];\R^d)},~\forall v\in L^2([0,\ell_t];\R^d) .
\end{equation}
Actually, we have the following lemma.
  \begin{lemma}
\label{17-1}
For  any   $\ell \in \mS$ and    $v\in L^2([0,\ell_t];\R^d),$
we have
\begin{align*}
  \cD^v w_t&=\int_0^{\ell_t}  J_{\gamma_u,t}Q v_u \dif u,
\end{align*}
here   $\gamma_u$ is defined by
$
    \gamma_u=\inf\{t\geq 0,S_t(\ell)\geq u\}.
$
Hence,   we also have
  \begin{eqnarray*}
    \cD_u^jw_t=J_{\gamma_u,t}  Q \theta_j, ~\forall u\in [0,\ell_t],
  \end{eqnarray*}
  where $\cD_u^j$ denotes the Malliavin derivative with respect to the $j$th component of the noise at time $u.$
\end{lemma}
\begin{proof}
  We need to prove that for any $v\in L^2([0,\ell_t],\mR^d)$,
\begin{eqnarray}
\label{1607-1}
  && \cD^v w_t=\int_0^t
   J_{r,t}Q \dif
   \big( \int_0^{\ell_r}v_s\dif s\big)
   =\int_0^{\ell_t} J_{\gamma_u,t}Q v_u \dif u.
\end{eqnarray}
The first equality  in (\ref{1607-1}) follows from the formula of constant variations or Fubini's theorem; see \cite[Page 370, lines 1--5]{Zhang-2016-CMS} for example.
So we give a proof of the second equality  in (\ref{1607-1}).
Obviously, we have
\begin{eqnarray}
\label{12-1}
\int_0^t
   J_{r,t}Q \dif
   \big( \int_0^{\ell_r}v_s\dif s\big)
   =\sum_{r\leq t } J_{r,t} Q \int_{\ell_{r-}}^{\ell_r} v_s\dif s
\end{eqnarray}
Since  $\gamma_u=r,  u\in (\ell_{r-},\ell_r)$,  it holds that
\begin{eqnarray*}
  \sum_{r\leq t } J_{r,t}Q  \int_{\ell_{r-}}^{\ell_r} v_u\dif u=  \sum_{r\leq t } \int_{\ell_{r-}}^{\ell_r} J_{\gamma_u,t}Q v_u\dif u
  =\int_0^{\ell_t} J_{\gamma_u,t}Q v_u\dif u.
\end{eqnarray*}
Combining this with  (\ref{12-1}), we complete the proof of  the second equality  in (\ref{1607-1}).

\end{proof}

For any $s\leq t$ and $\ell\in \mS,$  define  the linear operator
$\cA_{s,t}: L^2([\ell_s,\ell_t];\mR^d)\rightarrow H  $
by
\begin{eqnarray}
\label{15-1}
  \cA_{s,t}v=\int_{\ell_s}^{\ell_t}
   J_{\gamma_u,t}Q v_u  \dif u, ~~v\in L^2([\ell_s,\ell_t];\mR^d).
\end{eqnarray}
The adjoint of $\cA_{s,t}$, $\cA_{s,t}^*:H\rightarrow  L^2([\ell_s,\ell_t];\mR^d)$, is given by
\begin{eqnarray*}
\cA_{s,t}^*\phi=\big(b_j \langle \phi,
J_{\gamma_u,t}e_j\rangle\big)_{j\in \cZ_0,
 u\in [\ell_s,\ell_t]}.
\end{eqnarray*}
The Malliavin matrix  $\cM_{s,t}:H\rightarrow H $ is defined by
\begin{eqnarray}
\label{200-3}
  \cM_{s,t}\phi &=& \cA_{s,t}\cA_{s,t}^*\phi
\end{eqnarray}
By a simple calculation, we have(c.f.  \cite[Lemma 2.2]{Zh14})
\footnote{Actually,  we have,
\begin{eqnarray*}
\langle \cM_{s,t}\phi,\phi\rangle
  &=& \sum_{j\in \cZ_0 }b_j^2\int_{\ell_s}^{\ell_t} \langle K_{\gamma_u,t}\phi, e_j\rangle^2 \dif u=\sum_{j\in \cZ_0 }b_j^2\sum_{r\in (s,t]}
  \int_{\ell_{r-}}^{\ell_r} \langle K_{\gamma_u,t}\phi, e_j\rangle^2 \dif u
 \\  &=& \sum_{j\in \cZ_0 }b_j^2\sum_{r\in (s,t]}
   \langle K_{r,t}\phi, e_j\rangle^2 (\ell_r-\ell_{r-})=\sum_{j\in \cZ_0 }b_j^2\int_s^t \langle K_{r,t}\phi, e_j\rangle^2 \dif \ell_r.
   \end{eqnarray*}

}
\begin{eqnarray}
\label{200-1}
  \langle \cM_{s,t}\phi,\phi\rangle
  =\sum_{j\in \cZ_0 }b_j^2\int_s^t \langle K_{r,t}\phi, e_j\rangle^2 \dif \ell_r.
\end{eqnarray}
\textbf{(For any function $f:[a,b]\rightarrow \mR,$ the integral $\int_a^b f(s)\dif \ell_s$ is interpreted as $\int_{(a,b]}f(s)\dif \ell_s.$ )}

In the rest of this subsection, we will provide  some estimates for $J_{s,t},J_{s,t}^{(2)},\cA_{s,t},$ $\cA_{s,t}^*,$$\cM_{s,t}$  and their Malliavin derivatives, which will be used in subsequent sections.


\begin{lemma}
\label{15-4}
There exists  a   constant  $\cC_0$  only  depending on   $ \nu$
 such that  for any $\xi,\phi,\psi\in H$  and $0\leq s\leq t\leq T,$
 $J_{s,t}$ and $J_{s,t}^{(2)}$ satisfy almost surely
 \begin{eqnarray}
  \label{87-10} \sup_{t\in [s,T]}\|J_{s,t}\xi\|^2 &\leq &
  \cC_0 \|\xi\|^2  e^{\cC_0  \int_{s}^T  \|w_r\|_1^{4/3}\dif r},
   \\ \label{87-11}
   \int_s^t \|J_{s,r}\xi\|_1^2 \dif r
   &\leq & \cC_0  \|\xi\|^2  e^{\cC_0  \int_{s}^t  \|w_r\|_1^{4/3}\dif r} ,
   \\  \label{87-12}
   \sup_{ t \in [s, T]}\|J^{(2)}_{s,t}(\phi,\psi)\|^2  & \leq &  \cC_0    \|\phi\|^2\|\psi\|^2
   e^{\cC_0   \int_{s}^T  \|w_r\|_1^{4/3}\dif r}.
 \end{eqnarray} 
Furthermore, for any $\kappa>0$ and $0\leq s\leq T,$  it also   holds that
    \begin{eqnarray}
    \label{1340-1}
    \begin{split}
    & \|J_{s,T}\xi\|^2
    \leq
     \cC_0 \exp\Big\{\frac{\nu \kappa}{120} \int_s^T \|w_s\|_1^{2}e^{-\nu (T-s)+8\mathfrak{B}_0\kappa (\ell_T-\ell_s)} \dif s
     \\ & \quad\quad\quad\quad\quad\quad \quad\quad  +C_\kappa  \int_s^T e^{ 2\nu (T-s)-16\mathfrak{B}_0 \kappa (\ell_T-\ell_s)}  \dif s \Big\}\|\xi\|^2,
     \end{split}
  \end{eqnarray}
  where $\cC_0$ is taken from (\ref{87-10})--(\ref{87-12}),  and $C_\kappa$ is a constant  depending on $\kappa, \nu.$
\end{lemma}

\begin{remark}
In the other places of this paper, we need to use the constants  appeared in  (\ref{87-10})--(\ref{87-12}), so we use a special notation $\cC_0$ to denote them.
\end{remark}
\begin{proof}
By (\ref{44-1}),
\begin{eqnarray}
\label{0927-2}
\begin{split}
 &  \langle B(\cK J_{s,t}\xi,w_t ),J_{s,t}\xi \rangle
  \leq C\|w_t\|_1 \|J_{s,t}\xi\|_{1/2}\|J_{s,t}\xi\|
  \\ & \leq  \frac{\nu}{4}\|J_{s,t}\xi\|_1^2+ C\|w_t\|_{1}^{4/3}\| J_{s,t}\xi\|^2,
  \end{split}
\end{eqnarray}
where $C=C(\nu).$  Therefore, applying the chain rule to $\|J_{s,t}\xi\|^2,$ one arrives at
\begin{eqnarray*}
  \dif \|J_{s,t}\xi\|^2 \leq -\nu \|J_{s,t}\xi\|^2_1\dif t+C\|w_t\|_{1}^{4/3}\| J_{s,t}\xi\|^2\dif t,
\end{eqnarray*}
which implies
\begin{eqnarray*}
 \|J_{s,t}\xi\|^2\leq  C \|\xi\|^2e^{C\int_s^t \|w_r\|_{1}^{4/3}\dif r},\quad \forall 0\leq s\leq t
\end{eqnarray*}
and
\begin{eqnarray*}
 &&  \nu \int_s^t  \|J_{s,r}\xi\|_1^2\dif r \leq
 \|\xi\|^2+  C\int_s^t \|w_r\|_{1}^{4/3}\| J_{s,r}\xi\|^2\dif r
  \\ &&\leq \|\xi\|^2+ C  \|\xi\|^2 \int_s^t \|w_r\|_{1}^{4/3} \dif r
  e^{C\int_s^t \|w_r\|_{1}^{4/3}\dif r}
  \leq C  \|\xi\|^2
  e^{C\int_s^t \|w_r\|_{1}^{4/3}\dif r}.
\end{eqnarray*}
The proof of (\ref{87-10}) and (\ref{87-11}) is complete.

Now we prove (\ref{87-12}). As in (\ref{0927-2}), we have
\begin{eqnarray*}
  \langle B(\cK J^{(2)}_{s,t}(\phi,\psi),w_t) , J^{(2)}_{s,t}(\phi,\psi)
  \rangle\leq \frac{\nu}{4}\|J^{(2)}_{s,t}(\phi,\psi) \|_1^2+ C\|w_t\|_{1}^{4/3}\| J^{(2)}_{s,t}(\phi,\psi) \|^2.
\end{eqnarray*}
Moreover,
\begin{eqnarray*}
 &&  \langle B(\cK  J_{s,t}\phi,J_{s,t}\psi), J^{(2)}_{s,t}(\phi,\psi) \rangle\leq C \|J^{(2)}_{s,t}(\phi,\psi)\|_1 \|J_{s,t}\phi\|_{1/2}  \|J_{s,t}\psi\|
 \\ &&\leq
 \frac{\nu}{4} \|J^{(2)}_{s,t}(\phi,\psi)\|_1^2+C \|J_{s,t}\phi\|_{1/2}^2  \|J_{s,t}\psi\|^2 .
\end{eqnarray*}
Hence, applying the chain rule to $\|J^{(2)}_{s,t}(\phi,\psi)\|^2,$
with the help of  (\ref{87-10})--(\ref{87-11}), we get
\begin{eqnarray*}
 &&  \|J^{(2)}_{s,t}(\phi,\psi)\|^2\! \leq  \! C e^{C\int_s^t \|w_r\|_1^{4/3}\dif r}\!  \int_s^t
  \! \[\|J_{s,r}\phi\|_{1/2}^2  \|J_{s,r}\psi\|^2+\|J_{s,r}\psi\|_{1/2}^2  \|J_{s,r}\phi\|^2 \]\! \dif r
  \\ &&\leq  C e^{C\int_s^t \|w_r\|_1^{4/3}\dif r} \|\phi\|^2\|\psi\|^2.
\end{eqnarray*}

(\ref{1340-1}) is  easily  obtained by Young's inequality and (\ref{87-10}).

\end{proof}

Recall that $P_N$ is the orthogonal projection from $H$ into $H_N=span\{e_j; j\in Z_{\ast}^2, |j|\leq N\}$ and $Q_N=I-P_N$.
For any $N\in \mN,t\geq 0$ and $\xi\in H,$  denote $\xi_t^h=Q_NJ_{0,t}\xi,
\xi_t^\iota =P_NJ_{0,t}\xi$ and $\xi_t=J_{0,t}\xi.$
\begin{lemma}
\label{1340-3}
 For any $t\geq 0$ and $\xi\in H,$  one has
 \begin{eqnarray*}
   \|\xi_t^h\|^2 \leq \exp\{-\nu N^2 t\}\|\xi\|^2
   +\frac{C\|\xi\|^2}{\sqrt{N}} \exp\{C\int_0^t \|w_s\|_1^{4/3}\dif s\} \sup_{s\in [0,t]}\|w_s\|,
 \end{eqnarray*}
 where $C$  is  a    constant   depending on
 $\nu.$
\end{lemma}
\begin{proof}
Note that $\|\xi_t^h\|^2_1\geq N^2\|\xi_t^h\|^2$  and that
\begin{eqnarray*}
  && \langle B(\cK \xi_t,w_t),\xi_t^h\rangle+\langle B(\cK w_t,\xi_t),\xi_t^h\rangle\leq
  C\|\xi_t^h\|_1 \|w_t\|_{1/2}  \|\xi_t\|
  \\ && \leq \frac{\nu}{4}\|\xi_t^h\|_1^2+C \|w_t\|_{1/2}^2  \|\xi_t\|^2,
\end{eqnarray*}
applying the chain rule  to $\|\xi_t^h\|^2,$  we find
\begin{eqnarray*}
  && \|\xi_t^h\|^2 \leq \exp\{-\nu N^2 t\}\|\xi\|^2+
  C \int_0^t \exp\{-\nu N^2 (t-s) \} \|w_s\|_{1/2}^2  \|\xi_s\|^2\dif s
  \\ &&\leq \exp\{-\nu N^2 t\}\|\xi\|^2
  \\ &&\ \   + C \exp\{C \int_0^t \|w_s\|_1^{4/3}\dif s\} \|\xi\|^2
   \sup_{s\in [0,t]} \|w_s\| \int_0^t \exp\{-\nu N^2 (t-s) \}\|w_s\|_{1}  \dif s
    \\ &&\leq \exp\{-\nu N^2 t\}\|\xi\|^2+
   C \exp\{C \int_0^t \|w_s\|_1^{4/3}\dif s\}\|\xi\|^2
   \sup_{s\in [0,t]} \|w_s\|
    \\ && \quad\quad\quad\quad \times \big(\int_0^t \exp\{-4 \nu N^2 (t-s) \} \dif s\big)^{1/4} \big(\int_0^t\|w_s\|_{1}^{4/3}  \dif s \big)^{3/4},
\end{eqnarray*}
where we have used  (\ref{87-10}) in the second  inequality. The above
inequality  implies the desired result.
\end{proof}

\begin{lemma}
\label{1340-2}
  Assume that $\xi_0^\iota=0$, then for any $t\geq 0,$
  \begin{eqnarray*}
   \|\xi_t^\iota\|^2 \leq   C \|\xi \|^2  \exp\{C\int_0^t \|w_s\|_1^{4/3}\dif s\}
 \sup_{s\in [0,t]}(\|w_s\|^{5/2}+1)
 \frac{1+t}{N^{1/4}},
  \end{eqnarray*}
   where $C$  is  a    constant   depending on
 $\nu.$ 
 Furthermore, combining the above  inequality with Lemma  \ref{1340-3}, for any $\xi\in H$ and $t\geq 0,$
  we have
    \begin{eqnarray*}
  &&  \|J_{0,t}Q_N \xi\|^2
 \\ && \leq C \big(e^{-\nu N^2 t}+\frac{1+t}{N^{1/4}} \big)
 \exp\{C\int_0^t \|w_s\|_1^{4/3}\dif s\}
 \sup_{s\in [0,t]}(\|w_s\|^{5/2}+1)\|\xi\|^2.
  \end{eqnarray*}
\end{lemma}
\begin{proof}
In view of (\ref{44-1}),  one has
\begin{eqnarray*}
  \langle B(\cK \xi_t^\iota,w_t),\xi_t^\iota\rangle
  \leq \frac{\nu}{4} \|\xi_t^\iota \|_1^2+C\|w_t\|_1^{4/3}\|\xi_t^\iota\| ^2
\end{eqnarray*}
(See (\ref{0927-2}) for similar arguments)  and
\begin{eqnarray*}
&&  \langle B(\cK w_t,\xi_t^h),\xi_t^\iota\rangle
 +  \langle B(\cK \xi_t^h, w_t),\xi_t^\iota\rangle
 \\ && \leq  C \|\xi_t^\iota \|_1 \|w_t\|\|\xi_t^h \|_{1/2}
 \leq \frac{\nu}{4} \|\xi_t^\iota \|_1^2+C\|w_t\|^2 \|\xi_t^h \|_{1/2}^2.
\end{eqnarray*}
Thus, applying the chain rule to $\|\xi_t^\iota\|^2,$  we have
\begin{eqnarray*}
 && \|\xi_t^\iota \|^2 \leq C \exp\{C\int_0^t \|w_s\|_1^{4/3}\dif s\}
 \int_0^t \|w_s\|^2 \|\xi_s^h \|_{1/2}^2\dif s
 \\  && \leq    C \exp\{C\int_0^t \|w_s\|_1^{4/3}\dif s\}
 \sup_{s\in [0,t]}\|w_s\|^2
 \int_0^t  \|\xi_s^h \|   \|\xi_s^h \|_{1} \dif s
  \\  && \leq    C \exp\{C\int_0^t \|w_s\|_1^{4/3}\dif s\}
 \sup_{s\in [0,t]}\|w_s\|^2
 \Big(\int_0^t  \|\xi_s^h \|^2 \dif s \Big)^{1/2}
 \Big( \int_0^t \|\xi_s^h \|_{1}^2  \dif s\Big)^{1/2}
   \\  && \leq    C \|\xi \| \exp\{C\int_0^t \|w_s\|_1^{4/3}\dif s\}
 \sup_{s\in [0,t]}\|w_s\|^2
 \Big(\int_0^t  \|\xi_s^h \|^2 \dif s \Big)^{1/2}
 \\  && \leq    C \|\xi \| \exp\{C\int_0^t \|w_s\|_1^{4/3}\dif s\}
 \sup_{s\in [0,t]}\|w_s\|^2
 \\ && \quad  \times \Big(\int_0^t \big[ \exp\{-\nu N^2 s\}\|\xi\|^2
   +\frac{C\|\xi\|^2}{\sqrt{N}} \exp\{C\int_0^s \|w_r\|_1^{4/3}\dif r\} \sup_{r\in [0,s]}\|w_r\| \big] \dif s \Big)^{1/2}
    \\  && \leq    C \|\xi \|^2  \exp\{C\int_0^t \|w_s\|_1^{4/3}\dif s\}
 \sup_{s\in [0,t]}\|w_s\|^2
 \Big(\frac{1}{N}+\frac{\sup_{s\in [0,t]}\|w_r\|^{1/2} \sqrt{t} }{N^{1/4}}\Big)
   \\  && \leq    C \|\xi \|^2  \exp\{C\int_0^t \|w_s\|_1^{4/3}\dif s\}
 \sup_{s\in [0,t]}(\|w_s\|^{5/2}+1)
 \frac{1+t}{N^{1/4}},
\end{eqnarray*}
where in the fourth  inequality we have used (\ref{87-11}) for the fourth  inequality and
Lemma \ref{1340-3} for the fifth inequality.
\end{proof}

Using the similar arguments as that in   \cite[Section~4.8]{HM-2006}  or   \cite[Lemma~A.6]{FGRT-2015}, we have the following lemma.
\begin{lemma}
\label{L:2.2}There is a   constant $C=C(\{b_j\}_{j\in \cZ_0},d)>0$ such that  for any
  $0\le s<t$ and  $\beta>0$, we have
  \begin{gather}
   \label{2.7}  \|\cA_{s,t}\|^2_{\cL(L^2([\ell_s,\ell_t];\R^d),
   {H})} \le  C \int_{\ell_s}^{\ell_t}\|J_{\gamma_u,t}\|_{\cL(H,{H})}^2\dd u,
     \\   \|\cA_{s,t}^*(\cM_{s,t}+\beta\I)^{-1/2}\|_{\cL({H},
     L^2([\ell_s,\ell_t];\R^d))} \le  1,\label{2.8}
   \\
   \|(\cM_{s,t}+\beta\I)^{-1/2}\cA_{s,t}\|_{\cL(L^2([\ell_s,\ell_t];\R^d),
   {H})} \le  1,\label{2.9}
   \\  \|(\cM_{s,t}+\beta\I)^{-1/2}\|_{\cL({H},{H})}
   \le  \beta^{-1/2}.\label{2.10}
  \end{gather}
\end{lemma}

\begin{lemma}
\label{15-5}
  For any $0\leq s\leq t, j\in \{1,\cdots,d\}$ and $u\in [0,\ell_t],$ we have
\begin{eqnarray*}
 \cD_u^jJ_{s,t}\xi =
 \left\{
 \begin{split}
 & J^{(2)}_{\gamma_u,t}(Q\theta_j,J_{s,\gamma_u}\xi),
 u\in [\ell_s,\ell_{t}]
 \\
 &
 J^{(2)}_{s,t}(J_{\gamma_u,s}Q\theta_j,\xi)
\quad \text{if } u< \ell_s.
 \end{split}
 \right.
\end{eqnarray*}
\end{lemma}
\begin{proof}
In view of Lemma  \ref{17-1},
the  proof  is the same as that in \cite[(4.29)]{HM-2006}.
We omit the details.
\end{proof}

As in  ~\cite{HM-2006}, if $A:\cH_1\rightarrow \cH_2$ is a random linear map between two Hilbert spaces, we denote by
$\cD_s^iA:\cH_1\rightarrow \cH_2 $ the random linear map defined by
\begin{eqnarray*}
  (\cD_s^iA)h=\langle \cD_s(Ah),\theta_i\rangle.
\end{eqnarray*}
\begin{lemma}\label{L:2.3}
The   operators $J_{s,t}$, $\cA_{s,t}$, and $\cA_{s,t}^*$
are Malliavin differentiable, and for any
$r>0$, the following inequalities~hold
\begin{align}
  \label{2.11}
    \[ \|\cD_r^iJ_{0,\sigma}\|_{\cL({H},{H})}
  & \le    C_\kappa    \exp\{\cC_0 \int_0^\sigma \|w_s\|_1^{4/3}\dif s \}\],
  \\  \label {2.12}
      \|\cD_r^i \cA_{0,\sigma}\|_{\cL(L^2([0,\ell_\sigma];\R^d),
   {H})}  & \le   C_\kappa   (\sigma+1)    \exp\{\cC_0 \int_0^\sigma \|w_s\|_1^{4/3}\dif s \},
  \\   \label{2.13}
    \|\cD_r^i\cA^*_{0,\sigma}\|_{\cL({H},
  L^2([0,\ell_\sigma];\R^d))}  & \le   C_\kappa  (\sigma+1)    \exp\{\cC_0 \int_0^\sigma \|w_s\|_1^{4/3}\dif s \},
\end{align}
where $\cC_0$ is the  same  constant as that appeared in     Lemma \ref{15-4},
$C_\kappa$ is a  constant  depending on $\kappa,\nu,\{b_j\}_{j\in \cZ_0},d$
(Recall that $\sigma$ depends on $\kappa.$).
\end{lemma}
\begin{proof}
  The inequality (\ref{2.11})
   is  derived from    Lemma \ref{15-4} and Lemma  \ref{15-5}.
(\ref{2.12}) is  obtained by  Lemmas  \ref{15-4}, \ref{15-5}, Cauchy-Schwarz inequality  and  the fact that
  \begin{eqnarray*}
  \cA_{0,\sigma}v=\int_0^{\ell_\sigma} J_{\gamma_u,\sigma}Q v(u)\dif u, \quad \ell_\sigma\leq \frac{\nu \sigma}{8\mathfrak{B}_0 \kappa  }.
  \end{eqnarray*}
The inequality (\ref{2.13}) is a consequence of (\ref{2.12}).
\end{proof}

\section{The invertibility of the Malliavin matrix $\cM_{0,t}.$}
\label{S:3}


Before stating the main results in  this section, we prepare two lemmas. Lemma \ref{31-1}
can be seen as the pure jump version of Theorem 7.1 in \cite{HM-2011}, which deal with the Wiener case.
 { In this section, we use $\Delta f(s)$ to denote $f(s)-f(s-)$.}

\begin{lemma}\label{Lemma 2024 Number jump}
Consider  a probability space $(\widetilde{\Omega},\widetilde{\mathcal{F}},\widetilde{\mathbb{P}})$,  and a given L\'evy process $\widetilde{L}(t),t\geq0$ on $(\widetilde{\Omega},\widetilde{\mathcal{F}},\widetilde{\mathbb{P}})$, which takes values in a topological vector space $(\mathcal{T},B(\mathcal{T}))$. Suppose that the L\'evy process $\widetilde{L}(t),t\geq0$ has a $\sigma$-finite  intensity measure $\widetilde{\nu}$.
 For any $G\in B(\mathcal{T})$, define $N^G_{\widetilde{L}}((t_1,t_2])=\sharp\{s\in(t_1,t_2]:~\Delta\widetilde{L}(s)\in G\}$. If $\widetilde{\nu}(G)=\infty$, then there exists a
 $\widetilde{\Omega}_0\in\widetilde{\mathcal{F}}$ with $\widetilde{\mathbb{P}}(\widetilde{\Omega}_0)=1$ such that for any $\widetilde{\omega}\in \widetilde{\Omega}_0$ and any $0\leq t_1<t_2$,
 \begin{eqnarray*}
   N^G_{\widetilde{L}}((t_1,t_2])(\widetilde{\omega})=\infty.
 \end{eqnarray*}
 Moreover, $\{s\in[0,\infty):~\Delta\widetilde{L}(s)(\widetilde{\omega})\in G\}$ is  dense on $[0,\infty)$.
\end{lemma}

\begin{proof}
Let $\mathbb{Q}$ denote the set of all rational number on $\mathbb{R}$. It is sufficient to prove that for any $0\leq t_1<t_2$ with $t_1, t_2\in \mathbb{Q}$, we have
$\widetilde{\mathbb{P}}(N^G_{\widetilde{L}}((t_1,t_2])=\infty))=1$. Let $K_n, n\geq 1$ be an increasing sequence of measurable subsets of $\mathcal{T}$ such that $K_n\uparrow \mathcal{T}$ and $\widetilde{\nu}(K_n)<\infty$. It is well known that $N^{G\cap K_n}_{\widetilde{L}}((t_1,t_2])$ is a Poisson random variable with parameter
$\widetilde{\nu}(G\cap K_n)(t_2-t_1)$. To prove $\widetilde{\mathbb{P}}(N^G_{\widetilde{L}}((t_1,t_2])=\infty))=1$, it suffices to show that for any positive integer $M$, $\widetilde{\mathbb{P}}(N^G_{\widetilde{L}}((t_1,t_2])>M))=1$.
Indeed,
\begin{eqnarray*}
&&\widetilde{\mathbb{P}}(N^G_{\widetilde{L}}((t_1,t_2])>M)=\lim_{n\rightarrow \infty}\widetilde{\mathbb{P}}(N^{G\cap K_n}_{\widetilde{L}}((t_1,t_2])>M)\nonumber\\
&=&1-\lim_{n\rightarrow \infty}\widetilde{\mathbb{P}}(N^{G\cap K_n}_{\widetilde{L}}((t_1,t_2])\leq M)\nonumber\\
&=&1-\lim_{n\rightarrow \infty}\sum_{m=0}^Mexp(-\nu(G\cap K_n)(t_2-t_1))\frac{(\widetilde{\nu}(G\cap K_n)(t_2-t_1))^m}{m!}=1,
\end{eqnarray*}
where we used the fact that $\widetilde{\nu}(G\cap K_n)(t_2-t_1)\rightarrow \infty$, as $n\rightarrow \infty$.
\end{proof}

Consider the probability space $(\tilde{\Omega},\tilde{\mathcal{F}},\tilde{\mathbb{P}})$, and a $\mR^d$-valued  L\'evy process $\tilde{L}(t)=( \tilde{L}^1 (t),
   \tilde{L}^2 (t),....,\tilde{L}^d (t)), t\geq 0$
   with a $\sigma$-finite intensity measure $\widetilde{\nu}$.
  For any $n\in \mathbb{N}$ and $1\leq i\leq  d$, let
  $$G_n^i=\{y\in\mathbb{R}^d\setminus\{0\}:\  { \max_{j\in \{1,\cdots,d \}\text{ with }j\neq i}|y_j|< \frac{|y_i|}{n}}\}.
  $$
  Assume that $\widetilde{\nu}$ satisfies the following condition:
  \begin{eqnarray*}
    \widetilde{\nu}(G_n^i )=\infty, \quad \forall n\in \mN\text{ and }1\leq i\leq d.
  \end{eqnarray*}
 By Lemma \ref{Lemma 2024 Number jump}, there exists a
$\Omega^i_n\in \tilde{\mathcal{F}} $ with $\tilde{\mathbb{P}}(\Omega^i_n)=1$ such that for any $\omega\in \Omega^i_n$, the set   $\{s\in[0,\infty):~\Delta \tilde L({\omega},s)\in G_n^i\}$ is  dense in $[0,\infty)$.
Let $\Omega_0 :=\cap_{i=1}^d \cap_{n\in\mathbb{N}}\Omega^i_n$, then $\tilde{\mathbb{P}}(\Omega_0 )=1$, and for any $\omega\in \Omega_0$ the set $\{s\in[0,\infty):~\Delta \tilde L({\omega},s)\in G_n^i\}$ is  dense in $[0,\infty)$ for any $1\leq i \leq d$ and $n\in\mathbb{N}$.
We stress that $\Omega_0$ only depends on   the  L\'evy process $\tilde{L}(t), t\geq0$.

\begin{lemma}
\label{31-1}
  If for some $\omega_0\in\Omega_0$,
the following three conditions are satisfied:
  \begin{itemize}
  \item[(1)] $a(\omega_0),b(\omega_0)\in[0,\infty)$ and $a(\omega_0)<b(\omega_0)$;

  \item[(2)] $g_i(\omega_0, \cdot):[a(\omega_0),b(\omega_0)]\rightarrow \mR,0\leq i\leq d,$ are continuous functions;

  \item[(3)]   \begin{eqnarray}
  \label{56-1}
    g_0(\omega_0,r)+\sum_{i=1}^d g_i(\omega_0,r)\tilde{L}^i (\omega_0,r)=0, \quad \forall r\in [a(\omega_0),b(\omega_0)].
  \end{eqnarray}
  \end{itemize}
  Then
    \begin{eqnarray*}
    g_i(\omega_0,r)=0,\quad \forall r\in [a(\omega_0),b(\omega_0)], \quad  0\leq i\leq d.
  \end{eqnarray*}

\end{lemma}
\begin{proof}
By (\ref{56-1}), it is sufficient to show that
\begin{eqnarray}\label{eq 2024 05}
  g_i(\omega_0,r)=0,\quad \forall r\in [a(\omega_0),b(\omega_0)],
\end{eqnarray}
for $i=1, \cdots, d$.
Let us prove
\begin{eqnarray*}
  g_1(\omega_0,r)=0,\quad r\in [a(\omega_0),b(\omega_0)].
\end{eqnarray*}
The proofs of the other  cases with $i=2, \cdots, d$ are similar.

The conditions (2) and (3) imply that
\begin{eqnarray}\label{eq 2024 06}
0  = \sum_{i=1}^d  \Delta \tilde L^i(\omega_0,r) g_i(\omega_0,r), \quad r\in [a(\omega_0),b(\omega_0)].
\end{eqnarray}

Fix $n\in\mathbb{N}$. For any $s\in \{s\in[0,\infty):~\Delta \tilde L({\omega_0},s)\in G_n^1\}\cap [a(\omega_0),b(\omega_0)],$
by (\ref{eq 2024 06}) and the definition of $G_n^1$, one has
\begin{eqnarray*}
0 & =& |\sum_{i=1}^d  \Delta \tilde L^i(\omega_0,s) g_i(\omega_0,s)|
\\ &\geq &  |\Delta \tilde L^1(\omega_0,s)| \cdot | g_1(\omega_0,s)|-\sum_{i=2}^d  | \Delta \tilde L^i(\omega_0,s)| \cdot | g_i(\omega_0,s)|
\\ &\geq & |\Delta \tilde L^1(\omega_0,s)| \cdot | g_1(\omega_0,s)|-\frac{d}{n} | \Delta \tilde L^1(\omega_0,s)| \cdot \sum_{i=2}^d| g_i(\omega_0,s)|,
\end{eqnarray*}
which implies
\begin{eqnarray*}
 | g_1(\omega_0,s)|
\leq \frac{d}{n}  \sum_{i=2}^d| g_i(\omega_0,s)|,
\end{eqnarray*}
where we have used the fact that $| \Delta \tilde L^1(\omega_0,s)|>0$.
Recall that  the definition of $\Omega_0$ implies that
the set   $\{s\in[0,\infty):~\Delta \tilde L({\omega_0},s)\in G_n^1\}$ is  dense on $[0,\infty)$ for any $n\in\mathbb{N}$.
By the continuity  of $g_i(\omega_0,\cdot), i=1,2,...,d$, we obtain
\begin{eqnarray*}
 \sup_{s\in [a(\omega_0),b(\omega_0)]}| g_1(\omega_0,s)|
\leq \frac{d}{n}  \sup_{s\in [a(\omega_0),b(\omega_0)]}\sum_{i=2}^d| g_i(\omega_0,s)|,\quad \forall n\in \mN.
\end{eqnarray*}
Since $n$ is arbitrary, we obtain that
\begin{eqnarray*}
g_1(\omega_0,s)=0,\quad \forall s\in [a(\omega_0),b(\omega_0)].
\end{eqnarray*}
The proof is complete.
\end{proof}

Recall the assumption (\ref{Condition nu S}): $\nu_S((0,\infty))=\infty$.  By Lemma \ref{Lemma 2024 Number jump},  for the process $S_t, t\geq 0$ , we have the following result.
\begin{lemma}
\label{25-1}
 \begin{eqnarray*}
   \mP^{\mu_\mS}\big(\ell:\{s:\Delta S_s(\ell)>0\} \text{ is dense in } (0,\infty)\big)=1.
 \end{eqnarray*}
\end{lemma}

Recall that
\begin{eqnarray*}
  \langle \cM_{0,\sigma}\phi,\phi\rangle
  =\sum_{j\in \cZ_0}b_j^2\int_0^\sigma \langle K_{r,\sigma}\phi, e_j\rangle^2 \dif \ell_r.
\end{eqnarray*}



The first objective  this section is to prove the following proposition.

\begin{proposition}
\label{1-66}
  For any $\alpha\in (0,1],N\in \mN$ and  $w_0\in H,$   one has
  \begin{eqnarray}
  \label{1-2}
    \mP\Big(\inf_{\phi\in \cS_{\alpha,N}}\langle
    \cM_{0,\sigma}\phi,\phi\rangle=0 \Big)=0,
  \end{eqnarray}
  where $\cS_{\alpha,N}=\{\phi:\|P_N\phi\|\geq \alpha, \|\phi\|=1\}.$
\end{proposition}
We will prove the following stronger result than
(\ref{1-2}) for later use:
  \begin{eqnarray}
  \label{86-1}
    \mP\Big(\omega=(\text{w},\ell): \inf_{\phi\in \cS_{\alpha,N}}\sum_{j\in \cZ_0}b_j^2\int_{\sigma/2}^\sigma \langle K_{r,\sigma}\phi, e_j\rangle^2 \dif \ell_r =0 \Big)=0.
  \end{eqnarray}

\begin{proof}

To prove (\ref{86-1}), we first make some preparations.

For every $i\in \cZ_0$ and $n\in \mN,$ set
  $$G_n^i=\big\{y=(y_j,j\in \cZ_0)\in\mathbb{R}^d\setminus\{0\}:\ { \max_{j\in \cZ_0 \text{ with }j\neq  i}|y_j|< \frac{|y_i|}{n}}\big\}.$$
Then, for any $i\in \cZ_0$ and $u>0$, we have
  \begin{eqnarray*}
  &&\int_{G_n^i }(2\pi u)^{-d/2} e^{-\frac{|y|^2}{2u}} \dif y=\int_{G_n^1 }(2\pi u)^{-d/2} e^{-\frac{|y|^2}{2u}} \dif y
  \\
  &=&C_d
  \int_{0}^\infty\dif z_1 \int_0^{z_1/n}\dif z_2\cdots \int_0^{z_1/n}\dif z_d
  ~(2\pi u)^{-d/2} e^{-\frac{\sum_{k=1}^dz_k^2}{2u}}\\
  &\geq&
  C_d n^{-d+1}(2\pi )^{-d/2}\int_{0}^\infty
  ~u^{-d/2} e^{-\frac{dz_1^2}{2u}} z_1^{d-1}  \dif z_1\\
  &=&C_d n^{-d+1}(2\pi )^{-d/2}\int_{0}^\infty
  ~e^{-\frac{dx^2}{2}} x^{d-1}  \dif x\\
  &=&C_{d,n}>0.
  \end{eqnarray*}
   Hence, for every $i\in \cZ_0$ and $n\in \mN,$
 \begin{eqnarray*}
  \nu_L(G_n^i)&=& \int_{G_n^i }\Big(\int_0^\infty
  (2\pi u)^{-d/2} e^{-\frac{|y|^2}{2u}} \nu_{S}(\dif u)\Big)\dif y
  \\ &=& \int_0^\infty \nu_{S}(\dif u) \int_{G_n^i }(2\pi u)^{-d/2} e^{-\frac{|y|^2}{2u}} \dif y
  \\ &\geq & C_{d,n}   \nu_{S}((0,\infty))=\infty,
  \end{eqnarray*}
  where $\nu_L$ is the L\'evy measure of $L(t)=W_{S_t}$ given in  (\ref{1129}).
Therefore, there exists a $\Omega_0^1 \in\mathcal{F}$ with $\mP(\Omega_0^1)=1$ and on the set $\Omega_0^1$, Lemma \ref{31-1} applies to the L\'evy process $\tilde{L}(t)=L(t)=W_{S_t}.$

By Lemma  \ref{25-1},  there exists a set $\mS_0\subseteq \mS$ with $\mP^{\mu_\mS}(\mS_0)=1$ and for any $\ell\in \mS_0,$
$\{s:\Delta S_s(\ell):=\ell_s-\ell_{s-}>0\}   \text{ is dense in } (0,\infty)$, which implies that if $f$ is a nonnegative continuous function on
some time interval $[a,b]$  and $\int_a^bf(s)\dif\ell_s=0$, then $f(s)=0, s\in[a,b]$. Denote $\Omega_0^2=\mW\times \mS_0\subseteq \Omega. $ Obviously, $\mP(\Omega_0^2 )=1.$

We are now in the position to prove (\ref{86-1}). Set
\begin{eqnarray}
\label{11-1}
  \mathcal{L}:=\Big\{\omega:  \inf_{\phi\in \cS_{\alpha,N} }\sum_{j\in \cZ_0}b_j^2\int_{\sigma/2}^\sigma \langle K_{r,\sigma}\phi,e_j\rangle^2\dif \ell_r =0  \Big\}
  \cap \Omega_0^1\cap \Omega_0^2.
\end{eqnarray}
In the following, we will prove $\mathcal{L}=\varnothing$, completing the proof of  (\ref{86-1}).

Assume that $\mathcal{L}\neq\varnothing$ and $\omega=(\text{w},\ell)$ belongs to the event $\mathcal{L}$.
Then, for some $\phi$ with
\begin{eqnarray}
\label{1-3}
\| P_N\phi\|\geq \alpha,
\end{eqnarray}one has
$\int_{\sigma/2}^\sigma \langle K_{r,\sigma}\phi,e_j\rangle^2\dif \ell_r =0$ for $j\in \cZ_0$.
By the property of  $\ell\in \mS_0$ stated above and   the  continuity of   $\langle K_{r,\sigma}\phi,e_j\rangle$  with respect to $r$,
it holds that
\begin{eqnarray}
\label{7-1}
  \sup_{r\in [\sigma/2,\sigma]}\left|\langle K_{r,\sigma}\phi,e_j\rangle \right|=0, \quad \forall j\in \cZ_0.
\end{eqnarray}
With the help of  (\ref{w-1}),  $\varrho_t:=\langle K_{t,\sigma}\phi,e_j\rangle$
  satisfies the following equation:
  \begin{eqnarray*}
	&& \partial_t  \varrho _t+c_j \varrho_t +\langle K_{t,\sigma}\phi,B(\cK w_t,e_j)+B(\cK e_j,w_t)\rangle =0,
  \\ && \varrho_\sigma=\langle\phi,e_j\rangle, \nonumber
\end{eqnarray*}
where $c_j$ is a constant depending on $j.$
Combining the above equation with (\ref{7-1}),   we deduce that for any $t\in [\sigma/2,\sigma],$
\begin{eqnarray*}
  \langle K_{t,\sigma}\phi,B(\cK w_t,e_j)+B(\cK e_j,w_t)\rangle =0.
\end{eqnarray*}
Let $v_t=w_t-\sum_{i\in \cZ_0} b_i { W_{S_t}^i}  e_i.$  Then, the above equation becomes
\begin{eqnarray*}
  \langle K_{t,\sigma}\phi,\tilde B(v_t+\sum_{i\in \cZ_0}b_i  { W_{S_t}^i}  e_i,e_j)\rangle =0.
\end{eqnarray*}
That is
\begin{eqnarray*}
  f(t)+\sum_{i\in \cZ_0} { W_{S_t}^i}  b_i \langle K_{t,\sigma}\phi, \tilde B(e_i,e_j)\rangle = 0, ~\forall t\in[\sigma/2,\sigma],
\end{eqnarray*}
where  $f(t):= \langle K_{t,\sigma}\phi,\tilde B(v_t,e_j)\rangle$ is a continuous stochastic process.
By the assumption (\ref{11-1}), the above equality and the fact that
 Lemma \ref{31-1} holds for   $\omega\in \Omega_0^1$,
%
one arrives at that
\begin{eqnarray*}
 \langle K_{t,\sigma}\phi, \tilde B(e_i,e_j)\rangle = 0, \quad \forall t\in [\sigma/2,\sigma].
\end{eqnarray*}
Checking through the above   arguments,  we actually proved that
\begin{eqnarray}
\label{1657}
\begin{split}
 & \langle K_{t,\sigma}\phi, e_j \rangle = 0, ~\forall t\in [\sigma/2,\sigma]  \\ &  \Rightarrow
  \langle K_{t,\sigma}\phi, \tilde B(e_i,e_j)\rangle= 0, \forall i\in \cZ_0\text{ and }t\in [\sigma/2,\sigma].
  \end{split}
\end{eqnarray}

Define the set $\cZ_n\subseteq \mZ_*^2$ recursively:
\begin{eqnarray*}
  \cZ_n=\{i+j \big| j\in \cZ_0,i \in \cZ_{n-1} \text{ with } \langle i^{\perp},j\rangle\neq 0, |i|\neq |j| \},
\end{eqnarray*}
where $i^{\perp}=(i_2,-i_1).$
Assume that we have proved
\begin{eqnarray*}
 \langle K_{t,\sigma}\phi, e_j \rangle = 0, \forall j\in \cZ_{n-1} \text{ and }t\in [\sigma/2,\sigma].
\end{eqnarray*}
Then, by (\ref{1657}), it follows that
\begin{eqnarray*}
 \langle K_{t,\sigma}\phi, \tilde B(e_j,e_i) \rangle= 0, \forall j\in \cZ_{n-1}, i\in \cZ_0 \text{ and }t\in [\sigma/2,\sigma].
\end{eqnarray*}
It is easy to verify that $\cZ_m$ is symmetric for any $m\geq 0,$ i.e.
$\cZ_m=-\cZ_m.$
Also by the definition of $\cZ_n,$  one can see  that
\begin{eqnarray*}
\{ e_j, j\in \cZ_n \} \subseteq  \text{span}  \{ \tilde B(e_i,e_j):  j\in \cZ_0, i\in \cZ_{n-1}\}.
\end{eqnarray*}
Hence,
\begin{eqnarray*}
 \langle K_{t,\sigma}\phi, e_j \rangle = 0, \forall j\in \cZ_{n} \text{ and }t\in [\sigma/2,\sigma].
\end{eqnarray*}
By this  recursion,
\begin{eqnarray*}
 \langle K_{t,\sigma}\phi, e_j \rangle = 0, \quad \forall j\in \cup_{n=1}^\infty \cZ_{n}=\mZ_*^2 \text{ and }t\in [\sigma/2,\sigma].
\end{eqnarray*}
(Here, we have used  \cite[Proposition 4.4]{HM-2006}.)
Let $t\rightarrow \sigma$ to get  $\phi=0$,  which contradicts (\ref{1-3}) .
Therefore, $\mathcal{L}=\varnothing$.

The proof of (\ref{86-1}) is complete.


\end{proof}


Proposition \ref{1-66}  is not sufficient for the proof of Proposition \ref{3-11}, we need a stronger statement.
 For $\alpha\in (0,1], w_0\in H,N\in \mN,\mathfrak{R}>0$ and $\eps>0$, let\footnote{Note that $\cM_{0,t}$ is the Malliavin matrix of $w_t$, the solution of equation (\ref{0.1})  at time $t$ with  initial value $w_0.$
 Therefore, $\cM_{0,\sigma}$ also depends on $w_0$. }
\begin{eqnarray}
\label{1-5}
  X^{w_0,\alpha,N}=\inf_{ \phi \in \cS_{\alpha,N}}
  \langle \cM_{0,\sigma}\phi,\phi\rangle.
\end{eqnarray}
and denote
\begin{eqnarray}
\label{35-1}
  r(\eps,\alpha, \mathfrak{R},N)=\sup_{{\|w_0\|\leq  \mathfrak{R}}}\mP(  X^{w_0,\alpha,N}<\eps).
\end{eqnarray}

Based on (\ref{86-1}) and the dissipative property of Navier-Stokes system,  we obtain
the following  result  whose proof is given  in  Appendix \ref{180-1}.
\begin{proposition}
\label{3-8}
For  $\alpha \in (0,1], \mathfrak{R}>0$ and $N\in \mN,$
we have
\begin{eqnarray*}
  \lim_{\eps\rightarrow 0}r(\eps,\alpha,\mathfrak{R},N)=0.
\end{eqnarray*}
\end{proposition}

\section{Proof of  Proposition  \ref{3-11}. }
\label{54-1}
Let us take $f\in C_b^1(H)$ and $\xi\in H$ with $\|\xi\|=1.$  Compute the derivative of $ \mE_{w_0} f(w_{t})$ with respect to $w_0$ in the direction $\xi$:
\begin{eqnarray}
\label{3.2}
  && \nabla_\xi  \mE_{w_0} f(w_{t})= \mE \nabla
  f(w_{t})J_{0,t }\xi.
\end{eqnarray}
In the papers~\cite{HM-2006, HM-2011},
their  ideas of proof of  the  asymptotic strong Feller property is~to approximate the perturbation
  $J_{0,t}\xi$ caused by the variation of the initial condition with a variation, $\cA_{0,t}v=\cD^v w_t $, of the noise by an appropriate   process $v$. Denote by $\rho_t$ the  residual error between~$J_{0,t}\xi$    and~$\cA_{0,t}v$:
\begin{align*}
  \rho_t=J_{0,t}\xi-\cA_{0,t}v.
\end{align*}
The proof of Proposition \ref{3-11}  is much  more involved than that in ~\cite{HM-2006, HM-2011} since  we even  don't have $ \mE \|J_{0,t}\xi\|<\infty.$

Let us first explain the main ideas of  the proof of  Proposition  \ref{3-11}.
Let
$\kappa_0=\kappa_0(\nu,\{b_j\}_{j\in \cZ_0},\nu_S,d)$
be the constant appeared in the statement of
Lemma \ref{qu-2}.
Recall that, for any $\kappa \in (0,\kappa_0],$
the stopping times $\sigma_k$ are defined in (\ref{2121-4})--(\ref{qu-20}).
For any $\kappa\in (0,\kappa_0]$ and $n\in \mN,$ { we define the following random variables on $\mS:$}
\begin{eqnarray}
\label{58-1}
  X_n=\int_{\sigma_{n}}^{\sigma_{n+1}} e^{2\nu(\sigma_{n+1}-s)-16\mathfrak{B}_0\kappa (\ell_{\sigma_{n+1}}-\ell_s)}\dif s, \quad Y_n=\ell_{\sigma_{n+1}}-\ell_{\sigma_n}.
\end{eqnarray}
By the strong law of large numbers, Lemma \ref{qu-2}  and the definitions of $\sigma_k$, we have\footnote{Actually, by Lemma \ref{qu-2},  $\mE X_n^2\leq \mE \[(\sigma_{n+1}-\sigma_n)^2 e^{4\nu (\sigma_{n+1}-\sigma_n)}\]<\infty.  $}
\begin{eqnarray*}
  \lim_{n\rightarrow \infty}\frac{\sum_{i=0}^{n-1}X_i}{n}<\infty, \quad \text{almost surely},
\end{eqnarray*}
and
\begin{eqnarray*}
  \lim_{n\rightarrow \infty}\frac{\sum_{i=0}^{n-1}Y_i}{n}
  \leq    \frac{\nu }{8 \mathfrak{B}_0 \kappa }\lim_{n\rightarrow \infty}\frac{\sum_{i=0}^{n-1}(\sigma_{i+1}-\sigma_i)}{n}<\infty,\quad \text{almost surely.}
\end{eqnarray*}
Therefore, with probability one, we have
\begin{eqnarray}
\label{1609}
  \Theta:=\sup_{n\geq 1}\frac{\sum_{i=0}^{n-1}X_i}{n}+\sup_{n\geq 1}
  \frac{\sum_{i=0}^{n-1}Y_i}{n}<\infty.
\end{eqnarray}
For any $\Upsilon,M>0, f\in C_b^1(H)$ and  $w_0,w_0'\in B_H(\Upsilon):=\{w\in H, \|w\|<\Upsilon\}$,  one has
\begin{eqnarray}
 \nonumber  && | P_t f(w_0)-P_t f(w_0')|
  = | \mE f(w_t^{w_0})
  -\mE f(w_t^{w_0'})|
  \\  \nonumber && \leq | \mE f(w_t^{w_0}) I_{\{\Theta \leq M\}}
  -\mE f(w_t^{w_0'})I_{\{\Theta \leq M\}}|+ 2\|f\|_\infty \mP(\Theta \geq M )
  \\ \label{20-3} &&:=I_1+I_2,
\end{eqnarray}
where $w_t^{w_0}$ is the solution of equation (\ref{0.1}) with initial value $w_0.$


One can choose the constant $M$ sufficiently large, independent of the initial condition $w_0$ and time $t$, to make $I_2$ arbitrarily small.
The main difficulty lies in the    estimate of   $I_1.$
 Denote  $P_t^{M} f(w_0)=\mE \[f(w_t^{w_0})I_{\{\Theta\leq M\}}\].$
 For  any process   $v \in L^2(\Omega \times [0,\infty); \mR^d)
 =L^2(\mW \times \mS \times [0,\infty); \mR^d)
 $, we write
  \begin{eqnarray}
   \nonumber  && | \nabla_\xi P_t^{M} f(w_0) |=
    | \mE \nabla f(w_t)J_{0,t}\xi
    I_{\{\Theta\leq M\}} |
    \\ \label{1604-1}  &&
    =| \mE \[ \nabla f(w_t)\cD^v w_t
    I_{\{\Theta\leq M\}}\]+\mE \[ \nabla  f(w_t)\rho_t
    I_{\{\Theta\leq M\}}\] |.
    \end{eqnarray}
For any fixed $\ell\in \mS,$ the process $v=v^\ell$  in the above  will be chosen such that $v^\ell\in L^2(\mW \times [0,\ell_t]; \mR^d) $ and that $v^\ell$ is  Skorokhod integrable with respect to the Brownian motion $W$.
    Since  $\{\ell: \Theta(\ell)\leq M\}$ is independent of the Brownian motion $W_t$, it holds that
\begin{eqnarray}
  \nonumber && \mE \[ \nabla f(w_t)\cD^v w_t
    I_{\{\Theta\leq M\}}\]= \mE^{\mu_{\mS}}\big[ I_{\{\Theta\leq M\}} \mE^{\mu_{\mW}}\big( \nabla f(w_t)\cD^v w_t\big)
    \big]
    \\  \nonumber &&= \mE^{\mu_{\mS}}\big[ I_{\{\Theta\leq M\}} \mE^{\mu_{\mW}}\big(  f(w_t)\int_0^{\ell_t } v(s)\dif W(s) \big)
    \big]
    \\ \label{1604-2}  && =\mE \[  f(w_t)\int_0^{\ell_t} v(s)\dif W(s)
    I_{\{\Theta\leq M\}}\],
\end{eqnarray}
In the above, for any fixed $\ell\in \mS,$  the integral $\int_0^{\ell_t} v(s)\dif W(s)$ is interpreted as the Skohorod integral.
{ In a word,} the estimate of $I_1$ is obtained through some   gradient  estimates of
$\nabla_\xi P_t^{M} f(w_0).$
In order to do this, by (\ref{1604-1})--(\ref{1604-2}),
we need to select suitable direction $v$ and
do some moment estimates for   $\rho_t$   and  $\int_0^{\ell_t} v(s)\dif W(s).$
This will be done in subsections  \ref{40-1}--\ref{40-3}. In subsection \ref{17-27}, we complete the  proof of Proposition \ref{3-11}.



\subsection{The choice of $v.$}
\label{40-1}

%

In this section, we always assume that $\|\xi\|=1$.
For any $\kappa>0,$  recall that the stopping times $\sigma_k$ are defined in (\ref{2121-4})--(\ref{qu-20}).
For any $\ell\in \mS$ and $\kappa>0$,
we will define the  perturbation $v$ to be $0$
 on all intervals of the type $[\ell_{\sigma_{n+1}},\ell_{\sigma_{n+2}}],
 n\in 2\mN,$ and by some  $v_{{\sigma_n},{\sigma_{n+1}}} \in L^2([\ell_{\sigma_n},\ell_{\sigma_{n+1}}],H), n\in 2\mN$
on the remaining intervals.
For fixed $\ell \in \mS$ and $n\in 2\mN,$ define
the infinitesimal variation:
\begin{eqnarray}
\label{p-1}
\begin{split}
  & v_{{\sigma_n},{\sigma_{n+1}}}(r)=\cA^*_{\sigma_n,\sigma_{n+1}}
  {(\cM_{\sigma_n,\sigma_{n+1}}+\beta\mI)}^{-1}J_{\sigma_n,\sigma_{n+1}}\rho_{\sigma_{n}},
  ~ r\in [\ell_{\sigma_n},\ell_{\sigma_{n+1}}]
 \\
 & v_{{\sigma_{n+1}},{\sigma_{n+2}}  }(r)=0,~ r\in [\ell_{\sigma_{n+1}},\ell_{\sigma_{n+2}}].
  \end{split}
\end{eqnarray}
where $\rho_{\sigma_n} $ is the residual of the infinitesimal
 displacement at time $\sigma_n$. Set
 \begin{align}\label{000}
 v(r)=
\begin{cases}
 v_{{\sigma_n},{\sigma_{n+1}}}(r),  \quad r\in [\ell_{\sigma_n},\ell_{\sigma_{n+1}}] { \text{ and } n\in 2\mN}, \\
v_{{\sigma_{n+1}},{\sigma_{n+2}}  }(r), \quad\quad r\in [\ell_{\sigma_{n+1}},\ell_{\sigma_{n+2}}]  { \text{ and }  n\in 2\mN}.
\end{cases}
\end{align}

  Here and after, we use  $v_{a,b}$  to denote the function $v$ restricted on the interval $[\ell_a,\ell_b]$  and
  the constant $\beta$  in  (\ref{p-1})
  will be decided later.
  Obviously,    $\rho_0=J_{0,0}\xi-\cA_{0,0}v=\xi.$

Similar to \cite{HM-2006}, we have the following  recursions  for $\rho_{\sigma_n}.$

  \begin{lemma}
\label{28-5}
For any $\beta>0,$
  if we definethe direction $v$  according to  (\ref{000}),
    then  
\begin{eqnarray*}
\rho_{\sigma_{n+2} }&=&{J}_{\sigma_{n+1},\sigma_{n+2}}\beta
(\cM_{\sigma_n,\sigma_{n+1}}+\beta \mI)^{-1}{J}_{\sigma_n,\sigma_{n+1}}
\rho_{\sigma_{n}}, \quad \forall n\in 2\mN.
\end{eqnarray*}

\end{lemma}
\begin{proof}
By a straightforward  calculation,
\begin{eqnarray*}
 && \rho_{\sigma_{n+2} }= {J}_{0,\sigma_{n+2}}\xi-\cA_{0,\sigma_{n+2}}v
= {J}_{0,\sigma_{n+2}}\xi-
 \int_{0}^{\ell_{\sigma_{n+2}}  } J_{\gamma_u,\sigma_{n+2} } Q v_u \dif u
 \\ && = J_{\sigma_{n+1},\sigma_{n+2} }
 {J}_{0,\sigma_{n+1} }\xi-J_{\sigma_{n+1},
 \sigma_{n+2}} \int_{0}^{\ell_{\sigma_{n+1}}  }J_{\gamma_u,\sigma_{n+1} } Q v_u
 \dif u
 \\ & &= J_{\sigma_{n+1},\sigma_{n+2}} \rho_{\sigma_{n+1}}
\end{eqnarray*}
and
\begin{eqnarray*}
 \rho_{\sigma_{n+1} } &=& {J}_{0,\sigma_{n+1} }\xi- \int_{0}^{\ell_{\sigma_{n+1}} }
 J_{\gamma_u,\sigma_{n+1} } Q  v_{u} \dif u
 \\ &=& {J}_{0,\sigma_{n+1} }\xi-\int_{0}^{\ell_{\sigma_{n}}}
 J_{\gamma_u,\sigma_{n+1} }
 Q v_{u}  \dif u  -
 \int_{\ell_{\sigma_{n}}}^{\ell_{\sigma_{n+1}} }J_{\gamma_u,\sigma_{n+1}} Q v_{u}  \dif u
 \\ &=& J_{\sigma_{n},\sigma_{n+1}}\rho_{\sigma_{n}}  -\cA_{\sigma_{n},\sigma_{n+1}}
 \cA^*_{\sigma_{n},\sigma_{n+1}}{(\cM_{\sigma_{n},\sigma_{n+1}}+\beta\mI)}^{-1}
 J_{\sigma_{n},\sigma_{n+1}}\rho_{\sigma_{n}}
 \\ &=& \beta (\cM_{\sigma_{n},\sigma_{n+1}}+\beta\mI)^{-1}{J}_{\sigma_{n},\sigma_{n+1}}
 \rho_{\sigma_{n}},
\end{eqnarray*}
which yields the desired result.
\end{proof}

\subsection{The control of $\rho_{\sigma_n}.$}
\label{40-2}

Let  $A_\eps=A_{\eps,w_0,\alpha,N}=\{X^{w_0,\alpha,N}\geq \eps\},$
where  the random variable $X^{w_0,\alpha,N}$ is defined in (\ref{1-5}).
To provide an estimate for
$\|\rho_{\sigma_{n}}\|,$ we start with some preparations.
\begin{lemma}(c.f. \cite[Lemma 5.4]{HM-2011})
\label{3-3}
  For any positive constants  $\beta,\eps,\alpha\in(0,1],N\in \mN$
  and $\xi\in H,$  the following inequality holds with probability $1$:
  \begin{eqnarray}
  \label{3-1}
 \beta \|P_N(\beta\mI+\cM_{0,\sigma})^{-1}\xi\|
\leq   \| \xi\|
\big(\alpha \vee \sqrt{\frac{\beta}{\eps}}\big) I_{A_\eps} +\|\xi\|I_{A_\eps^c }.
\end{eqnarray}
\end{lemma}

\begin{proof}
On the event  $A_\eps^c$, the inequality (\ref{3-1})  obviously  holds.
On the event $A_\eps$, this inequality is proved in \cite[Lemma 5.14]{HM-2011}, so we omit the details.

%
%
%
%
%
%
\end{proof}

Let $\cR_{\sigma_n,\sigma_{n+1}}^\beta=\beta(\cM_{\sigma_n,\sigma_{n+1}}+\beta\mI)^{-1}.$
We have  the following estimate for $\cR_{\sigma_n,\sigma_{n+1}}^\beta.$
\begin{lemma}
\label{3-9}
    For any $\kappa>0,\delta\in (0,1],p\geq 1$ and $N\in \mN$, there exists a $\beta=\beta(\kappa,\delta,p,N)>0$
  such  that
    \begin{eqnarray}
    \label{3-4}
\mE\[\|P_N \cR_{\sigma_{n},\sigma_{n+1}}^\beta\|^p \big| \cF_{\sigma_n} \]
\leq \delta e^{\kappa \|w_{\sigma_{n}}\|^2},\quad \forall n\in \mN.
\end{eqnarray}
\end{lemma}
\begin{proof}
  Denote $\sigma_1$ by $\sigma.$ We here give a proof for the case $n=0$ and $p\geq 2.$ The other cases can be proved similarly.
   Let  $\mathfrak{R}=\mathfrak{R}_{\delta,\kappa}$ be a positive  constant  such that
 $
  \exp\{\kappa \mathfrak{R}^2  \}\geq \frac{1}{\delta}.
$
  We divide into the following two cases to prove (\ref{3-4}).

 \textbf{Case 1:} $\|w_0\|\geq \mathfrak{R}$.
  In this case,
  \begin{eqnarray*}
\mE\[\|P_N \cR_{0,\sigma}^\beta \|^p \big| \cF_0 \]\leq
1\leq \delta e^{\kappa \|w_{0}\|^2}.
  \end{eqnarray*}

 \textbf{Case 2:}
 $\|w_0\| \leq \mathfrak{R}$.
 For any  positive constants $\eps,\beta$ and $\alpha\in (0,1],$ by Lemma \ref{3-3}, we have
  \begin{eqnarray*}
 \mE\[\|P_N \cR_{0,\sigma}^\beta \|^p \big| \cF_0 \] & \leq&
 C_p \big(\alpha \vee \sqrt{\frac{\beta}{\eps}}\big)^p+
 C_p r(\eps,\alpha, \mathfrak{R},N),
  \end{eqnarray*}
  where $C_p$ is a constant only depending on $p,$
  and  $r(\eps,\alpha, \mathfrak{R},N)$ is defined in (\ref{35-1}).
 Choose now  $\alpha=\alpha({p})$ small enough such that
 \begin{eqnarray*}
 C_p \alpha^p \leq \frac{\delta}{2}.
 \end{eqnarray*}
 By Proposition \ref{3-8}, $\lim_{\eps \rightarrow 0}  r(\eps,\alpha, \mathfrak{R},N)=0$.
Pick a small constant  $\eps$ such that
 \begin{eqnarray*}
C_p r(\eps,\alpha, \mathfrak{R},N)\leq \frac{\delta}{2}.
 \end{eqnarray*}
 Finally, we choose  $\beta$ small enough so that
 \begin{align*}
  C_p( \sqrt{\beta/\eps})^p <\frac{\delta}{2}.
 \end{align*}
Putting the above steps together, we see that
$
  \mE\[\|P_N \cR_{0,\sigma}^\beta \|^p \big|\cF_0  \]  \leq \delta
e^{\kappa \|w_{0}\|^2} .
$
\end{proof}

By Lemma \ref{28-5}, 
for any   $n\in 2\mN, \beta>0$ and $N\in \mN,$   one easily sees that
\begin{eqnarray}
 \nonumber  && \rho_{\sigma_{n+2}}
 = J_{\sigma_{n+1},\sigma_{n+2}}
 \cR_{\sigma_n,\sigma_{n+1}}^\beta J_{\sigma_n,\sigma_{n+1}}\rho_{\sigma_n}
  \\  \nonumber && = J_{\sigma_{n+1},\sigma_{n+2}}Q_N\cR_{\sigma_n,\sigma_{n+1}}^\beta J_{\sigma_n,\sigma_{n+1}}\rho_{\sigma_n}
  +J_{\sigma_{n+1},\sigma_{n+2}}P_N\cR_{\sigma_n,\sigma_{n+1}}^\beta J_{\sigma_n,\sigma_{n+1}}\rho_{\sigma_n}
  \\  \label{4-1} && :=  \rho_{\sigma_{n+2}}^{(1)}+ \rho_{\sigma_{n+2}}^{(2)}.
\end{eqnarray}
The values  of $\beta$  and $N$  will be decided later.


To estimate $\|\rho_{\sigma_{n+2}}\|$, we first consider  the term  $ \rho_{\sigma_{n+2}}^{(1)}.$
For any $\kappa>0,\xi\in H$ and $n\in \mN$,  by  Lemma  \ref{1340-2} and Young's inequality, we have
       \begin{eqnarray}
       \nonumber
  && \|J_{\sigma_{n+1},\sigma_{n+2}}Q_N \xi\|^2
 \\ \nonumber && \leq C \big(e^{-\nu N^2 (\sigma_{n+2}-\sigma_{n+1})}+\frac{1+\sigma_{n+2}-\sigma_{n+1}}{N^{1/4}} \big)\sup_{r\in [\sigma_{n+1},\sigma_{n+2}]}(\|w_r\|^{5/2}+1)\|\xi\|^2
 \\ \nonumber &&
  ~\times \exp\Big\{\frac{\nu \kappa}{120} \int_{\sigma_{n+1}}^{\sigma_{n+2}} \|w_r\|_1^{2}e^{-\nu(\sigma_{n+2}-r)+8\mathfrak{B}_0 \kappa (\ell_{\sigma_{n+2}}-\ell_r) }\dif r
\\ \label{1340-4} &&  \quad\quad \quad\quad  +C_\kappa \int_{\sigma_{n+1}}^{\sigma_{n+2}} e^{2 \nu(\sigma_{n+2}-r)-16 \mathfrak{B}_0  \kappa (\ell_{\sigma_{n+2}}-\ell_r) }\dif r  \Big\} .
  \end{eqnarray}
here and below    $C$ is a constant depending on  $ \nu,\{b_j\}_{j\in \cZ_0},\nu_S,d$, ~$C_\kappa$ is a constant depending on  $\kappa,\nu,\{b_j\}_{j\in \cZ_0},\nu_S,d.$
 Applying
(\ref{1340-1}) with $s=\sigma_n$ and $T=\sigma_{n+1}$, also  with the help of   (\ref{1340-4}) and the expression of  $\rho_{\sigma_{n+2}}^{(1)}$, we have
\begin{eqnarray}
 \nonumber  && \| \rho_{\sigma_{n+2}}^{(1)}\|^{40}
  \leq  \| J_{\sigma_{n+1},\sigma_{n+2}}Q_N\|^{40} \| J_{\sigma_n,\sigma_{n+1}}\|^{40}\|\rho_{\sigma_n}\|^{40}
 \\  \nonumber && \leq
 C  \exp\Big\{\frac{\nu \kappa }{6} \int_{\sigma_{n}}^{\sigma_{n+1} } \|w_s\|_1^{2}e^{-\nu(\sigma_{n+1}-s)+8\mathfrak{B}_0\kappa (\ell_{\sigma_{n+1}}-\ell_s)}\dif s-\frac{\kappa}{6}\mathfrak{B}_0(\ell_{\sigma_{n+1}}-\ell_{\sigma_{n}}) \Big\}
 \\  \nonumber  && \times
 \exp\Big\{\frac{\nu \kappa }{6} \int_{\sigma_{n+1}}^{\sigma_{n+2} } \|w_s\|_1^{2}e^{-\nu(\sigma_{n+2}-s)+8\mathfrak{B}_0\kappa (\ell_{\sigma_{n+2}}-\ell_s)}\dif s-\frac{\kappa}{6}\mathfrak{B}_0(\ell_{\sigma_{n+2}}-\ell_{\sigma_{n+1}}) \Big\}
\\  \nonumber  && \times
 \sup_{s\in [\sigma_{n+1},\sigma_{n+2} ]}(\|w_s\|^{50}+1)
 \\   \nonumber && \times
 \Big(\exp\big\{-20\nu N^2 (\sigma_{n+2}-\sigma_{n+1})\big\} +
 \frac{(1+\sigma_{n+2}-\sigma_{n+1})^{20}   }{N^{5}}\Big)
  \\   \nonumber && \times \exp\big\{\frac{\kappa}{6}\mathfrak{B}_0(\ell_{\sigma_{n+1}}-\ell_{\sigma_{n}})  +\frac{\kappa}{6}\mathfrak{B}_0(\ell_{\sigma_{n+2}}-\ell_{\sigma_{n+1}})  \big\}
 \\  \nonumber && \times  \exp\Big\{C_\kappa \int_{\sigma_{n}}^{\sigma_{n+1}} e^{2\nu(\sigma_{n+1}-s)-16\mathfrak{B}_0\kappa (\ell_{\sigma_{n+1}}-\ell_s)}\dif s  +C_\kappa \int_{\sigma_{n+1}}^{\sigma_{n+2}} e^{2\nu(\sigma_{n+2}-s)-16\mathfrak{B}_0\kappa (\ell_{\sigma_{n+2}}-\ell_s)}\dif s  \Big\}
 \\  \nonumber  && \times  \| \rho_{\sigma_{n}}  \|^{40}
 \\  \label{36-1}  &&:= C U_nU_{n+1} V_nR_ne^{C_\kappa(Y_n+Y_{n+1})} e^{C_\kappa(X_n+X_{n+1})} \| \rho_{\sigma_{n}}  \|^{40},
\end{eqnarray}
where  $X_n,Y_n$ are defined  in (\ref{58-1}),  and
\begin{eqnarray*}
  && U_n=\exp\Big\{\frac{\nu \kappa }{6} \int_{\sigma_{n}}^{\sigma_{n+1} } \|w_s\|_1^{2}e^{-\nu(\sigma_{n+1}-s)+8\mathfrak{B}_0\kappa (\ell_{\sigma_{n+1}}-\ell_s)}\dif s-\frac{\kappa}{6}\mathfrak{B}_0(\ell_{\sigma_{n+1}}-\ell_{\sigma_{n}}) \Big\},
  \\ && V_n= \sup_{s\in [\sigma_{n+1},\sigma_{n+2} ]}(\|w_s\|^{50}+1),
  \\ && R_n=\exp\big\{-20\nu N^2 (\sigma_{n+2}-\sigma_{n+1})\big\} +
 \frac{(1+\sigma_{n+2}-\sigma_{n+1})^{20}   }{N^{5}}.
\end{eqnarray*}
For the second  term  $ \rho_{\sigma_{n+2}}^{(2)},$
using
(\ref{1340-1}) twice,
with the similar arguments as that for deducing  (\ref{36-1}),  we obtain
\begin{eqnarray}
\nonumber &&  \|\rho_{\sigma_{n+2}}^{(2)}\|^{40}
\\ \label{36-2}  &&\!  \leq C U_n U_{n+1}\|P_N\cR_{\sigma_n,\sigma_{n+1}}^\beta\|^{40} e^{C_\kappa(Y_n+Y_{n+1})+C_\kappa(X_n+X_{n+1}) }\| \rho_{\sigma_{n}}  \|^{40}.
\end{eqnarray}

Combining  (\ref{4-1}) with  (\ref{36-1})(\ref{36-2}), for any $n\in 2\mN,$ one arrives at
\begin{eqnarray*}
\| \rho_{\sigma_{n+2}}\|^{40}\leq C \theta_n e^{C_\kappa(Y_n+Y_{n+1})+C_\kappa(X_n+X_{n+1}) }\| \rho_{\sigma_{n}}\|^{40},
\end{eqnarray*}
where 
$\theta_n=U_nU_{n+1} V_nR_n+U_n U_{n+1}\|P_N\cR_{\sigma_n,\sigma_{n+1}}^\beta\|^{40}.$
Assume that $\|\xi\|=1,$
then $\|\rho_0\|=1.$
By recursion, for any $n\in \mN,$  we have
\begin{eqnarray}
\label{30-2}
 \| \rho_{\sigma_{2n+2}}\|^{40}
  \leq C^{ n+1} \big(\prod_{i= 0}^{n} \theta_{2i} \big)  e^{C_\kappa \sum_{i=0}^{2n+1}(X_i+Y_i)},
\end{eqnarray}
Recall that   $\Theta =\sup_{n\geq 1}\frac{\sum_{i=0}^{n-1}X_i}{n}+\sup_{n\geq 1}\frac{\sum_{i=0}^{n-1}Y_i}{n}.$
Thus, on the event
$$
\{\ell:\Theta(\ell) \leq M\},
$$
it holds that
\begin{eqnarray}
\label{30-5}
 \| \rho_{\sigma_{2n+2}}\|^{40}\leq C^{n+1 }  e^{C_\kappa M(2n+2)} \prod_{i=0}^{n}\theta_{2i}, \quad \forall n\in   \mN.
\end{eqnarray}
  Notice  that   $\theta_i$  depends on the parameters $N, \beta.$  We have the following estimates for $\theta_i, i\in \mN$.
\begin{lemma}
\label{30-6}
   For  any   $\kappa \in (0,\kappa_0]$   and  $\delta\in (0,1)$,
 there exist       constants   $\beta>0 $ and  $N\in \mN$  which  depend   on   $\kappa,\delta, \nu,\{b_j\}_{j\in \cZ_0},\nu_S,d$
 such that
  \begin{eqnarray*}
    \mE \[\theta_i   \big| \cF_{\sigma_{i}}\]\leq \delta \exp\{\kappa \|w_{\sigma_i}\|^2 \}, \quad \forall i\in \mN.
  \end{eqnarray*}
\end{lemma}
\begin{proof}
  Obviously, for any  $\kappa>0,\delta'>0,$ by Lemma \ref{qu-2},   there exists a $N=N(\kappa,\delta',\nu,\{b_j\}_{j\in \cZ_0},\nu_S,d)>0$ such that
  \begin{eqnarray*}
\big(\mE |  R_i|^3  \big| \cF_{\sigma_{i}}\big)^{1/3} \leq \frac{\delta'}{4},\quad \forall i\in  \mN.
  \end{eqnarray*}
Hence, also  with the help of   Lemma \ref{qu-2}, it holds that
   \begin{eqnarray}
  \nonumber
     && \mE \[U_iU_{i+1} V_iR_i  \big| \cF_{\sigma_{i}}\]
     \\   \nonumber &&  \leq  \big(\mE |U_i|^6 \big| \cF_{\sigma_{i}} \big)^{1/6}
     \big(\mE |U_{i+1}|^6 \big| \cF_{\sigma_{i}} \big)^{1/6}
    \big(\mE |  V_i|^3\big| \cF_{\sigma_{i}} \big)^{1/3}
    \big(\mE |  R_i|^3   \big| \cF_{\sigma_{i}}\big)^{1/3}
    \\ \label{30-3} && \leq C_\kappa \frac{\delta'}{4} e^{\kappa \|w_{\sigma_i}\|^2 }.
  \end{eqnarray}

  By Lemma \ref{3-9}, for any  $\delta'>0,\kappa>0$  and  the $N=N(\kappa,\delta',\nu,\{b_j\}_{j\in \cZ_0},\nu_S,d)$ given above,
  there exists a $\beta=\beta(\kappa,\delta',N)>0$
  such  that
  \begin{eqnarray*}
    \big(\mE \|P_N\cR_{\sigma_i,\sigma_{i+1}}^\beta\|^{60}\big| \cF_{\sigma_{i}} \big)^{2/3}\leq \frac{\delta' }{4}\exp\{\frac{\kappa}{3}  \|w_{\sigma_i}\|^2\},\quad \forall i\in \mN.
  \end{eqnarray*}
Hence,
    \begin{eqnarray*}
     && \mE \big(U_iU_{i+1} \|P_N\cR_{\sigma_i,\sigma_{i+1}}^\beta\|^{40} \big| \cF_{\sigma_{i}}\big)
     \\ &&  \leq  \big(\mE |U_i|^6 \big| \cF_{\sigma_{i}} \big)^{1/6}
    \big(\mE |  U_{i+1}|^6   \big| \cF_{\sigma_{i}}\big)^{1/6}\big(\mE \|P_N\cR_{\sigma_i,\sigma_{i+1}}^\beta\|^{60}\big| \cF_{\sigma_{i}} \big)^{2/3}
    \\ && \leq C_\kappa \frac{\delta'}{4} e^{\kappa \|w_{\sigma_i}\|^2 }.
  \end{eqnarray*}
Combining this with (\ref{30-3}) and setting $\delta'=\frac{2\delta}{C_\kappa}$, we complete the proof.
\end{proof}

For  any  $\kappa \in (0,\kappa_0]$, $M>0 $ and $\gamma_0>0$, the constant $\beta$ in (\ref{p-1}) is decided through  the following Lemma.
 Recall that, $ \cC_0=\cC_0(\nu)$ is the constant introduced in Lemma  \ref{15-4}.

\begin{lemma}
\label{30-1}
     For  any  $\kappa \in (0,\kappa_0]$,  $M>0, \gamma_0>0$,
there exists a  positive    constant  $\beta=\beta(\kappa,M,\gamma_0,\nu,\{b_j\}_{j\in \cZ_0},\nu_S,d)$
   such that
  if we define the direction $v$  according to  (\ref{p-1}),
    then  
  the following holds
  \begin{eqnarray}
%
     \nonumber && \mE_{w_0} \[ (1+\sigma_{2n+2}-\sigma_{2n})^8 \exp\{8 \cC_0 \int_{\sigma_{2n}}^{\sigma_{2n+2}}\|w_s\|_1^{4/3}\dif s \} \|\rho_{\sigma_{2n}}\|^4  I_{\{\Theta \leq M\}}\]
    \\ \label{1002}  && \leq C_{\kappa,M,\gamma_0} \exp\{ 4  \kappa a \|w_0\|^2-n \gamma_0\}
  \end{eqnarray}
  for  every  $ n\in \mN$ and $w_0\in H$,
  where $a=\frac{1}{1-e^{-1}},$   $C_{\kappa,M,\gamma_0}$ is a constant depending on
  $\kappa,M,\gamma_0$ and  $\nu,\{b_j\}_{j\in \cZ_0},\nu_S,d.$

\end{lemma}
\begin{proof}
%
%
%
%
As in other places of this paper,
  $C$  denotes  a constant that  may  depend  on  $ \nu,\{b_j\}_{j\in \cZ_0},\nu_S,d;$      $C_\kappa$  denotes   a constant that  may   depend  on  $\kappa,\nu,\{b_j\}_{j\in \cZ_0},\nu_S,d;$
  First, we have the following

\emph{\textbf{Claim.}} \emph{For  any  $\delta\in (0,1)$  and  $\kappa \in (0,\kappa_0]$,
 there exists     constants   $\beta=\beta(\kappa,\delta,\nu,\{b_j\}_{j\in \cZ_0},\nu_S,d)$ and $N=N(\kappa,\delta,\nu,\{b_j\}_{j\in \cZ_0},\nu_S,d)$
 such that
\begin{eqnarray}
\label{1003}
\mE \prod_{i=0}^{n} \theta_{2i}^{1/2} \leq \delta^{(n+1)/2}   e^{2 a\kappa \|w_0\|^2+C_\kappa n}, \quad \forall n\in \mN,
  \end{eqnarray}
  where  $a=\frac{1}{1-e^{-1}},$ $\theta_{i}$ is defined as in Lemma \ref{30-6}.
}
Indeed,
by Lemma \ref{30-6} and (\ref{30-7}), we can choose  a  $\beta>0$ and  a $N\in \mN$ such that
\begin{eqnarray*}
 && \mE  \! \prod_{i=0}^{n} \theta_{2i}^{1/2}
\! =\mE  \Big[\!\! \prod_{i=0}^{n}\theta_{2i}^{1/2} e^{-\frac{\kappa}{2} \sum_{i=0}^{n} \|w_{\sigma_{2i}}\|^2  } e^{\frac{\kappa}{2} \sum_{i=0}^{n}  \|w_{\sigma_{2i}}\|^2  }\Big]
 \\ &&\leq
 \big( \mE  \prod_{i=0}^{n}\theta_{2i} e^{-\kappa \sum_{i=0}^{n}  \|w_{\sigma_{2i}}\|^2  }
 \big)^{1/2}
 \big( \mE  e^{\kappa  \sum_{i=0}^{n} \|w_{\sigma_{2i}}\|^2  }\big)^{1/2}
 \\ &&\leq \delta^{(n+1)/2}   e^{(a+1)\kappa \|w_0\|^2+ C_\kappa n },
  \end{eqnarray*}
  which yields (\ref{1003}).


Now  we  are in a  position to  prove  (\ref{1002}).
  By  Young's inequality and (\ref{30-2}),  we have
\begin{eqnarray}
\nonumber && (1+\sigma_{2n+2}-\sigma_{2n})^8 \exp\{8 \cC_0 \int_{\sigma_{2n}}^{\sigma_{2n+2}}\|w_s\|_1^{4/3}\dif s\} \|\rho_{\sigma_{2n}}\|^4
\\  \nonumber && \leq
(1+\sigma_{2n+2}-\sigma_{2n})^8
 \\ \nonumber && \times \exp\{\frac{\nu \kappa }{6} \int_{\sigma_{2n}}^{\sigma_{2n+1} } \|w_s\|_1^{2}e^{-\nu(\sigma_{2n+1}-s)+8\mathfrak{B}_0\kappa (\ell_{\sigma_{2n+1}}-\ell_s)}\dif s-\frac{\kappa}{6}\mathfrak{B}_0(\ell_{\sigma_{2n+1}}-\ell_{\sigma_{2n}})  \}
\\ \nonumber && \times  \exp\{\frac{\nu \kappa }{6} \int_{\sigma_{2n+1}}^{\sigma_{2n+2} } \|w_s\|_1^{2}e^{-\nu(\sigma_{2n+2}-s)+8\mathfrak{B}_0\kappa (\ell_{\sigma_{2n+2}}-\ell_s)}\dif s-\frac{\kappa}{6}\mathfrak{B}_0(\ell_{\sigma_{2n+2}}-\ell_{\sigma_{2n+1}}) \}
\\ \nonumber && \times \exp\big\{C_\kappa X_{2n}+ C_\kappa X_{2n+1}+\frac{\kappa}{6}\mathfrak{B}_0 (Y_{2n}+Y_{2n+1})\big\}
\\  \label{42-1} &&\times
C^{n/10}\big( \prod_{i=0}^{n-1} \theta_{2i}^{1/10}\big)  e^{\sum_{i=0}^{2n-1}C_\kappa (X_i+Y_i) },
\end{eqnarray}
  where $X_i=\int_{\sigma_{i}}^{\sigma_{i+1}} e^{2\nu(\sigma_{i+1}-s)-16 \mathfrak{B}_0\kappa (\ell_{\sigma_{i+1}}-\ell_s)}\dif s,Y_i=\ell_{\sigma_{i+1}}-\ell_{\sigma_i}.$
    Notice  that  the right hand of   (\ref{42-1})  depends on
  $\beta$ and $N$  through $\theta_i.$
We will determine  the values of $\beta$ and $N$ such that (\ref{1002})  hold.

{

For  $i\in \{2n,2n+1\}$,    set
\begin{eqnarray*}
\zeta_i=  \mE \exp\big\{\frac{\nu \kappa}{6}   \int_{\sigma_{i}}^{\sigma_{i+1} } \|w_s\|_1^{2}e^{-\nu(\sigma_{i+1}-s)+8\mathfrak{B}_0\kappa (\ell_{\sigma_{i+1}}-\ell_s)}\dif s -\frac{\kappa}{6} \mathfrak{B}_0(\ell_{\sigma_{i+1}}-\ell_{\sigma_{i}})\big\}.
\end{eqnarray*}
By Lemma  \ref{qu-2},
  we have
  \begin{eqnarray}
 \label{2121-1} && \Big(\mE (1+\sigma_{2n+2}-\sigma_{2n})^{\frac{120}{7}} \Big)^{7/15}\leq C_\kappa,
        \\ \label{2121-2}   &&  \big( \mE \zeta_i^6 \big)^{1/6}
 \leq C_\kappa \mE \exp\big\{\kappa e^{-1}\|w_{\sigma_{i}}\|^2/6\big\}, \quad i\in \{2n, 2n+1\}.
  \end{eqnarray}
  By  (\ref{1003}), for any    $\delta\in (0,1)$  and  $\kappa \in (0,\kappa_0]$,
 there exist    constants   $\beta$ and $N$
 such that
  \begin{eqnarray}
  \label{2121-3}
  \big( \mE  \prod_{i=0}^{n-1} \theta_{2i}^{1/2}\big)^{1/5}
\leq    \delta^{n/10}   e^{2 a\kappa \|w_0\|^2/5+C_\kappa (n-1) /5 }.
  \end{eqnarray}
  Notice that  on the event  $\{\ell:\Theta(\ell) \leq M\}$,
it holds that
\begin{eqnarray}
\label{2121-5}
\sum_{i=0}^{2n+1} (X_i+Y_i)\leq  M(2n+2).
\end{eqnarray}

    Hence,  by   H\"older's inequality, using  (\ref{42-1})--(\ref{2121-5}) and (\ref{30-7}),  we have
  \begin{eqnarray*}
  && \mE \Big[ (1+\sigma_{2n+2}-\sigma_{2n} )^8 \exp\{\cC_0 \int_{\sigma_{2n}}^{\sigma_{2n+2}}\|w_s\|_1^{4/3}\dif s\} \|\rho_{\sigma_{2n}}\|^4I_{\{\theta\leq M\}}\Big]
  \\ && \leq \mE \Big[(1+\sigma_{2n+2}-\sigma_{2n})^8  \zeta_{2n}\zeta_{2n+1}\big(\prod_{i=0}^{n-1} \theta_{2i}^{1/10}\big) e^{C_\kappa M(2n+2)}C^{n/10} \Big]
  \\ &&\leq  C_{\kappa,M}^n  \Big(\mE (1+\sigma_{2n+2}-\sigma_{2n})^{\frac{120}{7}} \Big)^{7/15}
  (\mE \zeta_{2n}^6)^{1/6}  (\mE \zeta_{2n+1}^6)^{1/6} (\mE \prod_{i=0}^{n-1} \theta_{2i}^{1/2} )^{1/5}
  \\ &&    \leq  C_{\kappa,M}^n\delta^{n/10}e^{ a\kappa \|w_0\|^2} \mE \exp\big\{\kappa e^{-1}\|w_{\sigma_{2n}}\|^2/6\big\} \mE \exp\big\{\kappa e^{-1}\|w_{\sigma_{2n+1}}\|^2/6\big\}
  \\ &&\leq C_{\kappa,M}^n\delta^{n/10}e^{4 a\kappa \|w_0\|^2}, \quad\forall n\in \mN,
  \end{eqnarray*}
  where $ C_{\kappa,M}\geq 1$ is a constant
  depending  on $\kappa,M$ and  $\nu,\{b_j\}_{j\in \cZ_0},\nu_S,d.$
  Choose  $\delta=\delta(\kappa,M,\gamma_0, \nu,\{b_j\}_{j\in \cZ_0},\nu_S,d)$ sufficiently small so that
  \begin{eqnarray*}
      C_{\kappa,M}^n\delta^{n/10} \leq  e^{- n \gamma_0}, \quad \forall n\in \mN.
  \end{eqnarray*}
  Then, we adjust  the values  of $\beta$ and $N$ according to (\ref{2121-3}).   The proof of  (\ref{1002})  is complete.

  }

\end{proof}

\subsection{The control of $
\int_0^{\ell_t}v(s)\dif W(s).$ }
\label{40-3}
For any $M,t>0$ and  $n\in \mN,$
the aim of this subsection is to give an estimate for the moment of the stochastic integral:
 \begin{equation*}
 \mE  \Big[  \Big| \int_{\ell_{\sigma_{2n}} }
  ^{ (\ell_{\sigma_{2n+1}}\wedge \ell_t)\vee \ell_{\sigma_{2n}}  } v(s)\dd W(s)\Big|^2
    I_{\{\Theta \leq M\} }
  \Big].
  \end{equation*}
 We start with an  estimate on the moments of $\rho_t.$
\begin{lemma}
\label{L:3.3}
  For any $ \kappa \in (0,\kappa_0], M>0,\gamma_0>0,$
  let $\beta$ be the constant chosen according to  Lemma  \ref{30-1}. Then,
  for any $w_0\in  H, n\in \mN$ and $t\ge0$,  one has
  \begin{equation}
  \label{3.6}
    \E_{w_0}\big[ \|\rho_t\|^4   I_{\{ \Theta \leq M\}}I_{\{t\in
    [{\sigma_{2n}},{\sigma_{2n+2}})\}}\big]   \le  C_{\kappa,M,\gamma_0}
     \exp\left\{4 a \kappa   \|w_0\|^2- n\gamma_0   \right\},
  \end{equation}
  where  $C_{\kappa,M,\gamma_0}  $
  is a constant depending on
  $\kappa,M,\gamma_0$ and  $\nu,\{b_j\}_{j\in \cZ_0},\nu_S,d.$
\end{lemma}
\begin{proof}
  From the construction we have
 $$
 \rho_t=\begin{cases} J_{\sigma_{2n},t}\rho_{\sigma_{2n}}-
 \cA_{\sigma_{2n},t}v_{{\sigma_{2n}},t}, &
 \text{for }t \in [{\sigma_{2n}},{\sigma_{2n+1}}],
  \\ J_{\sigma_{2n+1},t}\rho_{\sigma_{2n+1}}, &
  \text{for }t\in [{\sigma_{2n+1}},{\sigma_{2n+2}}] \end{cases}	
 $$ for any  $n\ge0$.
  { Using \eqref{p-1} and inequalities   \eqref{2.8}\eqref{2.10}}, we get
    \begin{equation} \label{3.8}
    \|v_{\sigma_{2n},\sigma_{2n+1}}\|_{L^2([\ell_{\sigma_{2n}},\ell_{\sigma_{2n+1}}];
    \R^d)} \le
    \beta^{-1/2}\|J_{\sigma_{2n},\sigma_{2n+1}} \rho_{\sigma_{2n}}\|.
\end{equation}
Hence,   by Lemma \ref{15-4}, (\ref{2.7}) 
and the definition of $\sigma_{2n+1}$, for any $t\in [\sigma_{2n},\sigma_{2n+1}]$,
   \begin{align}
    \|\rho_t\|&\le  \|J_{\sigma_{2n},t}\rho_{\sigma_{2n}}\|+
    \|\aA_{\sigma_{2n},t}v_{\sigma_{2n},t}\|  \nonumber
    \\  \nonumber  &\le  \|J_{\sigma_{2n},t}\rho_{\sigma_{2n}}\|+
    \|\aA_{{\sigma_{2n}},t}\|_{\cL(L^2([\ell_{\sigma_{2n}},\ell_t];\R^d), H)}
    \|v_{{\sigma_{2n}},t}\|_{L^2([\ell_{\sigma_{2n}},\ell_{\sigma_{2n+1}}];\R^d)}
    \\  \nonumber  &\le  \|J_{\sigma_{2n},t}\rho_{\sigma_{2n}}\|+
    \|\aA_{{\sigma_{2n}},t}\|_{\cL(L^2([\ell_{\sigma_{2n}},\ell_t];\R^d), H)}
    \|v_{{\sigma_{2n}},{\sigma_{2n+1}}}\|_{L^2([\ell_{\sigma_{2n}},
    \ell_{\sigma_{2n+1}}];\R^d)}
    \\ \nonumber
     &\le \|J_{\sigma_{2n},t}\rho_{\sigma_{2n}}\| +C \beta^{-1/2}(\ell_{\sigma_{2n+1}}-\ell_{\sigma_{2n}})^{1/2}
     \|J_{{\sigma_{2n}},{\sigma_{2n+1}}}\rho_{\sigma_{2n}}
     \| \sup_{s\in [{\sigma_{2n}},t]}\|J_{s,t}
      \|_{\cL(H,H)}
      \\ \label{21-1} &\leq C_{\kappa,\beta}  (1+\sigma_{2n+1}-\sigma_{2n}) \exp\big\{ \cC_0\int_{\sigma_{2n}}^{\sigma_{2n+1}}\|w_s\|_1^{4/3}\dif s \big\} \|\rho_{\sigma_{2n}}\|,
  \end{align}
  where $C_{\kappa,\beta}$ is a constant depending on $\kappa,\beta,\nu,\{b_j\}_{j\in \cZ_0},\nu_S,d.$
 For any
     $t\in [\sigma_{2n+1},\sigma_{2n+2}],$ also by  Lemma \ref{15-4},  it holds that
\begin{eqnarray}
\label{21-2}
\|\rho_t\|\le  \sup_{t\in [\sigma_{2n+1},\sigma_{2n+2}]}\|J_{\sigma_{2n+1},t}\rho_{\sigma_{2n+1}}\|\leq  \cC_0 \exp\{\cC_0 \int_{\sigma_{2n+1}}^{\sigma_{2n+2}}\|w_s\|_1^{4/3}\dif s \} \|\rho_{\sigma_{2n+1}}\|.
\end{eqnarray}

Combining (\ref{21-1}),(\ref{21-2}) with   (\ref{1002}), we complete the proof.

\end{proof}

\begin{lemma}
\label{L:3.4}
  For any $M>0,\gamma_0>0, \kappa \in (0,\kappa_0]$,
 let $\beta$ be the constant chosen according to Lemma  \ref{30-1}. Then,
  \begin{equation}
  \label{3.10}
  \begin{split}
  &  \mE_{w_0}  \Big[  \big| \int_{\ell_{\sigma_{2n}} }
  ^{ \ell_{\sigma_{2n+1}} } v(s)\dd W(s)\big|^2
    I_{\{\Theta \leq M\} }
  \Big]
     \\ &  \le   C_{\kappa,M,\gamma_0}  \exp\{2 \kappa a  \|w_0\|^2-\gamma_0   n /2\}
     \end{split}
  \end{equation}
  and
   \begin{equation}
  \label{3-10}
  \begin{split}
  &  \mE_{w_0}  \Big[  \big| \int_{\ell_{\sigma_{2n}}}
  ^{ (\ell_{\sigma_{2n+1}}\wedge \ell_t)\vee \ell_{\sigma_{2n}}  } v(s)\dd W(s)\big|^2
    I_{\{\Theta \leq M\} }
  \Big]
     \\ &  \le   C_{\kappa,M,\gamma_0}  \exp\{2 \kappa a \|w_0\|^2-\gamma_0  n/2 \}
     \end{split}
  \end{equation}
  for
  $n\geq 0, t\geq 0$  and $w_0\in  H$, here $C_{\kappa,M,\gamma_0} $
is a constant depending on
  $\kappa,M,\gamma_0$ and $\nu,\{b_j\}_{j\in \cZ_0},\nu_S,d.$
\end{lemma}
\begin{proof}  We only prove (\ref{3.10});
 the estimate   (\ref{3-10})
 is treated in a similar way.
  Using  the generalised It$\hat{\text{o}}$
  isometry (see Section~1.3 in~\cite{nualart2006}) and
  the fact that $v(t)=0$ for $t\in [\ell_{\sigma_{2n+1}},\ell_{\sigma_{2n+2}} ]$,
   we have
  \begin{align}
  \nonumber &\E_{w_0}  \[ \big| \int_{\ell_{\sigma_{2n}}}^{\ell_{\sigma_{2n+1}}} v(s)\dd  W(s)\big|^2  I_{\{\Theta \leq M\}}\]
   \\ & =  \E_{w_0}  \[ \int_{\ell_{\sigma_{2n}}}^{\ell_{\sigma_{2n+1}}}   |v(s)|_{\R^d}^2 \dd  sI_{\{\Theta \leq M\}}\]
     \nonumber
     \\&\quad+\E_{w_0} \[ \int_{\ell_{\sigma_{2n}}}^{\ell_{\sigma_{2n+1}}}
     \!\!\!\int_{\ell_{\sigma_{2n}}}^{\ell_{\sigma_{2n+1}}}
      \text{Tr}(\cD_sv(r)\cD_rv(s))\dd  s \dd  rI_{\{\Theta \leq M\}} \] \nonumber
   \\  \nonumber   & \le    \E_{w_0} \[ \int_{\ell_{\sigma_{2n}}}^{\ell_{\sigma_{2n+1}}}
     |v(s)|_{\R^d}^2 \dd  s I_{\{\Theta \leq M\}}\] \nonumber\\
     &\quad+\E_{w_0} \[
      \int_{\ell_{\sigma_{2n}}}^{\ell_{\sigma_{2n+1}}}
     \!\!\!\int_{\ell_{\sigma_{2n}}}^{\ell_{\sigma_{2n+1}}}
      |\cD_rv_{{\sigma_{2n}},{\sigma_{2n+1}}}(s)|_{\R^d\times \R^d}^2
     \dd  s\dd  r I_{\{\Theta \leq M\}}\]  \nonumber\\  & =L_1+L_2.
  \label{3.11}
  \end{align}
  Using \eqref{p-1},  \eqref{3.8},   Lemma \ref{15-4} and Lemma  \ref{30-1}, we have
     \begin{align}
  \nonumber  L_1&= \E_{w_0}   \int_{\ell_{\sigma_{2n}}}^{\ell_{\sigma_{2n+1}}}
   |v(s)|_{\R^d}^2 I_{\{\Theta \leq M\}} \dd  s
  \le
    \beta^{-1}\E_{w_0}  \|J_{\sigma_{2n},\sigma_{2n+1}}\rho_{\sigma_{2n}}\|^2I_{\{\Theta \leq M\}}
    \\  \nonumber & \leq   \beta^{-1}  \cC_0   \mE \Big[ e^{ \cC_0  \int_{\sigma_{2n}}^{\sigma_{2n+1}}  \|w_r\|_1^{4/3}\dif r}\| \rho_{\sigma_{2n}}\|^2I_{\{\Theta \leq M\}}\Big]
   \\ \label{3.12}
   &\le  C_{\kappa,M,\gamma_0} \exp\{  2 \kappa a \|w_0\|^2-n \gamma_0/2\}
  \end{align}
  where $\kappa\in (0,\kappa_0]$ and $C_{\kappa,M,\gamma_0}$ is constant depending on $\kappa,M,\gamma_0$ and $\nu,\{b_j\}_{j\in \cZ_0},\nu_S,d$
  that may change  from line to line.
   To estimate $L_2$,   we use the explicit form of   $\cD_r v.$
  Notice  that,   for any $r\in[\ell_{\sigma_{2n}}, \ell_{\sigma_{2n+1}}]$
  and~$i=1,\dots, d$,
\begin{align*}
   \cD_r^iv_{{\sigma_{2n}},{\sigma_{2n+1}}}
 &= \cD_r^i(\cA_{{{\sigma_{2n}},{\sigma_{2n+1}}}}^*)
 (\cM_{{{\sigma_{2n}},{\sigma_{2n+1}}}}
 +\beta \I)^{-1}J_{{{\sigma_{2n}},{\sigma_{2n+1}}}}\rho_{\sigma_{2n}}
    \\ &\quad + \cA_{{{\sigma_{2n}},{\sigma_{2n+1}}}}^*(\cM_{{{\sigma_{2n}},{\sigma_{2n+1}}}}
    +\beta\I)^{-1}
 \\ & \quad\quad\times \Big( \cD_r^i(\cA_{{{\sigma_{2n}},{\sigma_{2n+1}}}})
 \cA_{{{\sigma_{2n}},{\sigma_{2n+1}}}}^*
  +\cA_{{{\sigma_{2n}},{\sigma_{2n+1}}}}\cD_r^i(\cA_{{{\sigma_{2n}},{\sigma_{2n+1}}}}^*)
  \Big) \\ & \quad\quad\times
 (\cM_{{{\sigma_{2n}},{\sigma_{2n+1}}}}+\beta\I)^{-1}
 J_{{{\sigma_{2n}},{\sigma_{2n+1}}}} \rho_{\sigma_{2n}}
 \\ & \quad +\cA_{{{\sigma_{2n}},{\sigma_{2n+1}}}}^* (\cM_{{{\sigma_{2n}},
 {\sigma_{2n+1}}}}+\beta\I)^{-1}
 \cD_r^i(J_{{{\sigma_{2n}},{\sigma_{2n+1}}}}) \rho_{\sigma_{2n}}.
\end{align*}
By inequalities \eqref{2.8}--\eqref{2.10}, Lemma \ref{15-4} and Lemma \ref{L:2.3}, we have
\begin{align*}
  &  \|\cD_r^iv_{{\sigma_{2n}},{\sigma_{2n+1}}} \|_{L^2([\ell_{\sigma_{2n}},\ell_{\sigma_{2n+1}}];
   \R^d)}
 \\  & \le   \beta^{-1}\| \cD_r^i(\cA_{\sigma_{2n},\sigma_{2n+1}})
  \|_{\cL(L^2([\ell_{\sigma_{2n}},\ell_{\sigma_{2n+1}} ];\R^d),
  {H})}
  \|J_{\sigma_{2n},\sigma_{2n+1}} \rho_{\sigma_{2n}} \|
    \\ &\quad +2\beta^{-1}\| \cD_r^i(\cA_{\sigma_{2n},\sigma_{2n+1}}^*)
    \|_{\cL({H},L^2([\ell_{\sigma_{2n}},\ell_{\sigma_{2n+1}}];\R^d))}
  \|J_{\sigma_{2n},\sigma_{2n+1}} \rho_{\sigma_{2n}} \|  \\ &\quad
  + \beta^{-1/2} \| \cD_r^i(J_{\sigma_{2n},\sigma_{2n+1}})
   \rho_{\sigma_{2n}} \|
   \\ &\leq C_{\kappa,\beta}  (1+\sigma_{2n+1}-\sigma_{2n}) \exp\{2 \cC_0 \int_{\sigma_{2n}}^{\sigma_{2n+1}}\|w_s\|_1^{4/3}\dif s \} \|\rho_{\sigma_{2n}}\|,
\end{align*}
where  $C_{\kappa,\beta}$ is a  constant depending on $\kappa,\beta$ and $\nu,\{b_j\}_{j\in \cZ_0},\nu_S,d.$
By Lemma \ref{30-1}  and the fact that $ \ell_{\sigma_{2n+1}}-\ell_{\sigma_{2n}} \leq \frac{\nu}{8 \mathfrak{B}_0 \kappa }(\sigma_{2n+1}-\sigma_{2n})$, it follows that
\begin{align*}
 &  \E_{w_0} \int_{\ell_{\sigma_{2n}}}^{\ell_{\sigma_{2n+1}}} \!\!\! \int_{\ell_{\sigma_{2n}}}^{\ell_{\sigma_{2n+1}}}  |
  \cD_rv_{{\sigma_{2n}},{\sigma_{2n+1}} }(s)|_{\R^d\times \R^d}^2I_{\{\Theta \leq M\}} \dd
   s\dd  r
   \\ & \leq  \mE \Big[(\ell_{\sigma_{2n+1}}-\ell_{\sigma_{2n}})
   C_{\kappa,\beta}  (1+\sigma_{2n+1}-\sigma_{2n})^2  \exp\{4 \cC_0 \int_{\sigma_{2n}}^{\sigma_{2n+1}}\|w_s\|_1^{4/3}\dif s\} \|\rho_{\sigma_{2n}}\|^2
   \Big]
   \\ &    \leq C_{\kappa,M,\gamma_0} \exp\{ 2 \kappa a \|w_0\|^2-n \gamma_0/2 \}.
\end{align*}
Combining the above estimate with (\ref{3.11}) and (\ref{3.12}), we complete the proof.
\end{proof}

\subsection{Proof of Proposition  \ref{3-11}(continued).}
\label{17-27}


For $\xi\in H$, let $v$ be the process chosen as in (4.7). In order to treat the term $I_1$ in (\ref{20-3}), we
 observe that
  \begin{eqnarray*}
    && | \nabla_\xi P_t^{M} f(w_0) |=
    | \mE \nabla f(w_t)J_{0,t}\xi
    I_{\{\Theta\leq M\}} |
    \\ &&
    =| \mE \[ \nabla f(w_t)\cD^v w_t
    I_{\{\Theta\leq M\}}\]+\mE \[ \nabla  f(w_t)\rho_t
    I_{\{\Theta\leq M\}}\] |
    \\  &&
    = \mE \[  f(w_t)\int_0^{\ell_t}  v(s)\dif W(s)
    I_{\{\Theta\leq M\}}\]+\mE \[ \nabla  f(w_t) \rho_t
    I_{\{\Theta\leq M\}}\]
    \\ &&:=I_{11}+I_{12}.
    \end{eqnarray*}
    For  any  $M>0,$ $\kappa \in (0,\kappa_0]$
  and   $\gamma_0>0$,
   we set the value of   $\beta$ in  (\ref{p-1}) according to  Lemma  \ref{30-1}.
For the term $I_{11},$ by    Lemma \ref{L:3.4} in subsection \ref{40-3},  we have
    \begin{eqnarray*}
    I_{11}&\leq &
    \|f\|_\infty   \sum_{n=0}^{\infty }\mE
    \big|\int_{\ell_{\sigma_{2n}}  }^{(\ell_{\sigma_{2n+1}} \wedge \ell_t)\vee \ell_{\sigma_{2n}}}v(s)\dif W_s
    I_{\{\Theta\leq M\}}  \big|
    \\ &\leq &  \|f\|_\infty   \sum_{n=0}^{\infty}\Big(\mE
    \big|\int_{\ell_{\sigma_{2n}}  }^{(\ell_{\sigma_{2n+1}} \wedge \ell_t)\vee \ell_{\sigma_{2n}}} v(s)\dif W_s
    \big|^2
    I_{\{\Theta \leq M\} } \Big)^{1/2}
    \\ & \leq  &  \|f\|_\infty   \sum_{n=0}^{\infty} C_{\kappa,M,\gamma_0}\exp\{\kappa a \|w_0\|^2-\gamma_0  n/4 \}.
  \end{eqnarray*}
  Now consider the term $I_{12}.$ By Lemma  \ref{L:3.3} in subsection \ref{40-3},  we have
  \begin{eqnarray*}
   I_{12} &\leq & \sum_{n=0}^\infty \|\nabla f\|_\infty
    \E_{w_0}\big[ \|\rho_t\|   I_{\{ \Theta \leq M\}}I_{\{t\in
    [{\sigma_{2n}},{\sigma_{2n+2}})\}}\big]
   \\ &\leq  &  \sum_{n=0}^\infty \|\nabla f\|_\infty C_{\kappa,M,\gamma_0}
     \exp\left\{a \kappa \|w_0\|^2-\gamma_0  n/4 \right\}.
  \end{eqnarray*}
  Combining the estimates of $I_{11},I_{12}$, for any $\xi$ with $\| \xi\| $=1,  we conclude that
\begin{eqnarray}
\label{180-3}
  | \nabla_\xi P_t^{M} f(w_0) |\leq C_{\kappa,M,\gamma_0} \big(\|f\|_\infty+\|\nabla f\|_\infty \exp\left\{a \kappa  \|w_0\|^2\right\}\big),
\end{eqnarray}
where $C_{\kappa,M,\gamma_0}$ is a constant  independent  of $t.$
Let $\gamma(s)=sw_0+(1-s)w_0'$. Then, by  (\ref{180-3}),
\begin{eqnarray}
\nonumber && \mE f(w_t^{w_0})I_{\{\Theta\leq M\}}
  -\mE f(w_t^{w_0'})I_{\{\Theta\leq M\}}
  = \int_0^1 \langle \nabla  P_t^{M}f(\gamma(s)) ,w_0-w_0'\rangle \dif s
  \\  \label{20-2} &&\leq C_{\kappa,M,\gamma_0} \|w_0-w_0'\|\big(\|f\|_\infty+\|\nabla f\|_\infty \sup_{s\in [0,1]}\exp\{a \kappa  \|\gamma(s)\|^2\}\big).
\end{eqnarray}

Combining the estimates (\ref{20-2}) with (\ref{20-3}),
 for any $\kappa\in (0,\kappa_0]$, $w_0,w_0'\in B_H(\Upsilon), f\in C_b^1(H)$  and   $M\geq 1,t>0,$  we obtain that
 \begin{eqnarray*}
  && | P_t f(w_0)-P_t f(w_0')|
  \\ && \leq C_{\kappa,M, \gamma_0 } \|w_0-w_0'\|\big(\|f\|_\infty+\|\nabla f\|_\infty \exp\{ a \kappa   \Upsilon^2\}\big)
   + 2\|f\|_\infty   \mP(\Theta \geq M ).
\end{eqnarray*}

For any  bounded and Lipschitz continuous function $f$ on $H,$
by the arguments  in \cite[Page 1431]{KPS10},
there exists a sequence $(f_n)$ satisfies $(f_n)\subseteq C_b^1(H)$
and $\lim_{n\rightarrow  \infty}f_n(x)=f(x)$ pointwise.
In addition, $\|f_n\|_\infty\leq \|f\|_\infty$ and $\|\nabla f_n\|_\infty\leq Lip(f),$ where  $Lip(f)=\sup_{x\neq y}\frac{|f(x)-f(y)|}{\|x-y\|}.$
Therefore, for any $t\geq 0,$ one has
\begin{eqnarray}
\nonumber && | P_t f(w_0)-P_t f(w_0')|=\lim_{n\rightarrow \infty}| P_t f_n(w_0)-P_t f_n(w_0')|
\\ \nonumber &&\leq  \lim_{n\rightarrow \infty} \Big[C_{\kappa,M, \gamma_0 } \|w_0-w_0'\|\big(\|f_n\|_\infty+\|\nabla f_n \|_\infty \exp\{ a \kappa   \Upsilon^2\}\big)
   + 2\|f_n\|_\infty   \mP(\Theta \geq M )\Big]
   \\ \nonumber && \leq  C_{\kappa,M, \gamma_0 } \|w_0-w_0'\|(\|f\|_\infty+Lip(f)  \exp\big\{ a \kappa   \Upsilon^2\big\}
   + 2\|f\|_\infty   \mP(\Theta \geq M )
   \\ \label{02-1}  &&:=J_1+J_2.
\end{eqnarray}
For  any $\eps>0,$
 by (\ref{1609}),   we  can  set  a $M>0$ such that  $J_2<\frac{\eps}{2}.$
Obviously, there exists a $\delta>0$ such that  for any $w_0,w_0'\in B_H(\Upsilon)$ with $\|w_0-w_0'\|<\delta$,
$
  J_1<\frac{\eps}{2}.
$
Combining these with (\ref{02-1}), one arrives at
that for any positive constants   $\eps,\Upsilon,$ there exists a
$\delta>0,$ such that
\begin{eqnarray*}
| P_t f(w_0)-P_t f(w_0')|<\eps, ~~\forall t\geq 0 \text{ and }
w_0,w_0'\in B_H(\Upsilon) \text{ with  } \|w_0-w_0'\|<\delta.
\end{eqnarray*}
This  completes  the proof of Proposition  \ref{3-11}. 

\section{Proof of weak irreducibility}
\label{54-2}
We start with the following lemma.
\begin{lemma}
\label{16-1}
For any $T>0,\eps>0$  and non-zero reals numbers $b_i,i\in \cZ_0,$  one has
\begin{eqnarray*}
\mP(\sup_{t\in [0,T]}\|\sum_{i\in\cZ_0}b_iW^i_{S_t}e_i\|<\eps)\geq p_0>0.
\end{eqnarray*}
where $p_0=p_0({T,\eps,\{b_i\}_{i\in \cZ_0}})$ is a constant.
\end{lemma}
\begin{proof}
For sufficiently big constant $M>0$, one sees  that
\begin{eqnarray*}
&& \mP(\sup_{t\in [0,T]}\|\sum_{i\in\cZ_0}b_iW^i_{S_t}e_i\|<\eps)
\\ &&\geq  \mP^{\mu_\mS}(S_T\leq M)
\mP^{\mu_\mW}(\sup_{t\in [0,M]}\|\sum_{i\in\cZ_0}b_iW^i_{t}e_i\|<\eps)
>0,
\end{eqnarray*}
which completes  the proof.

\end{proof}


Now we are in position to give a proof of Proposition \ref{16-6}.

\begin{proof}
The proof of this proposition is  the same as that in
 \cite[Lemma 3.1]{EM01}.
 For the reader's convenience, we provide the proof here.
 { Define $v_t=w_t- \eta_t$,  where $ \eta_t=\sum_{i\in\cZ_0}b_iW^i_{S_t}e_i
  $}, $w_t$ is the solution to (\ref{0.1}) at time $t$. Then, $v_t$ satisfies
  \begin{eqnarray*}
    \frac{\partial v_t}{\partial t}=\nu \Delta (v_t+ \eta_t)
    +B(\cK  w_t,w_t)=\nu \Delta (v_t+ \eta_t)
    +B(\cK  w_t,v_t+ \eta_t)
  \end{eqnarray*}
  Taking the $L^2$-inner product of this equation with $v_t$ produces
  \begin{eqnarray*}
    \frac{1}{2}\frac{d}{dt}\|v_t\|^2 &=& -\nu \|\nabla v_t\|^2+
    \langle \nu \Delta \eta_t, v_t\rangle+
    \langle B(\cK w_t,  \eta_t), v_t\rangle
    \\
    &\leq &  -\nu \|\nabla v_t\|^2 +C_1 \|v_t\| \|\Delta  \eta_t\|
    +C_1 \|v_t\|\|\Delta  \eta_t\| \|\cK w_t \|_1
     \\
    &= &  -\nu \|\nabla v_t\|^2 +C_1  \|v_t\| \|\Delta  \eta_t\|
    +C_1 \|v_t\|\|\Delta  \eta_t\| \| (v_t+ \eta_t) \|
    \\ &\leq & -\nu \|\nabla v_t\|^2 +C_1 \|v_t\| \|\Delta  \eta_t\|
    +C_1 \|v_t\|\|\Delta  \eta_t\| \|v_t \|
    +C_1 \|v_t\|\|\Delta  \eta_t\| \| \eta_t\|
     \\ &\leq & -\frac{\nu}{2} \|\nabla v_t\|^2 +\frac{4 C_1^2 }{\nu } \|\Delta \eta_t\|^2
    +C_1 \|v_t\|^2 \|\Delta  \eta_t\|
    +\frac{4 C_1^2 }{\nu } \|\Delta  \eta_t\|^2  \| \eta_t \|^2,
  \end{eqnarray*}
  where $C_1=C_1(\nu)$ is a  constant.
  For any $T,\delta>0$, we define,
  \begin{eqnarray*}
    \Omega'(\delta,T)&=& \left\{
    g=(g_s)_{s\in [0,T]}\in D([0,T];H):\sup_{s\in [0,T]}\|\Delta g_s \|\leq \min \{\delta,
    \frac{\nu}{4C_1}\}\right\},
  \end{eqnarray*}
  where $\Delta$ stands for the Laplacian operator.
  If $\eta \in  \Omega'(\delta,T)$,  one has
  \begin{eqnarray*}
    \|v_t\|^2 \leq  \|v_0\|^2e^{-\frac{\nu}{2} t}
    + \frac{4 C_1^2 }{\nu  }\cdot \frac{2}{\nu}\cdot \Big[\min\left(\delta,\frac{\nu}{4C_1}\right)^4+
    \min\left(\delta,\frac{\nu}{4C_1}\right)^2\Big].
  \end{eqnarray*}
  Let $\mathcal{C}$ and $\gamma$ be given as in the statement of Proposition \ref{16-6}.
  As $\|w_0\|\leq \mathcal{C},$ there exists a  $T$
  and a $\delta$
  such that
  \begin{eqnarray*}
    \|v_T\|\leq \frac{\gamma}{2} \text{ and } \delta \leq \frac{\gamma}{2} .
  \end{eqnarray*}
  Putting everything together, one has
  \begin{eqnarray*}
    \|w_0\|\leq  \mathcal{C} \text{ and  }
     \eta \in  \Omega'(\delta,T) \Rightarrow
    \|w_T\|\leq \|v_T\|+\| \eta_T\|\leq \gamma.
  \end{eqnarray*}
  Combining this fact  with Lemma \ref{16-1}, we complete the proof.
\end{proof}

\appendix{}

\section{Proof of  Lemma   \ref{qu-2}. }

\label{ssss-1}

This section is organized as follows. In the subsection  \ref{sss-1}, we make some preparations.
  Then, we provide the proofs of (\ref{2727-1})--(\ref{30-7}) in subsection \ref{A-1-3}, the proof of (\ref{100-1}) is given in  subsection \ref{A-1-4}.

\subsection{Preparations}\label{sss-1}

For $\kappa>0,\eps\in (0,1]$ and $\ell\in \mS$, set
  \begin{eqnarray*}
    \ell_t^\eps=\frac{1}{\eps}\int_t^{t+\eps}\ell_s\dif s+\eps t,
  \end{eqnarray*}
  and
  \begin{eqnarray*}
  \sigma^\eps=  \sigma^\eps(\ell):
  = \inf\big\{t\geq 0:\nu t-8\mathfrak{B}_0\kappa \ell_t^\eps  >    1 \big\}.
  \end{eqnarray*}
Keeping in mind that $\ell$ is a $\rm c\grave{a}dl\grave{a}g$
 increasing function from $\mR^{+}$
 to $\mR^{+}$  with $\ell_0 = 0$, it is easy to see that the following lemma is valid.
\begin{lemma}
 \label{qu-5}
 For $\ell\in \mS$,
\begin{itemize}
  \item[(i)] $\ell_\cdot^\eps:[0,\infty)\rightarrow[0,\infty)$ is continuous and strictly increasing;

  \item[(ii)] for any $t\geq0$, $\ell_t^\eps$ strictly decreases to $\ell_t$ as $\eps$ decreases to 0.
 \end{itemize}
\end{lemma}

 With regard to stopping times $\sigma^\eps$ and $\sigma,$  the following moment estimates hold.
  \begin{lemma}\label{Lemma A2}
     There exists a  constant $\widetilde{\kappa}_0>0$ such that, for any $\kappa\in(0,\widetilde{\kappa}_0]$,
              \begin{eqnarray}
          \label{11-5}  \sup_{\eps\in(0,1]}\mE^{\mu_{\mS}} \exp\{10 \nu \sigma^\eps \} & \leq &  C_\kappa,
          \\  \label{11-6}
           \mE^{\mu_\mS} e^{10 \nu \sigma} &\leq & C_\kappa,
        \end{eqnarray}
           where  $C_\kappa$ is a constant depending on
  $\kappa$, $\nu,\{b_j\}_{j\in \cZ_0},\nu_S,d$.
  \end{lemma}

  \begin{proof}{
  We only give a proof for(\ref{11-5}); the proof of (\ref{11-6}) is similar.
  Let $c_\kappa =\int_{0}^\infty (e^{160 \mathfrak{B}_0 \kappa u}-1)\nu_{S}(\dif u).$ Then $c_\kappa<\infty$ for sufficiently small constant $\kappa$.
  For any $n\in \mN$ and $\eps\in (0,1),$ by the fact that $e^{160 \mathfrak{B}_0 \kappa S_t(\ell)-c_\kappa  t}=e^{160 \mathfrak{B}_0 \kappa \ell_t-c_\kappa  t}$ is a local martingale(c.f.\cite[Corollary 5.2.2]{Dav-2009})  and  that $\ell^\eps_n\leq \ell_{n+\eps}+\eps n$,  we have
  \begin{eqnarray}
   \nonumber  && \mP^{\mu_\mS}(\ell: \sigma^\eps(\ell)>n)\leq  \mP^{\mu_\mS}(\nu n-8\mathfrak{B}_0 \kappa \ell^\eps_n\leq 1)
    \\ \nonumber  &&= \mP^{\mu_\mS}(160 \mathfrak{B}_0 \kappa \ell^\eps_n\geq 20  \nu n-  20 )\\ \nonumber
    &&\leq  \mP^{\mu_\mS}(160 \mathfrak{B}_0 \kappa \ell_{n+\eps} \geq (20 \nu-160 \mathfrak{B}_0 \kappa \eps) n-  20 )
    \\ \nonumber  &&{ = \mP^{\mu_\mS}(160 \mathfrak{B}_0 \kappa \ell_{n+\eps} -c_\kappa(n+\eps) \geq (20 \nu-160 \mathfrak{B}_0 \kappa \eps-c_\kappa ) n- 20-c_\kappa \eps )}
    \\  \label{28-1}  && { \leq \exp\{- (20 \nu-160 \mathfrak{B}_0 \kappa \eps-c_\kappa ) n+20+c_\kappa \eps  \}.
    }
  \end{eqnarray}
By the Condition   \ref{14-2},  one has   $\lim_{\kappa \rightarrow 0}c_\kappa =0.$ Therefore,   there exists a  constant $\widetilde{\kappa}_0>0$ such that
  \begin{eqnarray*}
    160\mathfrak{B}_0 \kappa   +c_\kappa <10 \nu, \quad \forall \kappa \in (0, \widetilde{\kappa}_0].
  \end{eqnarray*}
 For any $\kappa\in (0, \widetilde{\kappa}_0]$ and $\eps\in(0,1]$, by (\ref{28-1})  and the above inequality,  we conclude that
  \begin{eqnarray*}
    && \mE^{\mu_\mS} e^{10 \nu \sigma^\eps}\leq 1+\sum_{n=0}^\infty e^{10 \nu (n+1) }\mP(\sigma^\eps \in ( n,n+1] )
    \\ &&{ \leq 1+\sum_{n=0}^\infty e^{10 \nu (n+1)}\exp\{- (20 \nu-160 \mathfrak{B}_0 \kappa \eps-c_\kappa ) n+20+c_\kappa \eps  \}}
    \\ && { \leq   1+\frac{\exp\{10\nu+20 +c_\kappa  \}}{1-\exp\{-10 \nu +160\mathfrak{B}_0 \kappa   +c_\kappa  \}}
    <\infty.}
  \end{eqnarray*}
  The proof is complete.
  }
  \end{proof}

  For any $\kappa\in(0,\widetilde{\kappa}_0]$, set
  \begin{eqnarray}\label{eq S1}
  \mathbb{S}_1=\{\ell\in\mathbb{S}:\sigma^1(\ell)<\infty\}.
  \end{eqnarray}
 We have the following lemma.
\begin{lemma}\label{lemma A3}
      $\mP^{\mu_{\mS}}(\mathbb{S}_1)=1$. And for any $\ell\in \mathbb{S}_1$, the following statements hold.
 \begin{itemize}
   \item[(1)] For any $\eps\in(0,1)$, $\sigma^\eps<\sigma^1<\infty$ and $\nu \sigma^\eps-8\mathfrak{B}_0\kappa \ell_{\sigma^\eps}^\eps=    1$;
  \item[(2)]$\sigma^\eps$ strictly decreases to $\sigma$ as $\eps$ decreases to 0;
\item[(3)] $\ell_{\sigma^\eps}^{\eps}$ strictly decreases to $\ell_{\sigma}$ as $\eps$ decreases to 0;
\item[(4)] $\nu \sigma-8\mathfrak{B}_0\kappa \ell_{\sigma}=1$;
%
  \item[(5)] $\limsup_{\eps\rightarrow 0}\int_{0}^{\sigma^\eps} e^{-\nu (\sigma^\eps-s) +8\mathfrak{B}_0\kappa (\ell_{\sigma^\eps}^\eps-\ell_s^\eps)   }  \dif \ell_s^\eps
  \leq\ell_\sigma.$
\end{itemize}
\end{lemma}

\begin{proof}
Lemma \ref{Lemma A2} implies that $\mP^{\mu_{\mS}} (\mathbb{S}_1)=1$.

By Lemma \ref{qu-5} and the definition of $\sigma^\eps$, it is easy to see that (1) holds. Moreover, for any $\ell \in   \mathbb{S}_1,$  $\sigma^\eps$ strictly decreases to a constant $a$ as $\eps$ decreases to 0.
On the other hand, for any $\eps\in(0,1]$,
 \begin{eqnarray}\label{2626-3}
 \sigma^\eps >a\geq \sigma.
 \end{eqnarray}
%
%
%
 For any $s<a$,
  by the definition of $\sigma^\eps$,
  \begin{eqnarray*}
 \nu s -8\mathfrak{B}_0\kappa \ell_s^\eps  \leq 1,\ \ \forall\eps\in(0,1].
  \end{eqnarray*}
  Letting $\eps\rightarrow 0$, by Lemma \ref{qu-5}, we get
    \begin{eqnarray*}
 \nu s -8\mathfrak{B}_0  \kappa  \ell_s \leq 1.
  \end{eqnarray*}
  Hence $\sigma\geq a$.
  Combining this with (\ref{2626-3}), we get (2).

  Using Lemma \ref{qu-5}, (2), the definition of $\ell^\eps$, and the fact that $\ell$ is increasing yields, it follows that for any $\eps\in(0,1]$,
  \begin{eqnarray*}
  \ell_{\sigma}<\ell^\eps_{\sigma}\leq \ell^\eps_{\sigma^\eps}\leq \ell_{\sigma^\eps+\eps}+\eps\sigma^\eps.
  \end{eqnarray*}
  Letting $\eps\rightarrow 0$, by (2) and the right { continuity}   of $\ell$, we get (3).
  Combining (1), (2) and (3), one gets (4).

  Combining (1) with the fact that for any $s\leq \sigma^\eps$, $\nu s-8\mathfrak{B}_0 \kappa \ell_s^\eps \leq 1$, one arrives at
 \begin{eqnarray*}
  &&\limsup_{\eps\rightarrow 0}\int_{0}^{\sigma^\eps} e^{-\nu (\sigma^\eps -s) +8\mathfrak{B}_0\kappa (\ell_{\sigma^\eps}^\eps-\ell_s^\eps)   }  \dif \ell_s^\eps\\
  &\leq&  \limsup_{\eps\rightarrow 0}
 e^{-1+1  }\int_{0}^{\sigma^\eps }   \dif \ell_s^\eps
 \leq  \limsup_{\eps\rightarrow 0} \ell_{\sigma^\eps }^\eps = \ell_{\sigma},
 \end{eqnarray*}
completing the proof of (5).

The proof of Lemma \ref{lemma A3} is complete.
\end{proof}

Let $\mathcal{H}_0=\text{span}\{e_k:k\in \cZ_0\}$ and  $D([0,\infty);\mathcal{H}_0)$ be the space of all $\rm c\grave{a}dl\grave{a}g$ functions taking values in $\mathcal{H}_0$. Keeping in mind that $d=|\cZ_0|<\infty$, it is well-known that, for any $w_0\in H$ and $g\in D([0,\infty);\mathcal{H}_0)$, there exists a unique solution $\Psi(w_0,g)\in C([0,\infty);H)\cap L^2_{loc}([0,\infty);V)$ to the following PDE:
    \begin{eqnarray*}
   \Psi(w_0,g)(t)
   &=&
   w_0+\nu \int_0^t\Delta (\Psi(w_0,g)(s)+g_s)\dif s
      \\ && ~ +\int_0^t B(\cK \Psi(w_0,g)(s) + \cK g_s,\Psi(w_0,g)(s)+g_s)\dif s.
  \end{eqnarray*}
Here $V=\{h\in H:\|h\|_1<\infty\}.$

{  We denote $\eta_t^\eps=Q(W_{\ell_t^\eps}-W_{\ell_0^\eps}),\eta_t=QW_{\ell_t},
v_t^\eps=\Psi(w_0,\eta^\eps)(t) $, and $ v_t=\Psi(w_0,\eta)(t)$.
}It is easy to see that $v_t+\eta_t$ is the unique solution $w_t$ to
(\ref{0.1}), i.e., $w_t=v_t+\eta_t$, and for any $\ell \in \mS$ and $\eps\in(0,1]$, $w^\eps_t:=v^\eps_t+\eta^\eps_t$ is the
solution of the following PDE:
  \begin{eqnarray*}
   w_t^{\eps}=w_0+\int_0^t \[ \nu \Delta w_s^{\eps} +B(\cK w_s^{\eps},w_s^{\eps})\] \dif s+Q(W_{\ell_t^\eps}-W_{\ell_0^\eps}).
  \end{eqnarray*}

 Recall $\mathbb{S}_1$ introduced in \eqref{eq S1}. We have
\begin{lemma} \label{lem A4}
For any $\ell\in \mathbb{S}_1$ and $\mathrm{w}\in  \mW$, the following statements hold:
  \begin{eqnarray}\label{eq lem A4 01}
  \begin{split}
    & \lim_{\eps \rightarrow 0}
      \int_0^{\sigma^\eps }  e^{-\nu (\sigma^\eps-s)+8\mathfrak{B}_0\kappa (\ell_{\sigma^\eps}^\eps-\ell_s^\eps) }  \| w_s^{\eps}\|_1^2
    \dif s
    \\ & =
    \int_0^{\sigma }  e^{-\nu (\sigma-s)+8\mathfrak{B}_0\kappa (\ell_{\sigma}-\ell_s) }  \| w_s\|_1^2
    \dif s,
    \end{split}
  \end{eqnarray}
and
     \begin{eqnarray}\label{eq lem A4 02}
\lim_{\eps \rightarrow 0}\|w_{\sigma^\eps}^{\eps}-w_{\sigma}\|^2=0.
  \end{eqnarray}
\end{lemma}

\begin{proof}
To prove this lemma, we first need some a priori estimates for $\Psi$.

By the chain rule and (\ref{44-1}), there exists a constant $C=C(\nu)>0$ such that, for any $w_0\in H$, $g\in D([0,\infty);\mathcal{H}_0)$ and $t\geq0$,
\begin{eqnarray*}
&&\|\Psi(w_0,g)(t)\|^2+\nu\int_0^t\|\Psi(w_0,g)(s)\|^2_1ds\\
&\leq&
\|w_0\|^2+\nu\int_0^t\|g_s\|_1^2\dif s+C\int_0^t\|g_s\|_2\|\Psi(w_0,g)(s)\|^2ds+C \int_0^t\|g_s\|_1^2\|\Psi(w_0,g)(s)\|ds.
\end{eqnarray*}
Applying the Gronwall lemma and using the fact that, for any $\alpha>0$, there exists a constant $C_\alpha$ such that $\|h\|_\alpha\leq C_\alpha\|h\|,\ \forall h\in \mathcal{H}_0$, there exists a constant $C>0$ such that, for any $T>0$,
\begin{eqnarray}\label{eq prio 01}
\begin{split}
& \sup_{t\in[0,T]}\|\Psi(w_0,g)(t)\|^2+\nu\int_0^T\|\Psi(w_0,g)(s)\|^2_1ds
\\
& { \leq
C\Big(\|w_0\|^2+\int_0^T(1+\|g_s\|^2)ds\Big)e^{C\int_0^T(1+\|g_s\|^2)ds}.}
\end{split}
\end{eqnarray}

 For any $g^1, g^2\in D([0,\infty);\mathcal{H}_0)$, put $\Psi^1(t)=\Psi(w_0,g^1)(t)$ and $\Psi^2(t)=\Psi(w_0,g^2)(t)$, simplifying the notation.
Using similar arguments as above,
\begin{eqnarray*}
&&\|\Psi^1(t)-\Psi^2(t)\|^2+\nu\int_0^t\|\Psi^1(s)-\Psi^2(s)\|^2_1ds\\
&\leq&
\nu\int_0^t\|g^1(s)-g^2(s)\|^2_1ds\\
   &&+
   2\int_0^t\langle B(\cK \Psi^1(s)+\cK g^1_s,\Psi^1(s)+g^1_s)-B(\cK \Psi^2(s)+\cK g^2_s,\Psi^2(s)+g^2_s), \Psi^1(s)-\Psi^2(s)\rangle ds\\
&=&
\nu\int_0^t\|g^1(s)-g^2(s)\|^2_1ds\\
   &&+
   2\int_0^t\langle B(\cK \Psi^1(s)+\cK g^1_s,\Psi^1(s)+g^1_s)-B(\cK \Psi^1(s)+\cK g^1_s,\Psi^2(s)+g^2_s), \Psi^1(s)-\Psi^2(s)\rangle ds\\
   &&+
   2\int_0^t\langle B(\cK \Psi^1(s)+\cK g^1_s,\Psi^2(s)+g^2_s)-B(\cK \Psi^2(s)+\cK g^2_s,\Psi^2(s)+g^2_s), \Psi^1(s)-\Psi^2(s)\rangle ds\\
   &\leq&
\nu\int_0^t\|g^1(s)-g^2(s)\|^2_1ds
   +
   C\int_0^t\|\Psi^1(s)-\Psi^2(s)\|\|\Psi^1(s)+g^1_s\|\|g^1(s)-g^2(s)\|_2ds\\
   &&+
   C\int_0^t\|\Psi^1(s)-\Psi^2(s)\|_1\|\Psi^2(s)+g^2_s\|_1\|\Psi^1(s)-\Psi^2(s)\| ds\\
   &&+
   C\int_0^t\|g^1(s)-g^2(s)\|_1\|\Psi^2(s)+g^2_s\|_1\|\Psi^1(s)-\Psi^2(s)\| ds\\
   &\leq&
C\Big(1+\sup_{s\in[0,t]}\|\Psi^1(s)+g^1_s\|^2\Big)\int_0^t\|g^1(s)-g^2(s)\|^2ds\\
   &&+
   C\int_0^t\Big(1+\|\Psi^2(s)+g^2_s\|^2_1\Big)\|\Psi^1(s)-\Psi^2(s)\|^2ds
   +
   \frac{\nu}{2}\int_0^t\|\Psi^1(s)-\Psi^2(s)\|^2_1ds.
\end{eqnarray*}
 Rearranging terms and using the Gronwall lemma, we arrive at, for any $T\geq0$,
 \begin{eqnarray}\label{eq prio 02}
 \begin{split}
&\sup_{t\in[0,T]}\|\Psi^1(t)-\Psi^2(t)\|^2+\frac{\nu}{2}\int_0^T\|\Psi^1(s)-\Psi^2(s)\|^2_1ds\\
   &\leq
C\Big(1+\sup_{s\in[0,T]}\|\Psi^1(s)+g^1_s\|^2\Big)
\int_0^T\|g^1(s)-g^2(s)\|^2ds\exp\Big\{C\int_0^T\big(
1+\|\Psi^2(s)+g^2_s\|^2_1\big)ds\Big\}.
\end{split}
\end{eqnarray}

For any $(\mathrm{w}, \ell) \in  \mW\times  \mS$, from the definitions of $\eta_t^\eps$ and $\eta_t$, it is easy to see that, for any $T\geq0$,
\begin{eqnarray}
\label{47-1}
\sup_{\eps\in(0,1]}\sup_{t\in[0,T]}\Big(\|\eta_t^\eps(\mathrm{w}, \ell)\|+\|\eta_t(\mathrm{w}, \ell)\|\Big)\leq { C} \sup_{t\in[0,\ell_{T+1}+T]}\|\mathrm{w}_t\|<\infty,
\end{eqnarray}
and
\begin{eqnarray*}
\lim_{\eps\rightarrow0}\int_0^T\|\eta_t^\eps(\mathrm{w}, \ell)-\eta_t(\mathrm{w}, \ell)\|^2dt=0.
\end{eqnarray*}
Combining the above two estimates with (\ref{eq prio 01}) and (\ref{eq prio 02}),
there exists a constant $C$ dependent on $\|w_0\|, T, \sup_{t\in[0,\ell_{T+1}+T]}\|\mathrm{w}_t\|$ such that
\begin{eqnarray}\label{eq 2023 1025 01}
\begin{split}
  & \sup_{\eps\in(0,1]}\Big(\sup_{t\in[0,T]}\|w^\eps_t\|^2+\int_0^T\|w^\eps_t\|^2_1dt\Big)(\mathrm{w}, \ell)
  \\ & \quad \quad +\Big(\sup_{t\in[0,T]}\|w_t\|^2+\int_0^T\|w_t\|^2_1dt\Big)(\mathrm{w}, \ell)\leq C,
  \end{split}
  \end{eqnarray}
  and
  \begin{eqnarray}\label{eq 2023 1025 02}
  \lim_{\eps\rightarrow0}\Big(\sup_{t\in[0,T]}\|v^\eps_t-v_t\|^2
  +\int_0^T\|w^\eps_t-w_t\|^2_1dt\Big)(\mathrm{w}, \ell)=0.
  \end{eqnarray}
Notice that, for any $\ell\in \mathbb{S}_1$ and $\mathrm{w}\in  \mW$,
\begin{eqnarray}\label{eq 2023 1025 03}
   \nonumber  && \|w_{\sigma^\eps}^{\eps}-w_{\sigma}\|\leq  \|v_{\sigma^\eps}^{\eps}-v_{\sigma}\|+\|\eta_{\sigma^\eps}^\eps-\eta_\sigma\|
   \\  \nonumber &&  \leq \|v_{\sigma^\eps}^{\eps}-v_{\sigma^\eps}\|+\|v_{\sigma^\eps}-v_{\sigma}\|
   +\|\eta_{\sigma^\eps}^\eps-\eta_\sigma\|
\\ && \leq
    \sup_{t\in[0,\sigma^1]}\|v_{t}^{\eps}-v_{t}\|
    +
    \|v_{\sigma^\eps}-v_{\sigma}\|+\|Q(W_{\ell_{\sigma^\eps}^\eps}-W_{\ell_0^\eps})-QW_{\ell_\sigma}\|.
\end{eqnarray}
Applying Lemmas \ref{qu-5} and \ref{lemma A3}, (\ref{eq 2023 1025 01})--(\ref{eq 2023 1025 03}), and the fact that $v_t$ is continuous
in $H$, we deduce (\ref{eq lem A4 01}) and (\ref{eq lem A4 02}), completing the proof of Lemma \ref{lem A4}.
\end{proof}

We also have the following estimate on $ w_t^{\eps}$.
  \begin{lemma}
  \label{qu-1}
  There exists a     positive constant  $ C$ which only depends on
  $\nu,\{b_j\}_{j\in \cZ_0},\nu_S$ and $d=|\cZ_0|$
 such that, for any $\kappa\in(0,\widetilde{\kappa}_0],\eps\in (0,1]$
 and $\ell\in \mS_1$(see \eqref{eq S1}),
  \begin{eqnarray}
\label{2626-1}
\begin{split}
 & \mE^{\mu_{\mW}} \exp\Big\{ \kappa \|w_{\sigma^\eps}^{\eps}\|-\kappa \| w_0 \|^2
   e^{-\nu \sigma^\eps  +8\mathfrak{B}_0\kappa \ell_{\sigma^\eps }^\eps}
   \\ & \quad\quad \quad \quad  +\nu
   \kappa  \int_0^{\sigma^\eps }  e^{-\nu (\sigma^\eps -s)+8\mathfrak{B}_0\kappa (\ell_{\sigma^\eps }^\eps-\ell_s^\eps) }  \| w_s^{\eps}\|_1^2
    \dif s
    \\ & \quad\quad \quad \quad - \kappa \mathfrak{B}_0 \int_{0}^{\sigma^\eps} e^{-\nu (\sigma^\eps  -s) +8\mathfrak{B}_0\kappa (\ell_{\sigma^\eps }^\eps-\ell_s^\eps)   }  \dif \ell_s^\eps   \Big\}
    \leq C.
    \end{split}
\end{eqnarray}
 Here $\widetilde{\kappa}_0$ is a constant appeared in Lemma \ref{Lemma A2}.
\end{lemma}
\begin{proof} Now we fix $\kappa\in(0,\widetilde{\kappa}_0],\eps\in (0,1]$
 and $\ell\in \mS_1$.

Let $\gamma^\eps$ be the inverse function of $\ell^\eps.$
  By a change of variable, for $t\geq \ell_{0}^{\eps},$
  $Y_t^{\eps}:=w_{\gamma_t^\eps}^{\eps},t\in [\ell_{0}^\eps,
 \infty)$ satisfies the following stochastic equation
  \begin{eqnarray}
  \label{100-3}
   Y_t^{\eps}= w_{0}+\int_{\ell_{0}^\eps}^t
    \[\nu \Delta Y_s^{\eps}+B(\cK Y_s^{\eps}, Y_s^{\eps})\]
   \dot   \gamma_s^\eps \dif s+Q (W_t-W_{\ell_{0}^\eps }).
  \end{eqnarray}

  By  It\^o's formula  we have
  \begin{eqnarray*}
   \dif \| Y_t^\eps\|^2
   = -2\nu  \| Y_t^\eps\|_1^2  \dot  \gamma_t^\eps\dif t+
  2 \langle  Y_t^\eps ,Q\dif W_t   \rangle
   +\mathfrak{B}_0 \dif t,
  \end{eqnarray*}
 and
\begin{eqnarray*}
   &&\dif \kappa \| Y_t^\eps\|^2e^{\nu \gamma_t^\eps-8\mathfrak{B}_0\kappa t   }\\
 &=& e^{\nu \gamma_t^\eps-8\mathfrak{B}_0\kappa t  } \[-2 \nu  \kappa
   \| Y_t^\eps\|_1^2  \dot  \gamma_t^\eps\dif t+
  2 \kappa  \langle  Y_t^\eps ,Q\dif W_t   \rangle
   +\kappa \mathfrak{B}_0\dif t\]
   \\ && \quad \quad  +  \kappa \| Y_t^\eps\|^2e^{\nu \gamma_t^\eps-8\mathfrak{B}_0\kappa t  }
   \big(\nu \dot \gamma_t^\eps-8\mathfrak{B}_0\kappa \big) \dif t
   \\ &\leq&
  - \nu
   \kappa  e^{\nu \gamma_t^\eps-8\mathfrak{B}_0\kappa t}  \| Y_t^\eps\|_1^2
   \dot  \gamma_t^\eps\dif t
   +\kappa \mathfrak{B}_0 e^{\nu \gamma_t^\eps-8\mathfrak{B}_0\kappa t  }  \dif t
+2 \kappa   e^{\nu \gamma_t^\eps-8\mathfrak{B}_0\kappa t   }
   \langle  Y_t^\eps ,Q\dif W_t   \rangle
   \\ &&- 8\mathfrak{B}_0\kappa^2  \| Y_t^\eps\|^2e^{\nu \gamma_t^\eps-8\mathfrak{B}_0\kappa t  } \dif t.
  \end{eqnarray*}
  Here we have used the inequality : $\|h\|_1\geq\|h\|,\ \forall h\in H$.
 Hence,
  \begin{eqnarray}
  \label{qu-3}
  \begin{split}
   &  \kappa \| Y_t^\eps\|^2  +\nu
   \kappa  \int_{\ell_0^\eps}^t  e^{-\nu (\gamma_t^\eps-\gamma_s^\eps)+8\mathfrak{B}_0\kappa (t-s) }  \| Y_s^{\eps}\|_1^2
   \dot  \gamma_s^\eps \dif s
   \\ & \leq   \kappa \| w_0 \|^2
   e^{-\nu\gamma_t^\eps+8\mathfrak{B}_0\kappa t}
   +\kappa \mathfrak{B}_0 \int_{\ell_0^\eps}^te^{-\nu (\gamma_t^\eps -\gamma_s^\eps) +8\mathfrak{B}_0\kappa (t-s)   }  \dif s + \tilde M_t,
   \end{split}
  \end{eqnarray}
where
  \begin{eqnarray*}
   && \tilde M_t =\tilde M_t^{\kappa,\eps}=2 \kappa  \int_{\ell_0^\eps}^{t}  e^{-\nu (\gamma_t^\eps-\gamma_s^\eps)+8\mathfrak{B}_0\kappa (t-s)   }
   \langle  Y^\eps_s ,Q\dif W_s    \rangle
  \\ && \quad\quad\quad\quad\quad\quad   -8\mathfrak{B}_0\kappa^2    \int_{\ell_0^\eps}^{t} \| Y^\eps_s\|^2e^{-\nu (\gamma_t^\eps-\gamma_s^\eps)+8\mathfrak{B}_0\kappa (t-s)   } \dif s.
  \end{eqnarray*}

 Next we prove that
\begin{eqnarray}
\label{2727-3}
  \mE^{\mu_{\mW}} \exp\{  \tilde M_{\ell^\eps_{\sigma^\eps} } \}\leq C.
\end{eqnarray}
Denote
\begin{align*}
  M_t &= 2\kappa  \int_{\ell_0^\eps}^t \langle  Y_s^\eps ,Q\dif W_s\rangle,
\quad \quad [M,M](s)= 4 \mathfrak{B}_0  \kappa^2  \int_{\ell_0^\eps}^t\|Y_s^\eps\|^2\dif s,
\\
N(s)&=M(s)-2 [M,M](s),
\quad\quad
  g(t,s)= e^{-\nu (\gamma_t^\eps-\gamma_s^\eps)+8\mathfrak{B}_0\kappa (t-s)   }.
\end{align*}
With these notations, one has $\tilde M_t=\int_{\ell_0^\eps}^t g(t,s)\dif N (s).$
 For any $K>0,$  by  the definition of  $\sigma^\eps$ and the fact that $ \nu \sigma^\eps -8\mathfrak{B}_0 \kappa \ell_{\sigma^\eps}^\eps = 1$ (see (1) of Lemma \ref{lemma A3}),
\begin{eqnarray}
\nonumber
&& \mP^{\mu_{\mW}} \Big(\tilde M_{\ell^\eps_{\sigma^\eps}} >  K \Big)=\mP^{\mu_{\mW}} \Big( \int_{\ell_0^\eps}^{\ell^\eps_{\sigma^\eps} }g(\ell^\eps_{\sigma^\eps},s)\dif N(s)>  K     \Big)
\\  \nonumber &&= \mP^{\mu_{\mW}} \Big(  e^{-\nu {\sigma^\eps} +8\mathfrak{B}_0\kappa \ell^\eps_{\sigma^\eps}  } \int_{\ell_0^\eps}^{\ell^\eps_{{\sigma^\eps}}}
e^{\nu \gamma_s^\eps-8\mathfrak{B}_0\kappa s}\dif N(s)> K     \Big)
\\  \nonumber &&=  \mP^{\mu_{\mW}} \Big(  \int_{\ell_0^\eps}^{\ell_{\sigma^\eps}^\eps}
e^{\nu \gamma_s^\eps-8\mathfrak{B}_0\kappa s}\dif N(s)>  eK    \Big)
\\  \nonumber &&=  \mP^{\mu_{\mW}} \Big( \int_{\ell_0^\eps}^{\ell^\eps_{{\sigma^\eps}}}
e^{\nu \gamma_s^\eps-8\mathfrak{B}_0\kappa s}\dif M(s)
\\  \nonumber && \quad\quad\quad\quad\quad\quad -\int_{\ell_0^\eps}^{\ell^\eps_{{\sigma^\eps}}}
e^{2 \nu \gamma_s^\eps - 16\mathfrak{B}_0\kappa s} 2  e^{- \nu \gamma_s^\eps+ 8\mathfrak{B}_0\kappa s} \dif  [M,M](s) >  eK       \Big)
\\  \nonumber &&\leq   \mP^{\mu_{\mW}} \Big(\int_{\ell_0^\eps}^{\ell^\eps_{{\sigma^\eps}}}
e^{\nu \gamma_s^\eps-8\mathfrak{B}_0\kappa s}\dif M(s)-\int_{\ell_0^\eps}^{\ell^\eps_{{\sigma^\eps}}}
e^{2 \nu \gamma_s^\eps-16\mathfrak{B}_0\kappa s} 2  e^{-1} \dif  [M,M](s)>eK    \Big)
\\  \label{2727-2} &&\leq \exp\big\{ -4  e^{-1 } eK    \big\}=e^{-4 K }.
\end{eqnarray}
In the first   inequality, we have used the fact:
 for any  $s=\ell_r^\eps (r\leq \sigma^\eps), - \nu \gamma_s^\eps+ 8\mathfrak{B}_0\kappa s=
 -\nu r+8\mathfrak{B}_0\kappa \ell_r^\eps \geq -1. $ {  For the last inequality we have used the following   fact (c.f. \cite[Theorem 5.2.9]{Dav-2009}):
 \begin{align*}
 & \mP^{\mu_{\mW}}\Big( \int_{\ell_0^\eps}^{\ell^\eps_{{\sigma^\eps}}}
e^{\nu \gamma_s^\eps-8\mathfrak{B}_0\kappa s}\dif M(s)- \frac{\alpha }{2}\int_{\ell_0^\eps}^{\ell^\eps_{{\sigma^\eps}}}
e^{2 \nu \gamma_s^\eps-16\mathfrak{B}_0\kappa s} \dif  [M,M](s) \geq \beta \Big)
\leq e^{-\alpha \beta}.
\end{align*}
}
The desired result   (\ref{2727-3}) follows immediately from  (\ref{2727-2}). The proof of (\ref{2727-3}) is complete.

Now replacing the $t$ in  (\ref{qu-3}) by $\ell^\eps_{\sigma^\eps},$   we obtain
\begin{eqnarray*}
   && \mE^{\mu_{\mW}}\Big[\exp\Big\{  \kappa \| Y^\eps_{\ell^\eps_{\sigma^\eps}}\|^2  + \nu
   \kappa  \int_{\ell_0^\eps}^{\ell^\eps_{\sigma^\eps}}  e^{-\nu ({\sigma^\eps}-\gamma_s^\eps)+ 8\mathfrak{B}_0\kappa (\ell^\eps_{\sigma^\eps}-s) }  \| Y_s^{\eps}\|_1^2
   \dot  \gamma_s^\eps \dif s\\
   &&\ \ \ \ \ \ \ -  \kappa \| w_0 \|^2
   e^{-\nu {\sigma^\eps}+8\mathfrak{B}_0\kappa \ell^\eps_{\sigma^\eps}}
   -\kappa \mathfrak{B}_0 \int_{\ell_0^\eps}^{\ell^\eps_{\sigma^\eps}}e^{-\nu ({\sigma^\eps} -\gamma_s^\eps) +8\mathfrak{B}_0\kappa (\ell^\eps_{\sigma^\eps}-s)   }  \dif s\Big\}\Big]
   \\ & & \leq   \mE^{\mu_{\mW}} \Big[\exp\Big\{\tilde M_{\ell^\eps_{\sigma^\eps}}\Big\}\Big].
\end{eqnarray*}
Combining the above inequality with $Y^\eps_{\ell^\eps_{\sigma^\eps}}=w_{\sigma^\eps}^\eps$,  $\gamma^\eps_s\big|_{s=\ell^\eps_r}=r$, (\ref{2727-3}), and
\begin{eqnarray*}
&& \int_{\ell_0^\eps}^{\ell^\eps_{\sigma^\eps}}  e^{-\nu ({\sigma^\eps}-\gamma_s^\eps)+8\mathfrak{B}_0\kappa (\ell^\eps_{\sigma^\eps}-s) }  \| Y_s^{\eps}\|_1^2
   \dot  \gamma_s^\eps \dif s
 = \int_{0}^{{\sigma^\eps}}  e^{-\nu ({\sigma^\eps}-r)+8\mathfrak{B}_0\kappa (\ell^\eps_{\sigma^\eps}-{\ell_r^\eps}) }  \|{w_r^\eps}\|_1^2
    \dif r,
   \\ &&  \int_{\ell_0^\eps}^{\ell^\eps_{\sigma^\eps}}e^{-\nu ({\sigma^\eps} -\gamma_s^\eps) +8\mathfrak{B}_0\kappa (\ell^\eps_{\sigma^\eps}-s)   }  \dif s
   =\int_{0}^{{\sigma^\eps}}e^{-\nu ({\sigma^\eps} -r) +8\mathfrak{B}_0\kappa (\ell^\eps_{\sigma^\eps}-{\ell_r^\eps})   }  \dif {\ell_r^\eps},
\end{eqnarray*}
we obtain the desired result (\ref{2626-1}).

The proof of Lemma \ref{qu-1} is complete.

\end{proof}

 \subsection{Proof of (\ref{2727-1})--(\ref{30-7}). }
 \label{A-1-3}

(\ref{2727-1}) was already proved  in Lemma  \ref{Lemma A2}, and it  obviously implies (\ref{2727-1-000}).
We now   prove
the  inequalities  (\ref{2626-2-000}) and (\ref{2626-2})   with $k=1$, and it is straightforward to extend it to the general case.
(\ref{2626-1}), $\ell_0=0,$ Fatou's Lemma and Lemmas \ref{lemma A3} and \ref{lem A4}
 imply   that
  \begin{eqnarray*}
  C &\geq &
  \mE^{\mu_{\mW}} \liminf_{\eps \rightarrow 0}  \exp\Big\{ \kappa \|w_{\sigma^\eps}^{\eps}\|^2-\kappa \| w_0 \|^2
   e^{-\nu \sigma^\eps  +8\mathfrak{B}_0\kappa \ell_{\sigma^\eps }^\eps}
   \\ &&  \quad\quad \quad \quad  +\nu
   \kappa  \int_0^{\sigma^\eps }  e^{-\nu (\sigma^\eps -s)+8\mathfrak{B}_0\kappa (\ell_{\sigma^\eps }^\eps-\ell_s^\eps) }  \| w_s^{\eps}\|_1^2
    \dif s
    \\ &&  \quad\quad \quad \quad -  \kappa \mathfrak{B}_0 \int_{0}^{\sigma^\eps} e^{-\nu (\sigma^\eps  -s) +8\mathfrak{B}_0\kappa (\ell_{\sigma^\eps }^\eps-\ell_s^\eps)   }  \dif \ell_s^\eps   \Big\}
    \\ &\geq &
    \mE^{\mu_{\mW}}   \exp\Big\{ \liminf_{\eps \rightarrow 0}  \kappa \|w_{\sigma^\eps}^{\eps}\|^2 -\lim_{\eps \rightarrow 0}  \kappa \| w_0 \|^2
   e^{-\nu \sigma^\eps  +8\mathfrak{B}_0\kappa \ell_{\sigma^\eps }^\eps}
   \\ &&  \quad\quad \quad \quad  +\liminf_{\eps \rightarrow 0}  \nu
   \kappa  \int_0^{\sigma^\eps }  e^{-\nu (\sigma^\eps -s)+8\mathfrak{B}_0\kappa (\ell_{\sigma^\eps }^\eps-\ell_s^\eps) }  \| w_s^{\eps}\|_1^2
    \dif s
    \\ &&  \quad\quad \quad \quad -\limsup_{\eps \rightarrow 0}  \kappa \mathfrak{B}_0 \int_{0}^{\sigma^\eps} e^{-\nu (\sigma^\eps  -s) +8\mathfrak{B}_0\kappa (\ell_{\sigma^\eps }^\eps-\ell_s^\eps)   }  \dif \ell_s^\eps   \Big\}
    \\ &\geq &
    \mE^{\mu_{\mW}}   \exp\Big\{  \kappa \|w_{\sigma}\|^2 - \kappa \| w_0 \|^2
   e^{-\nu \sigma  +8\mathfrak{B}_0\kappa \ell_{\sigma }}
   \\ &&  \quad\quad +
   \nu
   \kappa   \int_0^{\sigma }  e^{-\nu (\sigma -s)+8\mathfrak{B}_0\kappa (\ell_{\sigma}-\ell_s) }  \| w_s\|_1^2
    \dif s
   -   \kappa \mathfrak{B}_0 \ell_\sigma  \Big\}\\
       \\ &=&
    \mE^{\mu_{\mW}}   \exp\Big\{  \kappa \|w_{\sigma}\|^2 - \kappa \| w_0 \|^2
   e^{-1}
   \\ &&  \quad\quad +
   \nu
   \kappa   \int_0^{\sigma }  e^{-\nu (\sigma -s)+8\mathfrak{B}_0\kappa (\ell_{\sigma}-\ell_s) }  \| w_s\|_1^2
    \dif s
   -   \kappa \mathfrak{B}_0 \ell_\sigma  \Big\},
\end{eqnarray*}
which gives the desired result (\ref{2626-2-000}).

In the above inequality, taking expectation under the probability measure $\mP^{\mu_\mS}$, we get
       \begin{eqnarray*}
  \begin{split}
 & \mE   \Big[ \exp\Big\{ \kappa \|w_{\sigma}\|^2 - \kappa \| w_{0} \|^2
   e^{-1}
   \\ & \quad\quad \quad\quad  +
\nu \kappa   \int_{0}^{\sigma }  e^{-\nu (\sigma-s)+8\mathfrak{B}_0 \kappa (\ell_{{\sigma}}-\ell_s) }  \| w_s\|_1^2
    \dif s -  \kappa \mathfrak{B}_0 (\ell_{\sigma})      \Big\}  \Big]
  \leq C,
    \end{split}
\end{eqnarray*}
which is (\ref{2626-2}).

%

We consider now (\ref{100-2}) and  only    prove it for $k=0.$
Noticing that   $
   \nu \sigma-8\mathfrak{B}_0 \kappa  \ell_{\sigma} = 1
 $ (see (4) in Lemma \ref{lemma A3}), using  (\ref{2626-2}) with $k=1$, there exists a $\kappa_0>0$ such that
     \begin{eqnarray*}
 \mE  \exp\big\{ 2\kappa \|w_{\sigma}\|^2  -  2\kappa \mathfrak{B}_0  \ell_{\sigma}      \big\}
  \leq C\exp\{2\kappa e^{-1} \|w_{0}\|^2  \},~\forall  \kappa\in (0,\kappa_0].
\end{eqnarray*}
Thus, combining (\ref{2727-1-000}), we have
 \begin{eqnarray*}
    && \mE   \exp\big\{  \kappa   \|w_{\sigma}\|^2\big\}
     \leq    \big(\mE   \exp\big\{  2\kappa   \|w_{\sigma}\|^2-2\kappa \mathfrak{B}_0  \ell_{\sigma}  \big\} \big)^{1/2}
    \big(\mE   \exp\big\{ 2\kappa \mathfrak{B}_0  \ell_{\sigma}  \big\} \big)^{1/2}
    \\ && \leq  C\exp\{ \kappa e^{-1} \|w_{0}\|^2  \}
    \big(\mE   \exp\{ \nu \sigma\} \big)^{1/2}
   \leq  C_\kappa \exp\{ \kappa e^{-1} \|w_{0}\|^2  \}
 \end{eqnarray*}
 which yields the desired result (\ref{100-2}) for the case $k=0.$
%

 Following the arguments in the proof of \cite[(4.7)]{HM-2006},
also with the help of   (\ref{100-2}), we arrive at   (\ref{30-7}).


\subsection{Proof of (\ref{100-1}).}
\label{A-1-4}

  Recall $\mathbb{S}_1$ introduced in (\ref{eq S1}). Fix any $\ell\in\mathbb{S}_1$ and
 recall $Y_t^{\eps}$ introduced in (\ref{100-3}). By It\^o's formula,
  \begin{eqnarray}\label{007-1}
  &&\|Y_s^{\eps}\|^{2n}+2\nu n\int_{\ell_0^\eps}^s\|Y_u^{\eps}\|^{2n-2}\|Y_u^{\eps}\|^2_1\dot  \gamma_u^\eps du\\
    & =&\|w_0\|^{2n}+ 2n \int_{\ell_0^\eps}^s \|Y_u^{\eps}\|^{2n-2}\langle Y_u^{\eps},Q dW_u\rangle
    +
    n \int_{\ell_0^\eps}^s \|Y_u^{\eps}\|^{2n-2}\mathfrak{B}_0 du\nonumber\\
    &&+
    2n(n-1)\int_{\ell_0^\eps}^s\|Y_u^{\eps}\|^{2n-4}
    \sum_{j\in \cZ_0} \langle Y_u^{\eps},e_j\rangle^2 b_j^2 \dif u,\ \ s \geq \ell_0^\eps.\nonumber
  \end{eqnarray}
   Taking expectations with respect to $\mP^{\mu_\mW}$ yields
   \begin{eqnarray*}
  && \mE^{\mu_\mW }\|Y_s^{\eps}\|^{2n}
    \leq\|w_0\|^{2n}
    +
    C_n \int_{\ell_0^\eps}^s \mE^{\mu_\mW }\|Y_u^{\eps}\|^{2n-2} du
   \\ &&  \leq
    \|w_0\|^{2n}
    +
    \frac{1}{2}\sup_{u\in[\ell_0^\eps,t]}\mE^{\mu_\mW }\|Y_u^{\eps}\|^{2n}
    +
    C_nt^n, \quad \forall s\in [\ell_0^\eps,t].
  \end{eqnarray*}
 Rearranging terms, we arrive at, for any $t\geq \ell_0^\eps$,
    \begin{eqnarray*}
  \sup_{s\in[\ell_0^\eps,t]}\mE^{\mu_\mW }\|Y_t^{\eps}\|^{2n}
    \leq 2\|w_0\|^{2n}
    +
    C_nt^n.
  \end{eqnarray*}
  Using (\ref{007-1}) and the BDG inequality, we have
    \begin{eqnarray*}
  &&\mE^{\mu_\mW }\Big(\sup_{s\in[\ell_0^\eps,t]}\|Y_s^{\eps}\|^{2n}\Big)\\
    &\leq&\|w_0\|^{2n}+ 2n \mE^{\mu_\mW }\sup_{s\in[\ell_0^\eps,t]}\Big|\int_{\ell_0^\eps}^s \|Y_r^{\eps}\|^{2n-2}\langle Y_r^{\eps},Q dW_r\rangle\Big|
    +
    C_n \mE^{\mu_\mW }\int_{\ell_0^\eps}^t \|Y_r^{\eps}\|^{2n-2}dr\\
    &\leq&
    \|w_0\|^{2n}+ C_n \Big(\mE^{\mu_\mW }\sup_{s\in[\ell_0^\eps,t]}\Big|\int_{\ell_0^\eps}^s \|Y_r^{\eps}\|^{2n-2}\langle Y_r^{\eps},Q dW_r\rangle\Big|^2\Big)^{1/2}
    +
    C_n \mE^{\mu_\mW }\int_{\ell_0^\eps}^t \|Y_r^{\eps}\|^{2n-2}dr\\
    &\leq&
    \|w_0\|^{2n}+ C_n \Big(\mE^{\mu_\mW }\int_{\ell_0^\eps}^t \|Y_r^{\eps}\|^{4n-2}dr\Big)^{1/2}
    +
    C_n \mE^{\mu_\mW }\int_{\ell_0^\eps}^t \|Y_r^{\eps}\|^{2n-2}dr\\
   &\leq&   {
    \|w_0\|^{2n}+ C_n \big(\|w_0\|^{4n-2}t+t^{2n}\big)^{1/2}
    +
    C_n \big(\|w_0\|^{2n-2}t+t^{n}\big)
    }
    \\
   &\leq&  {
    C_n(1+t)\|w_0\|^{2n}
    +
    C_n(1+t^{n}).
    }
  \end{eqnarray*}
Using the fact that $w_{t}^{\eps}=Y_{\ell_t^\eps}^{\eps}$, we arrive at
  \begin{eqnarray*}
     \mE^{\mu_{\mW}} \sup_{s\in [0,t]}\|w_s^\eps\|^{2n}
    =
    \mE^{\mu_\mW }\sup_{s\in[\ell_0^\eps,{\ell_t^\eps}]}\|Y_s^{\eps}\|^{2n}
    \leq C_n(1+\ell_t^\eps)  \|w_0\|^{2n}+
    C_n\big(1+(\ell_t^\eps)^{n}\big).
  \end{eqnarray*}
  {
  By \eqref{eq 2023 1025 02}, Fatou's lemma,  (ii) of Lemma \ref{qu-5}
  and
  $$w_s=(w_s-\eta_s)-(w_s^\eps-\eta_s^\eps)+w_s^\eps+(\eta_s-\eta_s^\eps)=v_s-v_s^\eps+w_s^\eps
  +(\eta_s-\eta_s^\eps),\forall \eps>0,$$ one has
\begin{eqnarray*}
  && \mE^{\mu_{\mW}} \sup_{s\in [0,t]}\|w_s\|^{2n}
 \\ && \leq
   C_n \mE^{\mu_{\mW}}\Big[ \liminf_{\eps\rightarrow0} \sup_{s\in [0,t]}\|w_s^\eps\|^{2n}+\liminf_{\eps\rightarrow0}
    \sup_{s\in [0,t]}\|v_s-v_s^\eps\|^{2n}+\liminf_{\eps\rightarrow0}\sup_{s\in [0,t]}|\eta_s-\eta_s^\eps|^{2n}\Big]
   \\
   & &\leq
   C_n \liminf_{\eps\rightarrow0}\mE^{\mu_{\mW}} \sup_{s\in [0,t]}\|w_s^\eps\|^{2n}
   + C_n\liminf_{\eps\rightarrow0} \mE^{\mu_{\mW}} \sup_{s\in [0,t]}\big(|\eta_s|^{2n}+|\eta_s^\eps|^{2n}\big)
\\ && \leq  C_n(1+\ell_t)  \|w_0\|^{2n}+
    C_n\big(1+(\ell_t)^{n}\big).
\end{eqnarray*}
}
Since the above estimate holds for any fixed $\ell\in \mS_1$, ${\mu_\mS}(\mathbb{S}_1)=1$ (see Lemma \ref{lemma A3}) and since $\sigma$ only depends on $\ell\in \mS$, we first replace the $t$ by $\sigma$, and then take expectations with respect to $\mu_\mS$, to obtain
\begin{eqnarray*}
  \mE \sup_{s\in [0,\sigma]}\|w_s\|^{2n}
&\leq&  C_n(1+\mE^{\mu_\mS}\ell_\sigma)  \|w_0\|^{2n}+
    C_n\big(1+\mE^{\mu_\mS}(\ell_\sigma)^{n}\big)\\
&\leq&C_n\Big(1+\mE^{\mu_\mS}\Big(\frac{\nu \sigma }{8\kappa \mathfrak{B}_0 }\Big)\Big)  \|w_0\|^{2n}+
    C_n\big(1+\mE^{\mu_\mS}\Big(\frac{\nu \sigma }{8\kappa \mathfrak{B}_0 }\Big)^{n}\big)\\
&\leq&
C_{n,\kappa}(1+ \|w_0\|^{2n}).
\end{eqnarray*}
Here we have used the fact that $\ell_{\sigma}\leq \frac{\nu \sigma }{8\kappa \mathfrak{B}_0 }$(see (4) of Lemma \ref{lemma A3})
and (\ref{2727-1}), $C_{n,\kappa}$ is a constant depending on $n,\kappa$ and $\nu, \mathfrak{B}_0.$

The proof of (\ref{100-1}) is complete.

%

\section{Proof of Proposition  \ref{3-8}}
\label{180-1}

We will prove Proposition  \ref{3-8} by contradiction.

Suppose that   Proposition  \ref{3-8} { were}  not true,  then   there exist  sequences $\{w_0^{(k)} \} \subseteq  B_H(\mathfrak{R}), \{\eps_k\}\subseteq (0,1)$  and a positive number $\delta_0$  such that
\begin{eqnarray}
\label{D-1}
\lim_{k\rightarrow \infty} \mP(  X^{w_0^{(k)},\alpha,N}<\eps_k)\geq \delta_0>0  \text{ and  } \lim_{k\rightarrow \infty}\eps_k=0.
\end{eqnarray}
Our  aim is to find something
which contradicts (\ref{D-1}).

Since $H$ is a Hilbert space, there exists a subsequence $\{w_0^{(n_k)},k\geq 1\}$ of  $\{w_0^{(k)},k\geq 1\}$
such that    $w_0^{(n_k)}$ converges weakly to  some  $ w_0^{(0)}\in H.$
We still denote this subsequence  by $\{w_0^{(k)},k\geq 1\}.$
Let $w_t^{(k)}$ denote the solution of equation (\ref{0.1}) with $w_t|_{t=0}=w_0^{(k)}(k\geq 0)$.
In the following  equation
    \begin{align}
    \label{9-2}
	\partial_t   J_{s,t}\xi-  \nu \Delta J_{s,t}\xi -\tilde{B}(w_t,J_{s,t}\xi)&=0,
  \\  J_{s,s}\xi&= \xi, \nonumber
\end{align}
when $w_t$ is replaced by $w_t^{(k)}$, we denote its solution by $J_{s,t}^{(k)}\xi$.
We denote  the adjoint of  $J_{s,t}^{(k)}$ by $K_{s,t}^{(k)}.$

Recall that  for any $M\in \mN,$ ${H}_M=\{e_j: j\in \mZ_*^2  \text{ and } |j|\leq M\},$
${P}_M$  denotes the orthogonal projections from $H$ onto ${H}_M$
and   ${Q}_M u:=u-{P}_Mu, \forall u\in H.$
As before,  $C$  denotes   a constant depending
     $ \nu,\{b_j\}_{j\in \cZ_0},\nu_S,d,$
   $C_\mathfrak{R}$  denotes  a constant depending on $\mathfrak{R}$
   and   $ \nu,\{b_j\}_{j\in \cZ_0},\nu_S,d.$
   The values of the constants may change from line to line.

 In what follows, we will  give some  estimates  for $\|P_M w_t^{(k)}-P_Mw_t^{(0)}\|$ and $\|Q_Mw_t^{(k)}\|$  in Lemma  \ref{E-5},    and  for $\|J_{s,t}^{(k)}\xi-J_{s,t}^{(0)}\xi\|, s,t\in (0,\sigma]$ in  Lemma \ref{9-3}.  
 We finish  the   proof   of  Proposition  \ref{3-8} { at the end of the section.}

\begin{lemma}
\label{E-5}
For any $t\geq 0$, $k\in \mN$ and $M>\max\{|j|: j\in \cZ_0\},$ one has
  \begin{eqnarray}
  \label{E-1}
  \begin{split}
    \|Q_Mw_t^{(k)}\|^2& \leq e^{-\nu M^2 t }\|Q_Mw_0^{(k)}\|^2+
   \\ &
 \frac{C }{M^{1/2} }\big(\int_0^t \|w_s^{(k)}\|_1^{4/3}\dif s\big)^{3/4} \sup_{s\in [0,t]} \|w_s^{(k)}\|^{3}
 \end{split}
  \end{eqnarray}
  and
     \begin{eqnarray}
     \label{E-2}
     \begin{split}
      & \|P_M w_t^{(k)}-P_M w_t^{(0)} \|^2\leq
 C_\mathfrak{R}   \exp\big\{C \int_0^t (\|w_s^{(k)}\|_1^{4/3}+\|w_s^{(0)}\|_1^{4/3})\dif s\big\}
 \\ & \quad \quad \times \sup_{r\in [0,t]} (1+\|w_r^{(k)}\|^4+\|w_r^{(0)}\|^4 )
  \\ & \quad \quad \times \Big[ \|P_M w_0^{(k)}-P_M w_0^{(0)}\|^2+\frac{1+ t}{M^{1/2}} \Big].
     \end{split}
   \end{eqnarray}

\end{lemma}
\begin{proof}
  First, we give a proof of  (\ref{E-1}).
   By (\ref{44-1}),
  one has
  \begin{eqnarray*}
\langle  B (\cK w_t^{(k)},w_t^{(k)}),Q_M w_t^{(k)}\rangle & \leq &
C\|Q_M  w_t^{(k)}\|_1\|w_t^{(k)}\|_{1/2}\|w_t^{(k)}\|
\\ & \leq & \frac{\nu}{4} \|Q_M w_t^{(k)}\|_1^2+C \|w_t^{(k)}\|_1 \|w_t^{(k)}\|^3.
  \end{eqnarray*}
Thus, in view of the equation (\ref{0.1}),  we obtain
\begin{eqnarray*}
 \dif \|Q_Mw_t^{(k)}\|^2 &=&  -2 \nu \|Q_Mw_t^{(k)}\|_1^2+ \langle  B (\cK w_t^{(k)},w_t^{(k)}),Q_M w_t^{(k)}\rangle
 \\ &\leq & -\nu M^2  \|Q_M w_t^{(k)}\|^2\dif t+C \|w_t^{(k)}\|_1 \|w_t^{(k)}\|^3\dif t.
\end{eqnarray*}
 It follows that
\begin{eqnarray*}
&& \|Q_Mw_t^{(k)}\|^2  \leq  e^{-\nu M^2 t} \|Q_Mw_0^{(k)}\|^2
+C \int_0^t e^{-\nu M^2 (t-s) }  \|w_s^{(k)}\|_1 \|w_s^{(k)}\|^3\dif s
\\ & & \leq   e^{-\nu M^2 t} \|Q_Mw_0^{(k)}\|^2
\\ &&\quad +
C \big(\int_0^t e^{-4\nu M^2 (t-s) }  \dif s \big)^{1/4}
\big(\int_0^t \|w_s^{(k)}\|_1^{4/3}\dif s\big)^{3/4} \sup_{s\in [0,t]} \|w_s^{(k)}\|^{3}.
\end{eqnarray*}
This completes the   proof of  (\ref{E-1}).

Next, we will prove(\ref{E-2}).
One easily sees that
\begin{eqnarray}
\label{E-3}
\begin{split}
& \dif  \|P_M w_t^{(k)}-P_M w_t^{(0)}\|^2
\\ & = -\nu \|P_M w_t^{(k)}-P_M w_t^{(0)}\|_1^2\dif t
  \\ &\quad  +  \langle P_M w_t^{(k)}-P_M w_t^{(0)},
      B(\cK w_t^{(k)},w_t^{(k)})- B(\cK w_t^{(0)},w_t^{(0)})   \rangle.
      \end{split}
\end{eqnarray}
Clearly, we have
   \begin{eqnarray*}
     && \langle P_M w_t^{(k)}-P_M w_t^{(0)},
      B(\cK w_t^{(k)},w_t^{(k)})- B(\cK w_t^{(0)},w_t^{(0)})   \rangle
      \\ && = \langle P_M w_t^{(k)}-P_M w_t^{(0)},
      B(\cK w_t^{(k)}-\cK w_t^{(0)},w_t^{(k)})   \rangle
     \\ && \quad  + \langle P_M w_t^{(k)}-P_M w_t^{(0)},
    B(\cK w_t^{(0)},  w_t^{(k)}-w_t^{(0)})   \rangle
    \\ && :=I_1+I_2.
   \end{eqnarray*}
   For the term $I_1$, we have
   \begin{eqnarray*}
     I_1 &\leq &  C \|P_M w_t^{(k)}-P_M w_t^{(0)} \|_{1/2}
     \| P_M w_t^{(k)}- P_Mw_t^{(0)} \| \|w_t^{(k)}\|_1
     \\ && \quad \quad +C \|P_M w_t^{(k)}-P_M w_t^{(0)} \|_{1}
     \| Q_M w_t^{(k)}- Q_M w_t^{(0)} \|  \|w_t^{(k)}\|_{1/2}
     \\ &\leq &\frac{\nu}{6}\|P_M w_t^{(k)}-P_M w_t^{(0)} \|_{1}^2+C
     \|P_M w_t^{(k)}-P_M w_t^{(0)} \|^2 \| \|w_t^{(k)}\|_1^{4/3}
     \\ && \quad \quad +C \| Q_M w_t^{(k)}- Q_M w_t^{(0)} \|^2
     \|w_t^{(k)}\|^2_{1/2}.
   \end{eqnarray*}
    Obviously,
   \begin{eqnarray*}
     I_2& \leq&
      C\|P_M w_t^{(k)}-P_M w_t^{(0)} \|_{1}
     \| Q_M w_t^{(k)}- Q_M w_t^{(0)} \| \|w_t^{(0)}\|_{1/2}
     \\ &\leq & \frac{\nu}{6}\|P_M w_t^{(k)}-P_M w_t^{(0)} \|_{1}^2+C
     \| Q_M w_t^{(k)}- Q_M w_t^{(0)} \|^2 \|w_t^{(0)}\|_{1/2}^2
   \end{eqnarray*}
 Combining the estimates of $I_1$,$I_2$ with  (\ref{E-3}),(\ref{E-1}),
 taking into account the fact that  $\|Q_M w_0^{(k)}\|\leq \mathfrak{R},$ we obtain
 \begin{eqnarray*}
  && \|P_M w_t^{(k)}-P_M w_t^{(0)}\|^2 \leq
  \|P_M w_0^{(k)}-P_M w_0^{(0)}\|^2 e^{C\int_0^t \|w_s^{(k)}\|_1^{4/3} \dif s  }
  \\ &&\quad\quad +C  e^{C\int_0^t \|w_s^{(k)}\|_1^{4/3}\dif s}   \int_0^t \| Q_M w_s^{(k)}- Q_M w_s^{(0)} \|^2 (\|w_s^{(k)}\|_{1/2}^2 +\|w_s^{(0)}\|_{1/2}^2)\dif s
  \\
  && \leq
  \|P_M w_0^{(k)}-P_M w_0^{(0)}\|^2 e^{C\int_0^t \|w_s^{(k)}\|_1^{4/3} \dif s  }
  \\ &&\quad+ C e^{C\int_0^t \|w_s^{(k)}\|_1^{4/3}\dif s}   \big(\int_0^t \| Q_M w_s^{(k)}- Q_M w_s^{(0)} \|^8 \dif s\big)^{1/4}
  \\ && \quad \times  \Big(\big(\int_0^t \|w_s^{(k)}\|_{1}^{4/3}\dif s\big)^{3/4} +\big(\int_0^t \|w_s^{(0)}\|_{1}^{4/3})\dif s\big)^{3/4}\Big)\sup_{s\in [0,t]}(\|w_s^{(k)}\|+\|w_s^{(0)}\| )
  \\ &&\leq  e^{C\int_0^t \|w_s^{(k)}\|_1^{4/3}\dif s }\|P_M w_0^{(k)}-P_M w_0^{(0)}\|^2
  \\ && \quad+C_\mathfrak{R}  e^{C\int_0^t \|w_s^{(k)}\|_1^{4/3}\dif s+C\int_0^t
  \|w_s^{(0)}\|_1^{4/3}\dif s} \sup_{s\in [0,t]}(\|w_s^{(k)}\|+\|w_s^{(0)}\| )
  \\ &&  \times \Bigg(\int_0^t\Big[
  e^{-4\nu M^2 s }
  +\frac{1 }{M^2  }\big(\int_0^s \|w_r^{(k)}\|_1^{4/3}\dif r\big)^{3}  \sup_{s\in [0,t]} \|w_s^{(k)}\|^{12}
  \\ && \quad\quad \quad\quad \quad\quad \quad\quad  \quad\quad  +\frac{1 }{M^2  }\big(\int_0^s \|w_r^{(0)}\|_1^{4/3}\dif r\big)^{3}  \sup_{s\in [0,t]} \|w_s^{(0)}\|^{12}
  \Big]  \dif s\Bigg)^{1/4},
 \end{eqnarray*}
 which implies the desired result (\ref{E-2}).
  \end{proof}

  \begin{lemma}
  \label{9-3}
For any $0\leq s\leq t,$  $k\in \mN$ and $\xi\in H$ with $\|\xi\|=1,$ one has
   \begin{eqnarray*}
    && \|J_{s,t}^{(k)}\xi-J_{s,t}^{(0)}\xi\|^2
    \leq
   C\sup_{r\in [s,t]}\|w_r^{(k)}-w_r^{(0)}\|^2 \cdot  e^{C  \int_{s}^t
   (\|w_r^{(k)}\|_1^{4/3}
    +
    \|w_r^{(0)}\|_1^{4/3})\dif r},
    \end{eqnarray*}
    where $C$ is a  constant depending on  $ \nu,\{b_j\}_{j\in \cZ_0},\nu_S,d.$
  \end{lemma}
  \begin{proof}
    By the equation (\ref{9-2}), one easily  sees that
    \begin{eqnarray}
      \nonumber && \dif \|J_{s,t}^{(k)}\xi-J_{s,t}^{(0)}\xi\|^2
    \leq   -
      \nu\|J_{s,t}^{(k)}\xi-J_{s,t}^{(0)}\xi\|^2_1\dif t
      \\  \nonumber &&\quad +\langle J_{s,t}^{(k)}\xi-J_{s,t}^{(0)}\xi,
      {B}(\cK w_t^{(k)},J_{s,t}^{(k)}\xi)
      -{B}(\cK w_t^{(0)},J_{s,t}^{(0)}\xi) \rangle \dif t
      \\ \nonumber &&\quad+  \langle J_{s,t}^{(k)}\xi-J_{s,t}^{(0)}\xi,
      {B}(\cK J_{s,t}^{(k)}\xi, w_t^{(k)})
      -{B}(\cK J_{s,t}^{(0)}\xi, w_t^{(0)}) \rangle \dif t
      \\ \label{E-4} &&:=I_1(t)\dif t+I_2(t)\dif t+I_3(t)\dif t.
    \end{eqnarray}
 For the terms $I_2(t),I_3(t)$, we have
    \begin{eqnarray*}
      I_2(t)&= &
       \langle J_{s,t}^{(k)}\xi-J_{s,t}^{(0)}\xi,
      {B}(\cK w_t^{(k)}-\cK w_t^{(0)},J_{s,t}^{(k)}\xi)
       \rangle
       \\ &\leq &
      C  \|w_t^{(k)}- w_t^{(0)}\| \| J_{s,t}^{(k)}\xi-J_{s,t}^{(0)}\xi\|_1 \|J_{s,t}^{(k)}\xi\|_{1/2}
    \\ &\leq & C  \|w_t^{(k)}- w_t^{(0)}\|^2  \|J_{s,t}^{(k)}\xi\|_{1/2}^2+\frac{\nu}{6}
    \|J_{s,t}^{(k)}\xi-J_{s,t}^{(0)}\xi\|_1^2
    \end{eqnarray*}
    and
    \begin{eqnarray*}
      I_3(t) &= &
      \langle J_{s,t}^{(k)}\xi-J_{s,t}^{(0)}\xi,
      {B}(\cK J_{s,t}^{(k)}\xi-\cK J_{s,t}^{(0)}\xi, w_t^{(k)})
      +{B}(\cK J_{s,t}^{(0)}\xi, w_t^{(k)}- w_t^{(0)}) \rangle
      \\ &\leq & C \|w_t^{(k)}\|_1 \| J_{s,t}^{(k)}\xi- J_{s,t}^{(0)}\xi\|^{3/2}
      \| J_{s,t}^{(k)}\xi- J_{s,t}^{(0)}\xi\|_1^{1/2}
      \\ && +C\|w_t^{(k)} - w_t^{(0)}\|\| J_{s,t}^{(k)}\xi- J_{s,t}^{(0)}\xi\|_1\|J_{s,t}^{(0)}\xi \|_{1/2}
    \\ &\leq & \frac{\nu}{6}\| J_{s,t}^{(k)}\xi- J_{s,t}^{(0)}\xi\|_1^{2}
    +C \| J_{s,t}^{(k)}\xi- J_{s,t}^{(0)}\xi\|^{2} \|w_t^{(k)}\|_1^{4/3}
    \\ && +C \|w_t^{(k)}- w_t^{(0)}\|^2\|J_{s,t}^{(0)}\xi\|_{1/2}^2.
    \end{eqnarray*}
  Combining the above estimates of $I_2,I_3$ with (\ref{E-4}),   we obtain
    \begin{eqnarray*}
    && \|J_{s,t}^{(k)}\xi-J_{s,t}^{(0)}\xi\|^2
    \\ && \leq
   C \sup_{r\in [s,t]}\|w_r^{(k)}-w_r^{(0)}\|^2 \cdot  e^{C \int_{s}^t
   \|w_r^{(k)}\|_1^{4/3}
    \dif r}
    \int_s^t\[ \|J_{s,r}^{(k)}\xi\|_{1/2}^2+\|J_{s,r}^{(0)}\xi\|_{1/2}^2 \] \dif r.
    \end{eqnarray*}
By Lemma \ref{15-4}, H\"older's inequality and using the fact that $\|a\|_{1/2}^2\leq \|a\|\|a\|_1$, we obtain the desired result.

\end{proof}


\textbf{Now we are in a position to complete the   proof of   Proposition  \ref{3-8}.}

With the help of Lemmas   \ref{E-5}, \ref{9-3},  for any  $\kappa>0, r\in [\frac{\sigma}{2},\sigma]$ and  $\phi \in H$ with $\|\phi\|=1$, we have
\begin{eqnarray}
  \nonumber  && \|J_{r,\sigma}^{(k)}\phi-J_{r,\sigma}^{(0)}\phi \|\leq
  C e^{C \int_{\sigma/2}^\sigma
   (\|w_r^{(k)}\|_1^{4/3}
    +
    \|w_r^{(0)}\|_1^{4/3})\dif r} \cdot \sup_{r\in [\frac{\sigma}{2},\sigma]}\|w_r^{(k)}-w_r^{(0)}\|^2
    \\
       \nonumber  &&\leq C_\mathfrak{R}   \exp\big\{C  \int_0^\sigma (\|w_s^{(k)}\|_1^{4/3}+\|w_s^{(0)}\|_1^{4/3})\dif s\big\} \sup_{r\in [0,\sigma]} (1+\|w_r^{(k)}\|^4+\|w_r^{(0)}\|^4 )
        \\
           \nonumber  && \quad \times  \[e^{-\nu M^2\sigma/2}+\|P_M w_0^{(k)}-P_M w_0^{(0)}\|^2+\frac{1+ \sigma}{M^{1/2}} \]
   \\
   \nonumber  &&\leq {\Big[ C_\mathfrak{R}  \exp\big\{\frac{\nu \kappa}{6}  \int_0^\sigma (\|w_s^{(k)}\|_1^{2}+\|w_s^{(0)}\|_1^{2})e^{-\nu (\sigma-s)+8\mathfrak{B}_0 \kappa (\ell_{{\sigma}}-\ell_s) }\dif s\big\}}
   \\   \nonumber  && \quad\quad\quad\quad {\cdot  \sup_{r\in [0,\sigma]} (1+\|w_r^{(k)}\|^4+\|w_r^{(0)}\|^4 )\Big]}
   \\   \nonumber  && \quad \times{ \exp\{C_\kappa \int_0^\sigma e^{ 2\nu (\sigma-s)-16\mathfrak{B}_0 \kappa (\ell_{{\sigma}}-\ell_s) }\dif s  \}}
  \\   \nonumber  && \quad \times  \[e^{-\nu M^2\sigma/2}+\|P_M w_0^{(k)}-P_M w_0^{(0)}\|^2+\frac{1+ \sigma}{M^{1/2}} \]
 \\ \label{32-1} &&:={\Xi(k)} \cdot {\tilde  X} \cdot  \Upsilon(k,M), \quad \forall r\in [\frac{\sigma}{2},\sigma].
\end{eqnarray}
Note that
\begin{align*}
& \langle K_{r,\sigma}^{(k)}\phi,
 e_j\rangle^2 =\Big(\langle K_{r,\sigma}^{(0)}\phi,
 e_j\rangle  +\langle K_{r,\sigma}^{(k)}\phi-K_{r,\sigma}^{(0)}\phi,
 e_j\rangle\Big)^2
 \\ & \geq \frac{1}{2}\langle K_{r,\sigma}^{(0)}\phi,
 e_j\rangle^2-3\langle K_{r,\sigma}^{(k)}\phi-K_{r,\sigma}^{(0)}\phi,
 e_j\rangle^2  \geq
 \frac{1}{2}\langle K_{r,\sigma}^{(0)}\phi,
 e_j\rangle^2-3\| K_{r,\sigma}^{(k)}\phi-K_{r,\sigma}^{(0)}\phi\|^2,
 \end{align*}
and recall that $K_{r,\sigma}$ is the adjoint of $J_{r,\sigma}$.
It follows from (\ref{32-1})  that
\begin{eqnarray*}
 && \mP\left(\inf_{\phi\in \cS_{\alpha,N}}\sum_{j\in \cZ_0}\int_{0}^\sigma \langle K_{r,\sigma}^{(k)}\phi,
 e_j\rangle^2\dif S_r<\eps_k\right)
 \\ &&\leq \mP\Big(\frac{1}{2}\inf_{\phi\in \cS_{\alpha,N} }\sum_{j\in \cZ_0}\int_{\sigma/2}^\sigma \langle K_{r,\sigma}^{(0)}\phi,
 e_j\rangle^2\dif S_r<\eps_k
 \\ && \quad\quad\quad\quad\quad\quad   +3 d \sup_{\phi \in S_{\alpha,N}}\sup_{r\in [\sigma/2,\sigma]}\|K_{r,\sigma}^{(k)}\phi-K_{r,\sigma}^{(0)}\phi\|^2 S_\sigma \Big)
 \\ &&\leq   \mP\left(\inf_{\phi\in \cS_{\alpha,N}}\sum_{j\in \cZ_0}\int_{\sigma/2}^\sigma \langle
  K_{r,\sigma}^{(0)}\phi,
 e_j\rangle^2\dif S_r<2 \eps_k+
   \Xi(k)\cdot (6d S_\sigma \tilde X)\cdot     \Upsilon(k,M)\right)
 \end{eqnarray*}
 Therefore, for any $\mathcal{C}>0$, we deduce that
 \begin{eqnarray}
 \nonumber  && \mP\left(\inf_{\phi\in \cS_{\alpha,N}}\sum_{j\in \cZ_0}\int_{0}^\sigma \langle K_{r,\sigma}^{(k)}\phi,
 e_j\rangle^2\dif S_r<\eps_k\right)
 \\  \nonumber  &&\leq
 \mP\left(\inf_{\phi\in \cS_{\alpha,N}}\sum_{j\in \cZ_0} \int_{\sigma/2}^\sigma \langle K_{r,\sigma}^{(0)}\phi,
 e_j\rangle^2\dif S_r<2 \eps_k+
  \mathcal{C}^2  \Upsilon(k,M)\right)
  \\  \label{2032}  && \quad \quad\quad\quad \quad \quad\quad\quad  +\mP(  \Xi(k) \geq \mathcal{C}    )+ \mP(6d S_\sigma \tilde X> \mathcal{C} )
\end{eqnarray}
Letting $k\rightarrow \infty$ in (\ref{2032}), by (\ref{D-1}),  one sees that
\begin{eqnarray}
\nonumber
  && \delta_0\leq \mP\left(\inf_{\phi\in\cS_{\alpha,N}}\sum_{j\in \cZ_0} \int_{\sigma/2}^\sigma \langle K_{r,\sigma}^{(0)}\phi,
 e_j\rangle^2\dif S_r \leq
\mathcal{C}^2    (e^{-\nu M^2\sigma/2}+\frac{1+\sigma}{M^{1/2}})\right) \\  \label{B-1} && \quad \quad\quad\quad \quad \quad  +\frac{\sup_{k\geq 0} \mE \Xi(k)  }{\mathcal{C}}+ \mP(6d S_\sigma \tilde X> \mathcal{C} )
\end{eqnarray}
By Lemma \ref{qu-2} and the fact that $\sup_{k\in \mN}\|w_0^{(k)}\|\leq \mathfrak{R}$,  for any $\kappa\in (0,\kappa_0]$, one has  $\sup_{k}\mE \Xi(k)<\infty.$
In (\ref{B-1}),   first letting $M\rightarrow \infty$ and then letting $\mathcal{C} \rightarrow \infty$,
we conclude that
\begin{eqnarray}
\label{88-7}
  \delta_0 \leq \mP\left(\inf_{\phi\in\cS_{\alpha,N}}\sum_{j\in \cZ_0} \int_{\sigma/2}^\sigma \langle
  K_{r,\sigma}^{(0)}\phi,
 e_j\rangle^2\dif S_r=0\right).
\end{eqnarray}
On the other hand, since   $K_{r,t}^{(0)}$  is  the solution  of equation (\ref{w-1})
with  $w_t$   issued from $w_t|_{t=0}=w_0^{(0)},$
(\ref{86-1}) implies that
\begin{eqnarray*}
  && \mP \Big(\inf_{\phi\in \cS_{\alpha,N}}\sum_{j\in \cZ_0} \int_{\sigma/2}^\sigma \langle
  K_{r,\sigma}^{(0)}\phi,
 e_j\rangle^2\dif S_r=0 \Big)
  =0.
\end{eqnarray*}
This is in  conflict with (\ref{88-7}) and the proof is  complete.

\section{Acknowledgement}
This work is partially supported by National Key R\&D program of China (No. 2022
YFA1006001), National Natural Science Foundation of China (Nos. 12071123, 12131019, 12371151,
11721101). Xuhui Peng is also   supported
by   the science and technology innovation Program of Hunan Province (No. 2022RC1189).
 Jianliang Zhai's research is also supported by  the Fundamental Research Funds for the Central Universities(Nos. WK3470000031, WK0010000081).

\end{document}